\documentclass[preprint,12pt,3p,fleqn]{elsarticle}
\usepackage[utf8]{inputenc}
\usepackage{geometry}
\usepackage{amssymb,float}
\usepackage{amsmath,placeins,amsthm,}
\usepackage{amsfonts,bm}
\journal{Computers \& Fluids}
\usepackage{makeidx}
\usepackage{graphicx,comment}
\usepackage{psfrag}
\usepackage{epsfig,epstopdf,calc,color}
\usepackage{caption,subcaption}
\usepackage{hyperref}
\usepackage{url} 
\hypersetup{colorlinks, 
citecolor= blue, 
linkcolor= blue,
urlcolor= blue}

\theoremstyle{theorem}
\newtheorem{theorem}{Theorem}

\theoremstyle{theorem}
\newtheorem{remark}{Remark}

\theoremstyle{claim}

\theoremstyle{proposition}
\newtheorem{proposition}[theorem]{Proposition}

\theoremstyle{theorem}
\newtheorem{definition}{Definition}

\title{A provably stable and high-order accurate finite difference approximation for the incompressible boundary layer equations}

\begin{document}
\begin{frontmatter}

\author[label1]{Mojalefa P.~Nchupang\corref{cor1}}
\ead{nchmoj002@myuct.ac.za}

\author[label1]{Arnaud G.~Malan}
\ead{arnaud.malan@uct.ac.za}

\author[label2]{Fredrik  Laur\'{e}n}
\ead{fredrik.lauren@liu.se}

\author[label2,label3]{Jan Nordstr\"om}
\ead{jan.nordstrom@liu.se}

\address[label1]{Industrial CFD Research Group,~Department of Mechanical Engineering,~University of Cape Town,~Cape Town 7700,~South Africa}
\address[label2]{Department of Mathematics,~Applied Mathematics,~Link\"oping University,581 83 Link\"oping,~Sweden}
\address[label3]{Department of Mathematics and Applied Mathematics,~University of Johannesburg,~Auckland Park 2006,~South Africa}
\cortext[cor1]{Corresponding author}

\begin{abstract}
In this article we develop a high order accurate method to solve the incompressible boundary layer equations in a provably stable manner.~We first derive continuous energy estimates,~and then proceed to the discrete setting.~We formulate the discrete approximation using high-order finite difference methods on summation-by-parts form and implement the boundary conditions weakly using the simultaneous approximation term method.~By applying the discrete energy method and imitating the continuous analysis,~the discrete estimate that resembles the continuous counterpart is obtained proving stability.~We also show that these newly derived boundary conditions removes the singularities associated with the null-space of the nonlinear discrete spatial operator.~Numerical experiments that verifies the high-order accuracy of the scheme and coincides with the theoretical results are presented.~The numerical results are compared with the well-known Blasius similarity solution as well as that resulting from the solution of the incompressible Navier Stokes equations.
\end{abstract}
\end{frontmatter}
\section{Introduction}\label{intro}
Conservation laws describing fluid dynamics mathematically take the form of space-time partial differential equations (PDEs).~One such example is the unsteady incompressible Navier-Stokes (INS) equations.~In the recent years,~numerical simulations of incompressible flows have gained attraction due to their numerous applications.~These include biomedical engineering \cite{Hume2022,YullPark2007,Cerrolaza1996},~aircraft design \cite{Alam2015,Kurzin2009,Haddadpour2008},~and atmospheric-ocean modeling \cite{Marshall1997b,Teixeira2013}.~Traditional second order numerical schemes have been widely used to analyze and predict flow parameters such as velocities and pressure \cite{Malan2022,Mowat2014}.~These second order accurate approaches however numerically damp flow vortexes \cite{Changfoot2019} while requiring excessive element numbers in the boundary layers.~Further,~mainstream incompressible flow solution schemes augment the incompressible mass conservation equation $\nabla \cdot \textbf{u}$ to avoid the resulting singular coefficient matrix.~The two main augmentation approaches are the so-called pressure-based (projection scheme) \cite{Patankar2018} and density-based (artificial compressibility) methods \cite{Chorin1967}.~These approaches introduce the need for more boundary conditions which place additional constraints on pressure gradients at boundaries \cite{Malan2002,Malan2013}.~Finally,~the ubiquitous practice of upwinding convective terms when solving incompressible flows \cite{Merrick2018} adds both complexity and nonphysical dissipation to the flow solution.
\par The key contributions of this article address these concerns.~For this purpose we employ the celebrated incompressible boundary layer equations as a model problem and endeavour to prove the existence of a stable and high order accurate solution without any need for additional augmented pressure/density based equations and without the use of upwinding.~The discretization method is finite difference on summation-by-parts (SBP) form \cite{KREISS1974,Gustafsson2008,Gustafsson2013}.~Key however to numerical stability is the correct boundary condition imposition, for which we employ a penalty-like method called simultaneous approximation term (SAT) \cite{Carpenter1994}.~The augmentation of the SBP operators with the SAT technique allows the proofs of stability to be straightforwardly attainable.~The stability of the numerical approximations ties back to the well-posedness of the continuous mathematical model which fundamentally depends on the choice of boundary conditions \cite{Gustafsson2008,Gustafsson2013,Strikwerda1977}.~To guarantee a bounded and stable numerical solution of a linear problem,~a minimal number \cite{Nordstrom2020b,Kreiss1970a} and appropriate form of boundary conditions must be known.~Well-posedness of nonlinear PDEs is unclear and incomplete,~however,~the linear theory can to some extent be extended to nonlinear problems using linearization principles.~We will follow the detailed guideline in \cite{Nordstrom2017a} and its application to the INS equations \cite{Lauren2022,Nordstrom2018},~the shallow water equations \cite{Nordstrom2022,Nordstrom2015a},~and Euler equations \cite{Nordstrom2022a,Nordstrom2022b} to derive a provably stable and high-order accurate approximation scheme for the boundary layer equations.~The energy method \cite{Gustafsson2008,Nordstrom2017a} which is based on the principle of integration-by-parts is used as a basic tool to derive the desired boundary conditions that yields an estimate.~Furthermore,~the stable imposition of these boundary conditions using SAT eliminates the saddle point problem typically associated with the spatial operator of the incompressible flow equations leading to unique solutions \cite{Benzi2005,Nordstrom2020c,Lauren2021}.
\par To set up the continuous problem,~we consider a laminar incompressible flow aligned with a thin plate of length $l$.~We model the problem  using the laminar incompressible boundary layer (IBL) equations which are derived from the INS equations at high Reynolds number (Re) using dimensional analysis (see \cite{Frank2006b} for details).~Of note is that the continuity equation contains no explicit relationship between the pressure and the velocity gradients.~As noted above,~popular 2nd order methods deal with this by creating such a relationship artificially i.e.~by using artificial compressibility \cite{Chorin1967,Malan2002,Rahman2001},~the pressure projection scheme \cite{Vreman2014a,Matsui2021} or a combination of the two \cite{Malan2013,Oxtoby2012a}.~Staggering grids is another popular method used to enforce divergence and damp spurious oscillations in the solution domain \cite{OReilly2017a,Kress2003,Gustafsson2000}.~In this work we will only use the fundamental equations,~and demonstrate that the resulting scheme is both stable and accurate,~by using the SBP-SAT technique.~Importantly,~this is without the use of so-called upwinding for discretization.~Note that though this work employs high-order finite difference approximations,~the analysis also holds for any numerical approximations that can be written on SBP form.~ Examples include discontinuous Galerkin method \cite{Manzanero2019,Chan2018},~spectral element method \cite{Yamaleev2017},~finite element method \cite{Abgrall2020},~and finite volume method \cite{Ham2006,Nordstrom2003}. 
\par The rest of the paper is organized as follows;~we present the continuous analysis in Section $\ref{sec2}$ and derive the energy-stable boundary conditions.~In Section $\ref{sec3}$,~we impose these boundary conditions and deduce the penalty coefficients such that we get a bounded energy estimate.~We formulate the SBP-SAT semi-discrete approximation in Section $\ref{sec4}$ and mimic the continuous analysis to derive the discrete energy estimate that resemble the continuous one.~Moreover,~we study the null-space of the discrete spatial operator and investigate the effect of SAT boundary conditions on the positive definiteness of the resulting coefficient matrix.~Temporal discretization is considered in Section $\ref{sec5}$.~The computational results that verifies the high-order accuracy of the approximation scheme are presented in Section $\ref{sec6}$ and the Blasius boundary layer is considered as a validation model.~We also make comparison with the full INS equations and draw conclusion in Section $\ref{sec7}$.
\begin{center}
\begin{figure}[H]
\centering
\psfrag{t0}{$y$} 
\psfrag{t2}{$U_\infty$}
\psfrag{t3}{$0$}
\psfrag{t4}{$\delta(x)$}
\psfrag{t5}{$l$}
\psfrag{t6}{$x$}
\includegraphics[scale=0.7]{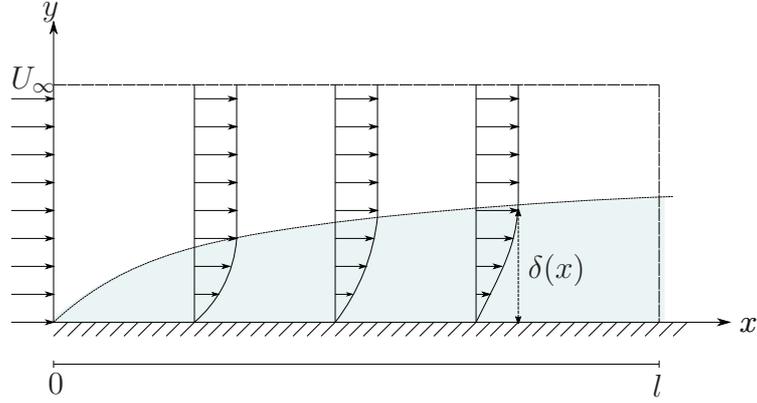}
\caption{Viscous fluid flowing over a thin plate of length $l$ leading to the formation of boundary layer.}
\label{fig1}
\end{figure} 
\end{center}
\section{The continuous problem}\label{sec2}
We consider a viscous fluid flowing past a thin plate of length $l$ with a uniform speed $U_\infty >0$.~The flow is laminar and incompressible with $\text{Reynolds number} \gg 1$,~leading to the development of boundary layer of thickness $\delta$ as depicted in Figure $\ref{fig1}$.~Let $\Omega \subset \mathbb{R}^2$ be the computational domain with Cartesian coordinates $\textbf{x}= (x,y)$ and boundary $\Gamma$.~We position the plate along the $x$-axis such that the leading edge is at the origin.~Further,~we denote the components of the velocity in the $x$-,~$y$-directions with $u$ and $v$,~respectively and the static pressure with $p$.~We start our analysis from the two-dimensional (2D) INS equations for Newtonian fluids under isothermal conditions
\begin{align}
u_t + uu_x + vu_y  &=-\frac{1}{\rho}p_x +  \frac{\mu}{\rho} ( u_{xx} + u_{yy}),  \nonumber \\
v_t + uv_x + vv_y  &= -\frac{1}{\rho}p_y + \frac{\mu}{\rho}( v_{xx} + v_{yy}),  \label{eq1} \\
u_x + v_y &= 0. \nonumber
\end{align}
Here,~$\rho> 0$ is the constant density and will be set to one in the subsequent analysis while $\mu>0$ denotes the constant dynamical viscosity.~The subscripts in $\eqref{eq1}$ denote the temporal and spatial partial derivatives.~At large Reynolds number (defined as $\text{Re} = \rho \frac{ U_\infty l}{\mu}$),~it can be shown using dimensionless scaling \cite{Frank2006b} that $\eqref{eq1}$ reduces to the IBL equations
\begin{align}
u_t + uu_x + vu_y  &=- p_x +\mu u_{yy}, \nonumber \\
0 &= -p_y ,  \label{eq2} \\
u_x + v_y &= 0. \nonumber
\end{align}
\par We begin the continuous analysis by writing $\eqref{eq2}$ as an initial-boundary value problem.~In matrix-vector form,~the system $\eqref{eq2}$ with boundary and initial conditions included can be written as
\begin{align}
\mathcal{I} \text{U}_t + \mathcal{D}(\text{U})\text{U} &= 0,  \quad \qquad \textbf{x} \in \Omega , \quad t>0, \nonumber \\
\mathcal{B}\text{U} &= \textbf{g}(\textbf{x},t), \quad \textbf{x} \in \Gamma, \quad t>0, \label{eq3}    \\
\mathcal{I} \text{U} &= \textbf{f}(\textbf{x},t), \quad \textbf{x} \in \Omega , \quad t=0 , \nonumber
\end{align}
where $\text{U} = [u,v,p]^T$ contain the dependent variables.~The continuous vector functions $\textbf{g}$,~$\textbf{f}$ are known and specifies boundary and initial data to the problem (we assume that they are compatible such that the solution is smooth).~Furthermore,~the exact form of the boundary operator $\mathcal{B}$ will be determined later.~Lastly,~in $\eqref{eq3}$,~$\mathcal{D}$ is the spatial nonlinear operator and is expressed as
\begin{align*}
\mathcal{D}(\text{U}) &= A(\text{U})\frac{\partial}{\partial x} + B(\text{U})\frac{\partial}{\partial y} - \mu \mathcal{I}\frac{\partial^2}{\partial y^2}    ,
\intertext{where}
\mathcal{I} &=\begin{bmatrix}
1 & 0 & 0 \\ 0&0&0 \\0&0&0
\end{bmatrix} ,\quad  A = \begin{bmatrix}
u & 0 & 1 \\0&0&0\\1&0&0
\end{bmatrix}, \quad B= \begin{bmatrix}
v & 0 & 0 \\0&0&1\\0&1&0
\end{bmatrix}.  \nonumber
\end{align*}
\noindent Furthermore,~we split the convective terms in terms of the conservative and non-conservative components using the flux splitting technique \cite{Nordstrom2018,Nordstrom2006a}
\begin{align}
AU_x = \frac{1}{2} \left[AU_x + (AU)_x - A_xU\right] , \quad BU_y = \frac{1}{2} \left[BU_y+ (BU)_y - B_yU\right]. \label{eq4}
\end{align}
\begin{remark}The flux splitting $\eqref{eq4}$ is crucial for the upcoming discrete analysis. \end{remark}
\noindent Noting that $A_x +B_y = \mathcal{I}(u_x + v_y) = 0$,~the skew-symmetric form of the governing system in $\eqref{eq3}$ becomes
\begin{align}
\mathcal{I} \text{U}_t + \dfrac{1}{2} \big[A\text{U}_x + (A\text{U})_x   + B\text{U}_y + (B\text{U})_y \big] -  \mu \mathcal{I}\text{U}_{yy} &=  0 ,\label{eq5}
\end{align}
i.e.~$\mathcal{D}(\text{U})\text{U} =  \dfrac{1}{2} \big[A\text{U}_x + (A\text{U})_x   + B\text{U}_y + (B\text{U})_y \big] -  \mu \mathcal{I}\text{U}_{yy}.$
\begin{center}
\begin{figure}[H]
\centering
\psfrag{d0}{$\Omega$} 
\psfrag{d1}{$y$}
\psfrag{d2}{$\Gamma_w$}
\psfrag{d3}{$\Gamma_s$}
\psfrag{d4}{$x$}
\psfrag{d5}{$\Gamma_e$}
\psfrag{d6}{$\Gamma_n$}
\psfrag{c1}{$\bm{n}^w = (-1,0)$}
\psfrag{c2}{$\bm{n}^s = (0,-1)$}
\psfrag{c3}{$\bm{n}^e = (1,0)$}
\psfrag{c4}{$\bm{n}^n = (0,1)$}
\includegraphics[scale=0.6]{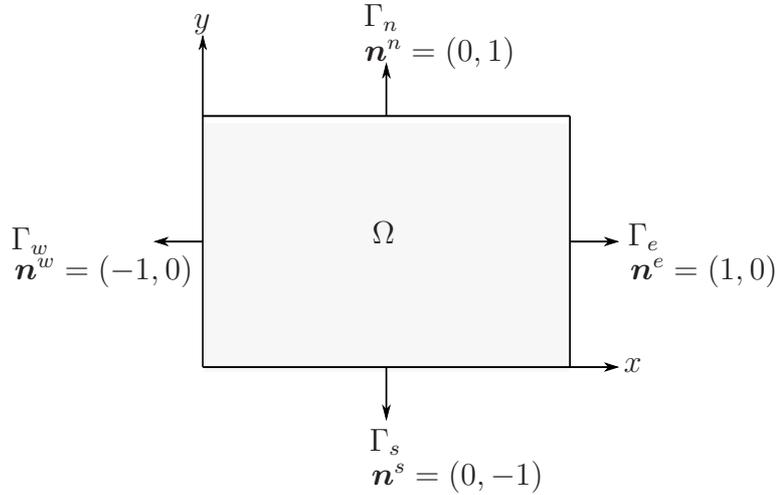}
\caption{Two-dimensional computational domain $\Omega$ showing outwards pointing normal vectors along the boundary $\Gamma$.}
\label{fig2}
\end{figure} 
\end{center}
\subsection{Boundedness}
Next,~we employ the energy method to derive the appropriate form of boundary conditions in $\eqref{eq3}$ that leads to an energy estimate.~The energy method (which involves multiplying $\eqref{eq5}$ with the $2\text{U}^T$ and integrating over the computational domain $\Omega$) applied to $\eqref{eq5}$ yields
\begin{align}
2\int \limits_\Omega U^T \mathcal{I} U_t dV + \int\limits_\Omega \text{U}^T \big[A\text{U}_x + (A\text{U})_x+ B\text{U}_y + (B\text{U})_y \big] dV&= 2\mu \int\limits_\Omega \text{U}^T\mathcal{I}\text{U}_{yy} dV, \label{eq6}
\end{align}
where $dV = dxdy$ is the volume element.~Let $||\text{U}||^2_\mathcal{I} = \int \limits_\Omega \text{U}^T \mathcal{I}\text{U}dV$ denote the $L_2$ semi-norm.~Then,~by using integration by parts (IBP) and the Divergence theorem $\int \limits_\Omega 
 U_{x_i} dV = \oint \limits_\Gamma U n_{x_i} ds$ to simplify $\eqref{eq6}$,~we obtain
\begin{align}
\frac{d}{dt}||\text{U}||_\mathcal{I}^2 &+ \int\limits_\Omega  (U^TAU_x - U_x^TAU )dV  + \int\limits_\Omega  (U^TBU_y - U_y^TBU )dV + 2\mu \int\limits_\Omega  (U_y^T\mathcal{I} U_y  )dV \label{eq7}\\&= -\oint\limits_{\Gamma} \big[ \text{U}^T(An_x + Bn_y)\text{U} - 2\mu  \text{U}^T \mathcal{I}\text{U}_yn_y \big] ds.  \nonumber
 \end{align}
Since $A$,~$B$ in $\eqref{eq7}$ are symmetric,~the non-conservative convective terms on the left-hand side (LHS) of $\eqref{eq7}$ vanishes and the energy rate becomes
 \begin{align}
\frac{d}{dt}||\text{U}||_\mathcal{I}^2 + 2\mu ||\text{U}_y||^2_\mathcal{I} &= -\oint\limits_{\Gamma} \big[ \text{U}^T(An_x + Bn_y)\text{U} - 2\mu  \text{U}^T \mathcal{I}\text{U}_yn_y \big] ds \label{eq8} \\& = -\oint\limits_{\Gamma} (u_{\bm{n}} u^2 + 2u_{\bm{n}}p - 2\mu uu_y n_y)ds. \nonumber
\end{align}
In $\eqref{eq7}$,~$\oint\limits_{\Gamma}(\cdot)ds$ is the boundary line integral with the infinitesimal line element $ds = \sqrt{dx^2 + dy^2}$ along boundary $\Gamma$.~Furthermore,~$u_{\bm{n}} = un_x + vn_y$ is the boundary normal velocity and $\bm{n} = (n_x,n_y)$ is the normal outward pointing unit vector as depicted in Figure \ref{fig2}.~Let $BT$ denote the boundary terms in $\eqref{eq8}$,~then
\begin{align}
\text{BT} =  -\oint\limits_{\Gamma} (u_{\bm{n}} u^2 + 2u_{\bm{n}}p - 2\mu uu_y n_y)ds.  \label{eq9}
\end{align}
A bounded energy estimate is guaranteed if BT is non-positive.~This can be achieved by imposing appropriate boundary conditions.~That is,~we need to establish the correct minimal number and form of boundary conditions \cite{Nordstrom2020b,Nordstrom2017a}.~To do that,~we rewrite the boundary terms $\eqref{eq9}$ in a matrix-vector form 
\begin{align}
\text{BT} &= -\oint \limits_{\Gamma} \textbf{q}^T \text{M} \textbf{q} ds = -\oint\limits_{\Gamma} \begin{bmatrix}
u\\v\\p\\ \mu u_y
\end{bmatrix}^T \begin{bmatrix}
u_{\bm{n}} &0&n_x &-n_y \\ 0 & 0 & n_y & 0 \\ n_x & n_y & 0 & 0 \\-n_y &0 & 0 & 0
\end{bmatrix}\begin{bmatrix}
u\\v\\p\\ \mu u_y
\end{bmatrix}ds.  \label{eq10}
\end{align}
\begin{proposition}The number of boundary conditions required to bound $\eqref{eq8}$ coincide with the number of negative eigenvalues of M. \end{proposition}
\begin{proof}
See \cite{Nordstrom2020b}.
\end{proof}
\noindent The eigenvalues of M are obtained by solving for the roots of the characteristic polynomial
\begin{align*}
\det(\text{M}-\lambda I_4) = ( \lambda^2 -n_y^2)(\lambda^2 - u_{\bm{n}}\lambda - 1) = 0,
\end{align*}
where $I_4$ is a unit matrix of size $4$.~We will consider each boundary separately,~starting with the north and south boundaries.~Noting that $(n_x,n_y) = (0,\pm 1)$,~the eigenvalues $\lambda_i$ and the associated eigenvectors $x_i$ are
\begin{align}
\lambda_{1} &= \frac{u_{\bm{n}}}{2} - \sqrt{\left(\frac{u_{\bm{n}}}{2}\right)^2 + 1}, \quad \lambda_{2} = - 1, \quad \lambda_{3} = + 1, \quad \lambda_{4} = \frac{u_{\bm{n}}}{2} + \sqrt{\left(\frac{u_{\bm{n}}}{2}\right)^2 + 1} ,\label{eq11} \\
x_1 &= \begin{bmatrix}
\lambda_1 \\ 0 \\ 0 \\ -n_y
\end{bmatrix}, \quad    x_2 = \begin{bmatrix}
0 \\ 1 \\ -n_y \\ 0
\end{bmatrix},  \quad x_3= \begin{bmatrix}
0 \\ 1 \\ n_y \\ 0
\end{bmatrix}, \quad x_4 = \begin{bmatrix}
\lambda_4 \\ 0\\ 0 \\ -n_y
\end{bmatrix}.  \nonumber
\end{align}
For the east and west boundaries,~$(n_x,n_y) = (\pm1,0)$,~we have
\begin{align}
\lambda_{1} &= \frac{u_{\bm{n}}}{2} - \sqrt{\left(\frac{u_{\bm{n}}}{2}\right)^2 + 1}, \quad \lambda_{2} = 0, \quad \lambda_{3} = 0, \quad \lambda_{4} = \frac{u_{\bm{n}}}{2} + \sqrt{\left(\frac{u_{\bm{n}}}{2}\right)^2 + 1} ,\label{eq12} \\
x_1 &= \begin{bmatrix}
\lambda_1 \\ 0 \\ n_x \\ 0
\end{bmatrix}, \quad x_2 = \begin{bmatrix}
0 \\ n_x \\ 0\\ 0
\end{bmatrix},  \quad x_3= \begin{bmatrix}
0 \\ 0\\ 0 \\ n_x
\end{bmatrix}, \quad x_4 = \begin{bmatrix}
\lambda_4 \\ 0\\ n_x \\ 0
\end{bmatrix}.  \nonumber
\end{align}
\begin{remark} We assume outflow $(u_\textbf{n} >0) $ at the north and east boundaries,~and inflow $(u_\textbf{n} \leq 0)$ at the west and south boundaries. \end{remark}
\begin{remark}
In $\eqref{eq11}$ and $\eqref{eq12}$,~$\lambda_1< 0$ and $\lambda_4>0$ for all $u_\textbf{n}$. 
\end{remark}
\noindent Therefore,~there are two negative eigenvalues $\lambda_1$,~$\lambda_2$ and two positive eigenvalues $\lambda_3$,$~\lambda_4$ in $\eqref{eq11}$ at the north and south boundaries.~This infers that precisely two boundary conditions must be prescribed at both the north and south boundaries (this is due to the presence of $\mu u_{yy}$ in the equations).~At the east and west boundaries,~there is only one negative eigenvalue $\lambda_1$  in $\eqref{eq12}$ inferring that only one boundary condition must be imposed at each boundary.
\par To determine the form of the boundary conditions that will lead to a finite energy estimate,~we return to $\eqref{eq10}$ and consider the eigenvalue decomposition of M
\begin{align}
\text{M} = \bar{\text{X}} \Lambda_\text{M}\bar{\text{X}}^T ,  \label{eq13}
\end{align}
where $\Lambda_\text{M} = \text{diag}(\lambda_1,\lambda_2,\lambda_3,\lambda_4)$ and $\bar{\text{X}} = \text{XN}$ are the scaled eigenvector matrices X whose columns are eigenvectors in $\eqref{eq11}$ and $\eqref{eq12}$ for respective boundaries.~The columns of X are scaled with the normalizing matrix N.~We further rearrange M as $\text{M} = \text{X}\Lambda \text{X}^T$ where $\Lambda = \text{N}^T \Lambda_\text{M} \text{N}$ is the scaled version of $\Lambda_\text{M}$.~For the north and south boundaries,~X and N are
\begin{align*}
\text{X}&= \begin{bmatrix}
\lambda_1 & 0 & 0 & \lambda_4\\ 0 & 1 & 1 & 0 \\ 0& -n_y & n_y & 0\\ -n_y & 0 &  0 & -n_y
\end{bmatrix}, \quad \text{N} = \text{diag}\left( \sqrt{\lambda^2_1 + 1},\sqrt{2},\sqrt{2},\sqrt{ \lambda^2_4 + 1 } \right)^{-1},
\end{align*}
and similarly for the east and west boundaries,~they are
\begin{align*}
\text{X}&= \begin{bmatrix}
\lambda_1 & 0 & 0 & \lambda_4\\ 0 & n_x & 0 & 0 \\ n_x& 0 & 0 & n_x\\ 0 & 0 &  n_x & 0
\end{bmatrix}, \quad \text{N} = \text{diag}\left( \sqrt{\lambda^2_1 + 1},1,1,\sqrt{ \lambda^2_4 + 1 } \right)^{-1}.
\end{align*}
Therefore,~by substituting  $\eqref{eq13}$ into $\eqref{eq10}$,~BT can be rewritten as
\begin{align}
\text{BT}&= -\oint \limits_{\Gamma} W^T \Lambda W ds ,   \label{eq14}
\end{align}
where $W = \text{X}^T \textbf{q}$.~For the north and south boundaries,~$W$ is
\begin{align*}
W = \begin{bmatrix}
\lambda_1 u - \mu u_y n_y \\ v - pn_y \\ v+pn_u \\ \lambda_4 u - \mu u_y n_y 
\end{bmatrix},
\end{align*}
and for the east and west boundaries,~$W$ is 
\begin{align*}
W = \begin{bmatrix}
\lambda_1 u + pn_x \\ v n_x  \\ \mu u_y n_x \\ \lambda_4 u + pn_x
\end{bmatrix}.
\end{align*}
Following \cite{Nordstrom2018},~we partition $\Lambda$ in terms of the positive,~zero,~and negative components i.e~$\Lambda = \text{diag}(\Lambda_+,\Lambda_0,\Lambda_-)$.~Similarly,~we write the corresponding variables as $W = [W_+,W_0,W_-]^T$ where $W_+,W_-$ are called the incoming and outgoing characteristics,~respectively.~The variable $W_0$ which is associated with $\Lambda_0$ is not interesting since $W_0^T \Lambda_0 W_0 = 0$ and it will be omitted in the subsequent derivations.~Noting that $\lambda_1,\lambda_2 < 0$ and $\lambda_3,\lambda_4>0$ at the north and south boundaries,~the positive and negative components of the matrix decomposition in $\eqref{eq14}$ are
\begin{align}
 W_+ &= \begin{bmatrix}v + pn_y  \\
\lambda_4 u  - \mu u_yn_y 
\end{bmatrix},\quad W_- = \begin{bmatrix}
\lambda_1 u - \mu u_yn_y \\  v - pn_y 
\end{bmatrix}, \label{eq15}\\ \Lambda_+ &= \begin{bmatrix}
\dfrac{\lambda_3}{2} & 0 \\ 0 & \dfrac{\lambda_4}{\lambda^2_4+1}
\end{bmatrix}, \quad \Lambda_- = \begin{bmatrix}
\dfrac{\lambda_1}{\lambda^2_1+1} & 0 \\ 0 & \dfrac{\lambda_2}{2}
\end{bmatrix}. \nonumber
\end{align}
Similarly,~noting that $\lambda_1< 0$ and $\lambda_4>0$ in $\eqref{eq12}$ for the east and west boundaries,~we have
\begin{align}
 W_+ &= \begin{bmatrix}
\lambda_4 u + pn_x 
\end{bmatrix}, \quad W_- = \begin{bmatrix}
\lambda_1 u + pn_x \end{bmatrix}, \quad \Lambda_+ = \begin{bmatrix}
\dfrac{\lambda_4}{\lambda^2_4+1}
\end{bmatrix}, \quad \Lambda_- = \begin{bmatrix}
\dfrac{\lambda_1}{\lambda^2_1+1}
\end{bmatrix} .\label{eq16}
\end{align}
Equation $\eqref{eq14}$ with the partition above then becomes
\begin{align}
\text{BT} = -\oint\limits_\Gamma \begin{bmatrix}
W_+ \\ W_-
\end{bmatrix}^T \begin{bmatrix}
\Lambda_+ & 0 \\ 0 & \Lambda_-
\end{bmatrix}\begin{bmatrix}
W_+ \\ W_-
\end{bmatrix} ds. \label{eq17}
\end{align}
We overcome the energy growth due the negative eigenvalues by specifying the boundary condition \cite{Nordstrom2017a} 
\begin{align}
W_- &= \mathcal{S}W_+ + \textbf{g},  \label{eq18}
\end{align}
i.e.~specifying the incoming characteristics in terms of the outgoing ones and data.~Here,~$\mathcal{S}$ is a matrix with the number of rows equal to the number of negative eigenvalues and the number of columns equal to the number of positive eigenvalues.
\begin{remark} The general form of the relation $\eqref{eq18}$ is $\mathcal{R}\big(W_- - \mathcal{S}W_+ \big) = \textbf{g}$ \cite{Nordstrom2023a} where $\mathcal{R}$ is an invertible matrix that combines the variables in $W_-$ and $\mathcal{RS}$ combines the variables in $W_+$ to implement the desired boundary conditions.~In this work,~we however choose $\mathcal{R}$ to be an identity matrix. \end{remark}
\noindent Substituting $\eqref{eq18}$ into $\eqref{eq17}$ leads to 
\begin{align}
\text{BT} = -\oint \limits_\Gamma \begin{bmatrix}
W_+ \\ \textbf{g}
\end{bmatrix}^T \begin{bmatrix}
\Lambda_+ + \mathcal{S}^T\Lambda_-\mathcal{S} & \mathcal{S}^T\Lambda_- \\ \Lambda_-\mathcal{S} & \Lambda_-
\end{bmatrix}\begin{bmatrix}
W_+ \\ \textbf{g}
\end{bmatrix} ds. \label{eq19}
\end{align}
By assuming a homogeneous form of the boundary condition in $\eqref{eq3}$ i.e.~$\textbf{g} = 0$ such that $W_- = \mathcal{S}W_+$,~BT further simplifies to
\begin{align}
\text{BT} = -\oint \limits_\Gamma 
W_+^T \left(
\Lambda_+ + \mathcal{S}^T\Lambda_-\mathcal{S} \right)W_+ds, \label{eq20}
\end{align}
which is non-positive if we can choose $\mathcal{S}$ such that
\begin{align}
\Lambda_+ + \mathcal{S}^T\Lambda_-\mathcal{S} \geq 0 . \label{eq21}
\end{align}
The non-homogeneous case is considered in \cite{Nordstrom2018,Nordstrom2017a} and that analysis will be omitted herein.~Based on $\eqref{eq19}$,$\eqref{eq20}$,~we observe that $\eqref{eq18}$ defines the general form that the boundary conditions in $\eqref{eq3}$ must have in order for $\text{BT}$ in $\eqref{eq9}$ to be non-positive
\begin{align}
\mathcal{B}\text{U} = W_- - \mathcal{S}W_+ = \textbf{g},  \label{eq22}
\end{align}
where $\mathcal{S}$ satisfies $\eqref{eq21}$.
\subsection{The energy stable boundary conditions}
Although the general form of boundary conditions $\eqref{eq22}$ yields a bounded energy norm for the continuous equation,~we need specific boundary conditions that are aligned with the physics of the original problem i.e.~we prescribe them according to the available boundary data.~By returning to $\eqref{eq9}$,~we propose a new set of energy stable boundary conditions satisfying $\eqref{eq21}$ and further show that they can be written in the general form $\eqref{eq22}$.
\par Starting at the south boundary which is aligned with the solid surface as depicted in Figure~$\ref{fig1}$,~we eliminate the contribution of the south boundary in BT by prescribing a no-slip velocity condition i.e.~$u = v = 0$.~The north boundary is considered next,~at which we should impose two boundary conditions.~It is important to note that at this boundary,~$u_{\bm{n}}>0$.~Next,~we turn to the vertical boundaries.~Since we assumed inflow $(u_{\bm{n}}  < 0)$ at the west boundary and outflow $(u_{\bm{n}}  > 0)$ at the east boundary,~we prescribe the following boundary condition; velocity at the west boundary and pressure at the east boundary.~In summary,~the proposed boundary conditions are
\begin{align}
u &=0 , \quad  v = 0, \qquad  \text{south boundary},  \nonumber \\
\frac{\alpha}{2} vu  -\mu u_y  &= 0, \quad p = p_\infty, \qquad   \text{north boundary}, \label{eq23}\\ 
u &= U_\infty, \quad \qquad  \qquad \text{west boundary}, \nonumber\\
p &= p_\infty, \qquad \qquad  \qquad \text{east boundary}, \nonumber
\end{align}
where $\alpha \in [0,1]$ is a constant that gives us the flexibility to impose either the Robin boundary condition $(\alpha = 1 )$ or the Neumann boundary condition $(\alpha = 1 )$ depending on the available data.~The above can be written in the form $\mathcal{B}\text{U} = \textbf{g}$ as
\begin{align}
\mathcal{B}_e\text{U}&= \begin{bmatrix}
0 & 0 & 1 
\end{bmatrix}  \begin{bmatrix}
u \\ v \\ p
\end{bmatrix} =\begin{bmatrix}
p_\infty
\end{bmatrix}, \qquad \mathcal{B}_n\text{U}= \begin{bmatrix}
\dfrac{\alpha}{2}  v - \mu \partial_y & 0 & 0 \\ 0 & 0 & 1 
\end{bmatrix} \begin{bmatrix}
u \\ v \\ p
\end{bmatrix}= \begin{bmatrix}
0 \\ p_\infty \end{bmatrix},  \label{eq24}\\
\mathcal{B}_w\text{U} &= \begin{bmatrix}
1 & 0 & 0 
\end{bmatrix} \begin{bmatrix}
u \\ v \\ p
\end{bmatrix} = \begin{bmatrix}
U_\infty
\end{bmatrix}  , \qquad \mathcal{B}_s\text{U} = \begin{bmatrix}
1 & 0 & 0 \\ 0 & 1 & 0
\end{bmatrix}  \begin{bmatrix}
u \\ v \\ p
\end{bmatrix} = \begin{bmatrix}
0 \\ 0
\end{bmatrix} .   \nonumber
\end{align}
Here,~the subscripts $e,n,w,s$ denotes the east,~north,~west and south boundaries respectively as shown in Figure~$\ref{fig2}$ and $p_\infty$ is the freestream pressure.
\par By strongly imposing the homogeneous form of $\eqref{eq24}$ in $\eqref{eq8}$,~most of the boundary terms vanishes and only the contribution from the north and east boundaries remains which carries appropriate signs since $u_{\bm{n}}>0$ at the outflow boundaries.~The energy rate $\eqref{eq8}$ becomes
\begin{align}
\frac{d}{dt} ||\text{U}||^2_\mathcal{I} + 2\mu || \text{U}_y||^2_\mathcal{I} = -\int \limits_{\Gamma_n} (1-\alpha )vu^2 dx - \int \limits_{\Gamma_e} u^3 dy \leq 0.  \label{eq25}
\end{align}
Notice that $\alpha  = 0$,~which prescribes the Neumann boundary condition in $\eqref{eq23}$,~leads to a more dissipative energy rate.~Finally,~temporal integration over a finite time domain $[0,T]$ and imposing the initial condition leads to the energy estimate
\begin{align}
||\text{U}||^2_\mathcal{I} + 2\mu \int \limits_0^T || \text{U}_y||^2_\mathcal{I} dt \leq ||\textbf{f}||_{\mathcal{I}}^2 .  \label{eq26}
\end{align}
\begin{remark} The bound is imposed only on the horizontal component of the velocity since it is the only flow variable in $\eqref{eq2}$ with the temporal derivative.~This is different compared to the fully INS equations where also the vertical velocity is bounded \cite{Nordstrom2018}.~In both models,~there is no bound on the pressure and it was shown in the cited work that we don't need one for boundedness.
\end{remark}
The estimate $\eqref{eq26}$ shows that the new boundary conditions $\eqref{eq23}$ are energy stable.~Next,~we compute $\mathcal{S}$ for each boundary satisfying $\eqref{eq21}$ and show that the boundary conditions $\eqref{eq23}$ can be written in the general form $\eqref{eq18}$.~We begin with the north and south boundaries.~By proposing $\mathcal{S}$ with the form
\begin{align*}
\mathcal{S} = \begin{bmatrix}
0 & s_1 \\ s_2 & 0
\end{bmatrix},
\end{align*} 
and substituting it together with the variables $W_+$,~$W_-$ in $\eqref{eq15}$ into $\eqref{eq18}$,~we get
\begin{align}
W_- - \mathcal{S}W_+ &= \begin{bmatrix}
\lambda_1u - \mu u_yn_y \\ v - pn_y
\end{bmatrix} - \begin{bmatrix}
0 & s_1 \\ s_2 & 0
\end{bmatrix} \begin{bmatrix}
v+ pn_y \\ \lambda_4 u - \mu u_y n_y
\end{bmatrix}  \label{eq27} \\& = \begin{bmatrix} (\lambda_1 - s_1\lambda_4)u + (s_1 -1)\mu u_y n_y \\
(1-s_2)v - (1+s_2)pn_y \end{bmatrix}   = \begin{bmatrix}
g_1 \\ g_2
\end{bmatrix}.   \nonumber
\end{align} 
Here,~$g_1$,~$g_2$ denotes boundary data.~To write the no-slip condition at the south boundary in the form $\eqref{eq18}$,~we seek $\mathcal{S}$ that transforms $\eqref{eq27}$ into
\begin{align*}
W_- - \mathcal{S}W_+ = R\begin{bmatrix}
u\\v
\end{bmatrix}= \begin{bmatrix}
0\\0
\end{bmatrix},  
\end{align*}
where $R$ is a non-singular matrix.~Setting $s_1 = 1$ and $s_2 = -1$ yields the desired results and consequently,~satisfies $\eqref{eq21}$ since 
\begin{align}
\Lambda_+ + \mathcal{S}^T\Lambda_-\mathcal{S}  = \begin{bmatrix}
\dfrac{\lambda_3}{2} & 0 \\ 0 & \dfrac{\lambda_4}{\lambda^2_4+1}
\end{bmatrix} + \begin{bmatrix}
0 &1 \\ -1 & 0
\end{bmatrix}^T \begin{bmatrix}
\dfrac{\lambda_1}{\lambda^2_1+1} & 0 \\ 0 & \dfrac{\lambda_2}{2}
\end{bmatrix}\begin{bmatrix}
0 &1 \\ -1 & 0
\end{bmatrix} = \begin{bmatrix}
0 &0 \\ 0 & 0
\end{bmatrix} .   \label{eq28}
\end{align}
 and $R = \text{diag}(\lambda_1 - \lambda_4, 2)$.
\par To show that the Robin velocity $(\alpha = 1 )$ and the Dirichlet pressure conditions at the north boundary can be written in the form $\eqref{eq18}$,~we set $s_1 = 0$ and $s_2 = 1$ in $\eqref{eq27}$ which leads to
\begin{align*}
W_- - \mathcal{S}W_+ = \begin{bmatrix}
\lambda_1u - \mu u_yn_y \\ - 2pn_y
\end{bmatrix}  = \begin{bmatrix}
1& 0 \\ 0 & -2
\end{bmatrix}\begin{bmatrix}
\lambda_1(v)u - \mu u_y \\ p
\end{bmatrix} = \begin{bmatrix}
g_1 \\ g_2
\end{bmatrix}.
\end{align*}  
Moreover,~$\eqref{eq21}$ is satisfied by this choice since
\begin{align}
\Lambda_+ + \mathcal{S}^T\Lambda_-\mathcal{S}  = \begin{bmatrix}
\dfrac{\lambda_3}{2} & 0 \\ 0 & \dfrac{\lambda_4}{\lambda^2_4+1}
\end{bmatrix} + \begin{bmatrix}
0 &0 \\ 1 & 0
\end{bmatrix}^T \begin{bmatrix}
\dfrac{\lambda_1}{\lambda^2_1+1} & 0 \\ 0 & \dfrac{\lambda_2}{2}
\end{bmatrix}\begin{bmatrix}
0 &0 \\ 1 & 0
\end{bmatrix} = \begin{bmatrix}
0 &0 \\ 0 & \dfrac{\lambda_4}{\lambda^2_4+1}
\end{bmatrix} \geq 0.  \label{eq29}
\end{align}
\par Next,~we turn to the east and west boundaries.~Similar to the horizontal boundaries,~we want to show that the west and east boundary conditions can be written in the form $\eqref{eq18}$.~We begin by substituting $W_+$,~$W_-$ in $\eqref{eq16}$ into $\eqref{eq18}$ to obtain
\begin{align}
W_- - \mathcal{S}W_+ = \begin{bmatrix} \lambda_1 u + pn_x \end{bmatrix} - \mathcal{S} \begin{bmatrix} \lambda_4 u + pn_x \end{bmatrix}= (\lambda_1 - \mathcal{S} \lambda_4)u + (1 - \mathcal{S})pn_x = \textbf{g} \label{eq30}.
\end{align}
where $\textbf{g}$ denotes data as before.~Starting with the west boundary,~we want to determine $\mathcal{S}$ such that it removes the pressure term from $\eqref{eq30}$ and only the velocity remains.~The obvious choice $\mathcal{S} = 1$ leads to
\begin{align*}
W_- - \mathcal{S}W_+ = (\lambda_1 - \lambda_4) u = \textbf{g}.
\end{align*}
Consequently,~this choice satisfies $\eqref{eq21}$ since 
\begin{align}
\Lambda_+ + \mathcal{S}^T\Lambda_-\mathcal{S}   = \dfrac{\lambda_4}{\lambda_4^2 + 1} + \dfrac{\lambda_1}{\lambda_1^2 + 1} = 0,  \label{eq31}
\end{align}
where $\Lambda_+$,~$\Lambda_-$ are given in $\eqref{eq16}$.
\par Similarly,~to write the pressure condition at the east boundary in the form $\eqref{eq18}$,~we need appropriate $\mathcal{S}$ satisfying $\eqref{eq21}$ to remove the velocity contribution in $\eqref{eq30}$.~Setting $\mathcal{S} = \frac{\lambda_1}{\lambda_4}$ yields
\begin{align*}
W_- - \mathcal{S}W_+ = \left( 1 - \frac{\lambda_1}{\lambda_4 }\right) p =\textbf{g},
\end{align*}
which  also satisfies $\eqref{eq21}$
\begin{align}
\Lambda^+ + \mathcal{S}^T\Lambda^-\mathcal{S} = \frac{\lambda_4}{\lambda_4^2+1} + \left( \frac{\lambda_1}{\lambda_4} \right)^2\frac{\lambda_1}{\lambda_1^2 + 1} = \frac{u_\textbf{n}(u_\textbf{n}^2 + 1)}{(\lambda_4^2+1)(\lambda_1^2+1)\lambda_4}\geq 0 ,\label{eq32}
\end{align}
since $u_\textbf{n}>0$ at the east boundary. 
\par To recap what we did,~we rotated the boundary matrix M in $\eqref{eq10}$ into the diagonal form using the matrix eigenvalue decomposition $\eqref{eq13}$.~This led us to establish the minimal number of boundary conditions required to bound an energy estimate which coincides with the number of negative eigenvalues of M.~We further defined the general form of boundary conditions $\eqref{eq18}$ in terms of the incoming and outgoing characteristics,~which results in an energy bound provided that $\eqref{eq21}$ holds.~By returning to the boundary term $\eqref{eq9}$,~we proposed a set of commonly used boundary conditions which cancels or limits the terms that add growth to the energy rate $\eqref{eq8}$.~Moreover,~we constructed $\mathcal{S}$ $\eqref{eq27}-\eqref{eq31}$ for each boundary satisfying $\eqref{eq21}$ and demonstrated that they can be written on the form $\eqref{eq18}$.~In the next Section,~we implement these derived boundary conditions weakly such that a stable scheme results.
\section{The weak imposition of the boundary conditions}\label{sec3}
In this section,~we implement the boundary conditions $\eqref{eq24}$ weakly and show that they lead to the energy estimate $\eqref{eq26}$.~For this purpose,~we will mimic the continuous analysis above.~Let's consider a weak formulation of $\eqref{eq3}$ which will also lay foundation for the upcoming discrete analysis
\begin{align}
\mathcal{I}\text{U}_t + \mathcal{D}(\text{U})\text{U} &= L \left( \Sigma  (\mathcal{B}\text{U} -\textbf{g}) \right),  \qquad  \textbf{x} \in \Omega , \quad t>0, \label{eq33}   \\
\mathcal{I}\text{U}&= \textbf{f}, \qquad \qquad \qquad  \qquad \textbf{x} \in \Omega , \quad t= 0 .  \nonumber
\end{align}
Here,~$\Sigma$ is a penalty coefficient matrix yet to be determined such that we get the energy estimate and $\mathcal{B}U$ is the boundary operator given in $\eqref{eq24}$.~The operator $L(\cdot)$ is called the lifting operator \cite{Box2005,Arnold2001} and it is defined for any two continuous vector functions $\psi$,~$\phi$ as
\begin{align}
\int \limits_{\Omega} \psi^T L(\phi) dxdy &= \oint\limits_\Gamma \psi^T\phi ds. \label{eq34}
\end{align}
\par By applying the energy method,~$\eqref{eq33}$ becomes
\begin{align}
\frac{d}{dt}||\text{U}||_\mathcal{I}^2 + 2\mu ||\text{U}_y||^2_\mathcal{I} = -\oint\limits_{\Gamma} (u_{\bm{n}} u^2 + 2u_{\bm{n}}p - 2\mu uu_y n_y)ds + 2 \oint\limits_{\Gamma} \text{U}^T \Sigma (\mathcal{B}\text{U} - \textbf{g} ) ds, \label{eq35}
\end{align}
which is identical to $\eqref{eq8}$ with an additional boundary term.~As before,~let BT denote the first boundary integral in $\eqref{eq35}$.~Similarly,~let the penalty boundary terms in $\eqref{eq35}$ be denoted by PT.~Our ambition is to deduce $\Sigma$ such that the weakly imposed boundary conditions $\eqref{eq23}$ are dissipative.~We propose the following penalty coefficients for each boundary
\begin{align}
\Sigma_s &= \begin{bmatrix}
-\dfrac{v}{2} + \mu \partial_y^T & 0  \\ 0 & 0  \\ 0 & - 1 
\end{bmatrix} , \quad
\Sigma_n =\begin{bmatrix}
1 & 0 \\ 0& 1\\ 0 & 0 
\end{bmatrix} \label{eq36}, \\
\Sigma_w &= \begin{bmatrix}
-\dfrac{u}{2} \\ 0  \\ -1 
\end{bmatrix}, \qquad \qquad \qquad \Sigma_e = \begin{bmatrix}
1 \\ 0 \\ 0
\end{bmatrix} \nonumber,
\end{align}
where $\partial^T_y$ is the $y-$partial derivative operating in the left direction.~The penalty terms in $\eqref{eq35}$ with the coefficients $\eqref{eq36}$ simplifies to
\begin{align}
\text{PT} = &+\int \limits_{\Gamma_e} 2up dy -\int \limits_{\Gamma_w} (u^3 + 2up) dy  + \int \limits_{\Gamma_n}  [ u ( \alpha vu - 2\mu u_y) + 2vp] dx  \label{eq37} \\ &-\int \limits_{\Gamma_s}  (vu^2 - 2\mu uu_y + 2vp) dx \nonumber
\end{align} 
Therefore,~by substituting $\eqref{eq37}$ into $\eqref{eq35}$,~most of the boundary terms vanish and only the dissipative terms remain.~The energy rate $\eqref{eq35}$ then becomes
\begin{align}
\frac{d}{dt}||\text{U}||_\mathcal{I}^2 + 2\mu ||\text{U}_y||^2_\mathcal{I} = -\int \limits_{\Gamma_n}   ( 1 - \alpha )vu^2  dx -\int \limits_{\Gamma_e} u^3dy \leq 0,\label{eq38}
\end{align}
which is identical to $\eqref{eq25}$ and temporal integration leads to the estimate $\eqref{eq26}$.
\par Next,~we return to $\eqref{eq35}$ and rewrite BT in form $\eqref{eq17}$.~Moreover,~we rewrite $\mathcal{B}$U in the general form $\eqref{eq22}$ and set $\textbf{g} = \textbf{0}$.~Equation $\eqref{eq35}$ with this notation becomes
\begin{align}
\frac{d}{dt}||\text{U}||_\mathcal{I}^2 + 2\mu ||\text{U}_y||^2_\mathcal{I} = &-\oint \limits_\Gamma \big( W_+^T\Lambda_+ W_+ + W_-^T\Lambda_- W_- \big) ds \label{eq39}\\&+ \oint \limits_\Gamma \big[ \text{U}^T\Sigma ( W_-^T - \mathcal{S}W_+) + ( W_-^T - \mathcal{S}W_+)^T\Sigma^T\text{U} \big] ds.\nonumber
\end{align}
The choice $ \text{U}^T \Sigma = W_-^T \Lambda $ \cite{Nordstrom2017a} transforms $\eqref{eq39}$ to
\begin{align}
\frac{d}{dt}||\text{U}||_\mathcal{I}^2 + 2\mu ||\text{U}_y||^2_\mathcal{I} = &-\oint \limits_\Gamma \begin{bmatrix}
W_+ \\ W_- 
\end{bmatrix}^T \begin{bmatrix}
\Lambda_+ &\Lambda_ \mathcal{S} \\ \mathcal{S}^T \Lambda_- & -\Lambda_-
\end{bmatrix} \begin{bmatrix}
W_+ \\ W_- 
\end{bmatrix}ds. \label{eq40}
\end{align} 
Adding and subtracting $\oint\limits_\Gamma W^T_+ \mathcal{S}^T\Lambda_-\mathcal{S} W_+ ds$ leads to the simplification
\begin{align}
\frac{d}{dt}||\text{U}||_\mathcal{I}^2 + 2\mu ||\text{U}_y||^2_\mathcal{I} = &-\oint \limits_\Gamma 
W_+^T (
\Lambda_+ + \mathcal{S}^T\Lambda_-\mathcal{S} )W_+ds   \nonumber \\&+ \oint \limits_\Gamma \begin{bmatrix}
W_+ \\ W_- 
\end{bmatrix}^T \begin{bmatrix}
\mathcal{S}^T\Lambda_-\mathcal{S} &-\Lambda_-\mathcal{S} \\ -\mathcal{S}^T \Lambda_- & \Lambda_-
\end{bmatrix} \begin{bmatrix}
W_+ \\ W_- 
\end{bmatrix}ds, \label{eq41} \\ =&-\oint \limits_\Gamma 
W_+^T (
\Lambda_+ + \mathcal{S}^T\Lambda_-\mathcal{S} )W_+ds \nonumber  \\&+ \oint \limits_\Gamma 
(W_-- \mathcal{S} W_+)^T \Lambda_-(W_-- \mathcal{S} W_+)ds,   \nonumber
\end{align} 
which is identical to $\eqref{eq20}$ plus an additional dissipative term.~The right-hand side (RHS) of $\eqref{eq41}$ is non-positive since we computed $\mathcal{S}$ that satisfies $\eqref{eq21}$ in $\eqref{eq28}-\eqref{eq32}$.~The estimate $\eqref{eq41}$ will be the target for the upcoming semi-discrete SBP-SAT approximation.
\section{The semi-discrete SBP-SAT formulation}\label{sec4}
In this section,~we derive a stable numerical approximation of $\eqref{eq3}$.~We approximate the spatial derivatives using finite difference operators on SBP form while keeping the temporal derivative continuous leading to a semi-discrete formulation.~The boundary conditions $\eqref{eq23}$ are imposed weakly using SAT methods which imitates the boundary imposition in $\eqref{eq33}$ discretely.~By mimicking the continuous analysis,~we will show that the newly formulated SBP-SAT approximation is stable.~We begin by discretizing the domain $\Omega$ using $N\times M$ equidistant grid points $ (x_i,y_j)$ where $i = 1,2, \hdots, N$,~$j = 1,2,\hdots,M$.~Let $\textbf{U} =(\textbf{u}^T,\textbf{v}^T,\textbf{p}^T)^T $ be an approximation of the solution for $\eqref{eq3}$ where $\textbf{u},\textbf{v},\textbf{p}$ respectively contains the discrete version of $u,v$,~and $p$,~projected on the Cartesian grid.~They are packaged as $NM \times 1$ vectors,~for example,~$\textbf{u} = (u_{11},\hdots ,u_{1M},\hdots ,u_{N,1},\hdots ,u_{NM})^T$ where $u_{ij} = u(x_i,y_j)$.~Let $E_1$ and $E_{N,M}$ be zero matrices with one only at the top-left  and right-bottom  corner,~respectively. 
\par We define the finite difference operators on SBP form,~next.
\begin{definition}The matrix $D  $ is a first derivative SBP operator of order $s$ if 
\begin{align*}
D\textbf{u} = \textbf{u}_x + \mathcal{O}(h)^s  \quad \text{and} \quad D = P^{-1}Q
\end{align*}
where $P$ is a positive definite and symmetric matrix while $Q$ is almost skew-symmetric and satisfies $Q+Q^T = -E_1 + E_N$.
\end{definition}
\noindent Here,~$P$ is diagonal and it approximates the continuous integral
\begin{align*}
\bm{1}^T P \textbf{u} \approx \int\limits_\Omega u dx
\end{align*}
where $\bm{1}^T = (1,\hdots,1)^T$.~The operator $D$ satisfies the principle of the integration-by-parts discretely since
\begin{align*}
\textbf{u}^TP(D\textbf{v}) &= \textbf{u}^TQ\textbf{v} = \textbf{u}^T \big( ( E_N-E_1) - Q^T \big)\textbf{v} = \textbf{u}^T (E_N-E_1) \textbf{v} - \textbf{u}^T  Q^T P^{-1}P\textbf{v} \\& = (uv)_N - (uv)_1- (D\textbf{u})^TP\textbf{v}.
\end{align*}
Moreover,~it is $2s$-order accurate in the interior stencil and $s$ accurate near the boundaries where $s$ is the order of the truncation error.~We approximate the second-order derivative by applying first derivative operator twice resulting in a wider stencil and the excepted global order of accuracy of the scheme is $s+1$ \cite{Mattsson2004,Svard2006}.
\par We extend the SBP operators to 2D domain using the Kronecker product \cite{Svard2014}.~The following Kronecker product's properties are essential and will be used frequently in the upcoming analysis:~$(A\otimes B)(C\otimes D) = AC \otimes BD$,~$(A\otimes B)^{-1} = A^{-1} \otimes B^{-1} $,~and $(A \otimes B)^T = A^T \otimes B^T$ where $A,B,C,D$ are matrices of appropriate sizes and $A$,~$B$ are nonsingular.~By using subscripts $x,$~$y$ to differentiate the operators operating in the $x,y-$directions,~the 2D SBP operators becomes
\begin{align}
D_x = P^{-1}_x Q_x \otimes I_M, \quad D_y = I_N \otimes P^{-1}_y Q_y. \label{eq42}
\end{align}
Here,~$I_N$ and $I_M$ are unit matrices of size $N\times N $ and $M \times M$,~respectively.~We further introduce the following notation to keep the derivations neater
\begin{align*}
\textbf{P} &= I_3 \otimes P_x \otimes P_y, \quad \textbf{D}_x = I_3 \otimes D_x, \quad \textbf{D}_y= I_3 \otimes \mathcal{D}_y, \\ \textbf{Q}_x &= I_3 \otimes Q_x \otimes P_y, \quad \textbf{Q}_y = I_3 \otimes P_x \otimes Q_y,
\end{align*}
and the last two block-matrices satisfies the SBP property
\begin{align}
\textbf{Q}_x + \textbf{Q}_x^T = I_3 \otimes (E_N - E_1) \otimes P_y  \quad \text{and} \quad \textbf{Q}_y + \textbf{Q}_y^T = I_3 \otimes P_x \otimes (E_M-E_1). \label{eq43}
\end{align}
Furthermore,~$\textbf{P}$ defines the discrete $L_2$ semi-norm $||\textbf{V}||^2_{\bm{\mathcal{I}}\textbf{P}} = \textbf{V}^T \bm{\mathcal{I}} \textbf{P} \textbf{V} $ for a $3NM\times 1$ vector $\textbf{V}$ and $\bm{\mathcal{I}}$ is a discrete analogue of $\mathcal{I}$ in $\eqref{eq3}$.
\par The SBP-SAT formulation approximating $\eqref{eq33}$ is
\begin{align}
\bm{\mathcal{I}}\textbf{U}_t + \bm{\mathcal{D}}( \textbf{U} )\textbf{U} &= \sum \limits_{k \in \{n,e,s,w\}}  \textbf{P}^{-1}  \bm{\Sigma}_k ( I_j \otimes  \mathbb{P}_k )(\bm{\mathcal{B}}_k\textbf{U} - G),\quad  j = \begin{cases} 1 \quad \text{for } k \in \{e,w\} \\ 2 \quad \text{for } k \in \{n,s\} \end{cases}.   \label{eq44}\\
\bm{\mathcal{I}}\textbf{U}(0) &= F.   \nonumber
\end{align}
Here,~$\bm{\mathcal{D}}$ is the discrete version of the spatial operator in $\eqref{eq3}$ with the splitting $\eqref{eq4}$ and it is given by
\begin{align}
 \bm{\mathcal{D}}( \textbf{U} )  =  \frac{1}{2} \big[ \textbf{A} \textbf{D}_x  + \textbf{D}_x \textbf{A}  +\textbf{B} \textbf{D}_y  + \textbf{D}_y \textbf{B} \big]   - \mu \bm{\mathcal{I}} \textbf{D}_y^2 , \label{eq45} 
\end{align}
where
\begin{align*}
\textbf{A} = \begin{bmatrix} \text{diag}(\textbf{u}) & \bm{0} & \textbf{I} \\ \bm{0} & \bm{0} & \bm{0} \\ \textbf{I} & \bm{0} & \bm{0} \end{bmatrix}, \quad \textbf{B} = \begin{bmatrix} \text{diag}(\textbf{v}) & \bm{0} & \bm{0} \\ \bm{0} & \bm{0} & \textbf{I} \\ \bm{0} & \textbf{I} & \bm{0} \end{bmatrix},
\end{align*}
and $\textbf{I}$,~$\bm{0}$ are unit and zero matrices of size $NM \times NM$.~The RHS of the governing equation in $\eqref{eq44}$ denotes the weakly imposed boundary conditions using the SAT method and is analogue to the RHS of $\eqref{eq33}$.~The discrete boundary operator and penalty coefficients along the $k$th boundary mimicking their continuous counterparts $\eqref{eq24}$,~$\eqref{eq36}$ are respectively denoted by $\bm{\mathcal{B}}_k$ and $\bm{\Sigma}_k$.~They are explicitly defined as
\begin{align*}
\bm{\mathcal{B}}_n = \begin{bmatrix}\dfrac{\alpha}{2} \text{diag}(\textbf{v}) - \mu D_y & \bm{0} & \bm{0} \\ \bm{0} & \bm{0} & \textbf{I}
\end{bmatrix}, \quad \bm{\mathcal{B}}_s = \begin{bmatrix} \textbf{I} & \bm{0} & \bm{0} \\ \bm{0} & \textbf{I} & \bm{0}
\end{bmatrix},\quad \bm{\mathcal{B}}_e = \begin{bmatrix}
\bm{0} & \bm{0} & \textbf{I}
\end{bmatrix}, \quad \bm{\mathcal{B}}_w = \begin{bmatrix}
\textbf{I} & \bm{0} & \bm{0}
\end{bmatrix},
\end{align*}
and
\begin{align*}
\bm{\Sigma}_n &= \begin{bmatrix}
\textbf{I} & \bm{0} \\ \bm{0} & \textbf{I} \\ \bm{0} & \bm{0}
\end{bmatrix}, \quad \bm{\Sigma}_s = \begin{bmatrix}
-\dfrac{1}{2}\text{diag}(\textbf{v}) +  \mu D_y^T & \bm{0} \\ \bm{0} & \bm{0} \\ \bm{0} & -\textbf{I}
\end{bmatrix}, \quad \bm{\Sigma}_e = \begin{bmatrix}
\textbf{I}  \\ \bm{0}  \\ \bm{0} 
\end{bmatrix} ,\quad \bm{\Sigma}_w = \begin{bmatrix}
-\dfrac{1}{2} \text{d}(\textbf{u}) \\ \bm{0}  \\ -\textbf{I} 
\end{bmatrix}.
\end{align*}
Furthermore,~the vectors $G_k$ and $F$ respectively contains the pointwise boundary and initial data.~Note that we only require initial data for the horizontal velocity since it is the only term containing temporal derivative.~Therefore,~$F$ has the form $F = [\vec{f},\vec{0},\vec{0}]^T$ where the elements of $\vec{f}$ are given by data $\textbf{f}$ in $\eqref{eq3}$ projected on the grid points.~Lastly,~in $\eqref{eq44}$,~the diagonal matrices $\mathbb{P}_k$ are the quadrature rules approximating the boundary line integral in $\eqref{eq8}$ and are explicitly given by
\begin{align*}
\mathbb{P}_k =  \begin{cases} P_x \otimes E_1   &\text{south boundary}, \\
E_N \otimes P_y   &\text{east boundary}, \\
P_x \otimes E_M   &\text{north boundary}, \\
E_1 \otimes P_y   &\text{west boundary}. \end{cases}
\end{align*}
With the notation above,~the SBP properties $\eqref{eq43}$ can now be written as $\textbf{Q}_x + \textbf{Q}_x^T = I_3 \otimes (\mathbb{P}_e - \mathbb{P}_w)$ and $\textbf{Q}_y + \textbf{Q}_y^T = I_3 \otimes (\mathbb{P}_n - \mathbb{P}_s)$.
\begin{remark}
Index $j$ in $\eqref{eq44}$ must be chosen such that the matrix multiplication is possible and it is equal to the number of boundary conditions prescribed per boundary.~In the case of vertical boundaries,~$j = 1$ since we are imposing exactly one boundary condition at the inflow and outflow boundaries.~Meanwhile,~for the north and south boundaries,~$j = 2$ since we are prescribing two boundary conditions at each boundary. \label{remark6}
\end{remark}
\subsection{Stability}
To derive the discrete energy estimate that resembles the continuous counterpart $\eqref{eq26}$,~we employ the discrete energy method.~By multiplying $\eqref{eq44}$ with $2\textbf{U}^T\textbf{P}$,~we obtain
\begin{align} 
2\textbf{U}^T \textbf{P} \bm{\mathcal{I}}\textbf{U}_t + \textbf{U}^T \big( \textbf{P}\textbf{A} \textbf{D}_x +  \textbf{Q}_x \textbf{A} \big) \textbf{U} + \textbf{U}^T \big(  \textbf{P}\textbf{B} \textbf{D}_y  + \textbf{Q}_y \textbf{B} \big)\textbf{U} - 2\mu\textbf{U}^T \bm{\mathcal{I}}\textbf{Q}_y \textbf{D}_y \textbf{U} &= \textbf{PT}, \label{eq46}
\end{align}
where
\begin{align*}
\textbf{PT} = \sum\limits_{k \in \{s,e,n,w\}}  2 \textbf{U}^T\Sigma_k (I_j \otimes \mathbb{P}_k) (\bm{\mathcal{B}}_k\textbf{U}-G_k).
\end{align*}
\par Next,~we simplify the terms in $\eqref{eq46}$ separately.~Starting with the nonlinear advection term in the $x$-direction and applying property $\eqref{eq43}$ on its conservative term,~we rewrite it as
\begin{align}
 \textbf{U}^T \textbf{P}\textbf{A}(\textbf{D}_x \textbf{U}) + \textbf{U}^T \textbf{Q}_x \textbf{A} \textbf{U}  &= \textbf{U}^T \textbf{P}\textbf{A} (\textbf{D}_x \textbf{U}) + \textbf{U}^T(I_3 \otimes ( \mathbb{P}_e - \mathbb{P}_w) )\textbf{A}\textbf{U}\nonumber \\&- \textbf{U}^T \textbf{Q}_x^T \textbf{P}^{-1} \textbf{P}\textbf{A} \textbf{U}  \label{eq47}\\
  &= \textbf{U}^T \textbf{P}\textbf{A}(\textbf{D}_x \textbf{U}) + \textbf{U}^T(I_3 \otimes ( \mathbb{P}_e - \mathbb{P}_w) )\textbf{A}\textbf{U} - (\textbf{D}_x \textbf{U})^T \textbf{P}\textbf{A}\textbf{U}  \nonumber \\&=  \textbf{U}^T(I_3 \otimes ( \mathbb{P}_e - \mathbb{P}_w) )\textbf{A} \textbf{U} \nonumber.
\end{align}
Notice that the non-conservative indefinite terms above cancels owing to the flux splitting $\eqref{eq4}$,~and only the boundary terms remains.~Similarly,~the advection terms in the $y$-direction simplifies to
\begin{align}
\textbf{U}^T \textbf{P}\textbf{B} (\textbf{D}_y \textbf{U}) + \textbf{U}^T \textbf{Q}_y \textbf{B} \textbf{U} = \textbf{U}^T(I_3 \otimes ( \mathbb{P}_n - \mathbb{P}_s) )\textbf{B}\textbf{U}. \label{eq48}
\end{align}
Next,~by using the SBP property $\eqref{eq43}$,~we simplify the viscous term in $\eqref{eq46}$ as
\begin{align}
2\mu \textbf{U}^T \bm{\mathcal{I}}\textbf{Q}_y (\textbf{D}_y  \textbf{U} )&= 2\mu  \textbf{U}^T (I_3 \otimes ( \mathbb{P}_n - \mathbb{P}_s) )\bm{\mathcal{I}}(\textbf{D}_y \textbf{U}) - 2\mu  \textbf{U}^T \textbf{Q}_y^T \textbf{P}^{-1}\textbf{P}\bm{\mathcal{I}}(\textbf{D}_y \textbf{U}) \label{eq49} \\&= 2\mu  \textbf{U}^T (I_3 \otimes ( \mathbb{P}_n - \mathbb{P}_s) )\bm{\mathcal{I}}(\textbf{D}_y \textbf{U}) - 2\mu (\textbf{D}_y \textbf{U})^T \textbf{P}\bm{\mathcal{I}}(\textbf{D}_y \textbf{U}) \nonumber.
\end{align}
In $\eqref{eq49}$,~we get both the boundary and the dissipative volume term.~By substituting $\eqref{eq47}$,~$\eqref{eq48}$,\\~$\eqref{eq49}$ into $\eqref{eq46}$,~the discrete energy rate become
\begin{align}
\frac{d}{dt} ||\textbf{U}||_{\bm{\mathcal{I}} \textbf{P}}^2 + 2\mu||\textbf{D}_y \textbf{U} ||^2_{\bm{\mathcal{I}}\textbf{P}} &= \textbf{BT} + \textbf{PT} , \label{eq50}
\end{align}
where

\begin{align}
\textbf{BT} = &-\textbf{U}^T \big[ (I_3 \otimes ( \mathbb{P}_e - \mathbb{P}_w) )\textbf{A} + (I_3 \otimes ( \mathbb{P}_n - \mathbb{P}_s) )\textbf{B}\big] \textbf{U} +2\mu\textbf{U}^T (I_3 \otimes ( \mathbb{P}_n - \mathbb{P}_s) )\bm{\mathcal{I}}(\textbf{D}_y\textbf{U}) \nonumber \\
= &-\big[  \textbf{u}^T  ( \mathbb{P}_e - \mathbb{P}_w) \text{d}(\textbf{u}) \textbf{u} +  \textbf{u}^T  ( \mathbb{P}_n - \mathbb{P}_s) \text{d}(\textbf{v}) \textbf{u}  + 2\textbf{u}^T  ( \mathbb{P}_e - \mathbb{P}_w) \textbf{p}  + 2\textbf{v}^T  ( \mathbb{P}_n - \mathbb{P}_s) \textbf{p} \label{eq51} \\&- 2 \mu \textbf{u}^T  ( \mathbb{P}_n - \mathbb{P}_s) D_y \textbf{u} \big]  \nonumber
\end{align}
which is the discrete version of $\eqref{eq35}$.~We will mimic the continuous analysis here to ensure the RHS of $\eqref{eq51}$ have an appropriate sign such that we obtain a discrete estimate.
\par To proceed,~we first rewrite $\textbf{BT}$ in $\eqref{eq51}$ in the form that resembles $\eqref{eq10}$.~Let the pair $(N^k_x,N^k_y)$ be the discrete boundary normals as defined below.
\begin{definition}The discrete outward pointing boundary normals are given by the pair $N^k = (N_x^k,N_y^k)$
\begin{align}
(N_x^s,N_y^s)  &= (\bm{0},I_N \otimes -E_1), \label{eq52}\\
(N_x^e,N_y^e)  &= (E_N \otimes I_M,\bm{0}) , \nonumber  \\
(N_x^n,N_y^n)  &= (\bm{0},I_N \otimes E_M), \nonumber \\
(N_x^w,N_y^w)  &= ( -E_1 \otimes I_M,\bm{0}) \nonumber
\end{align}
\end{definition}
\noindent Using $\eqref{eq52}$,~$\textbf{BT}$ $\eqref{eq51}$ can now be written such that it discretely imitates $\eqref{eq10}$
\begin{align}
\textbf{BT} &= -\sum_{k \in \{ s,e,n,w\}} \text{Q}^T (I_4 \otimes \mathbb{P}_k )\textbf{M}_k\text{Q} .\label{eq53}
\intertext{where}
\textbf{M}_k &=  \begin{bmatrix}
\text{diag}(\textbf{u}_{\bm{n}}^k)& \bm{0} & N_x^k & -N_y^k
\\ \bm{0} & \bm{0} & N_y^k & \bm{0} \\ N_x^k & N_y^k & \bm{0} & \bm{0}\\
-N_y^k & \bm{0} &\bm{0} & \bm{0}
\end{bmatrix} , \quad \text{Q} =   \begin{bmatrix}
\textbf{u} \\ \textbf{v} \\ \textbf{p} \\ \mu D_y \textbf{u} ,
\end{bmatrix} \nonumber
\end{align}
and $\textbf{u}_{\bm{n}}^k = N_x^k\textbf{u} + N_y^k\textbf{v} $ is the discrete boundary normal velocity.~Since all the matrices in $\eqref{eq53}$ are diagonal then there are $NM$ decoupled nonlinear equations.~However,~the number of nonzero entries is equal to the number of boundary grid points due to the normals  $\eqref{eq52}$.~By noting the similarity in the structures of M in $\eqref{eq10}$ and $\textbf{M}$ in $\eqref{eq53}$,~we adopt the similar matrix eigenvalue decomposition $\eqref{eq13}$ in the discrete sense
\begin{align}
\textbf{M}_k= \textbf{X}_k \bm{\Lambda}_k \textbf{X}_k^T.  \label{eq54}
\end{align}
Here,~$\bm{\Lambda}_k = \text{d}\big(\bm{\lambda}_1^k,\bm{\lambda}_2^k,\bm{\lambda}_3^k,\bm{\lambda}_4^k \big)$ is a $4NM \times 4NM$ diagonal matrix containing the eigenvalues of $\textbf{M}_k$ and $\textbf{X}_k$ is the associated eigenvector block-matrix on the $k$th boundary.~Vectors $\bm{\lambda}^k_i$ contains pointwise eigenvalues of $\textbf{M}_k$ which are obtained by projecting $\eqref{eq11}$ and $\eqref{eq12}$ on the north,~south and east,~west boundary grid points,~respectively.~For the north and south boundaries,~$\bm{\lambda}_i^k$ and $\bm{X}^k$ are
\begin{align}
\bm{\lambda}_1^k &= \dfrac{\textbf{u}_{\bm{n}}^k}{2} - \sqrt{\left( \dfrac{\textbf{u}_{\bm{n}}^k}{2} \right)^2 + 1}, \quad \bm{\lambda}_2^k = - \vec{1}, \quad \bm{\lambda}_3^k = \vec{1}, \quad \bm{\lambda}_4^k = \dfrac{\textbf{u}_{\bm{n}}^k}{2} + \sqrt{\left( \dfrac{\textbf{u}_{\bm{n}}^k}{2} \right)^2 + 1}, \quad  k\in \{n,s \},  \label{eq55}\\
\textbf{X}_k &= \begin{bmatrix}
\text{diag}(\bm{\lambda}^k_1) & \bm{0} & \bm{0} & \text{diag}(\bm{\lambda}^k_4) \\ \bm{0} & \textbf{I} & \textbf{I} & \bm{0} \\ \bm{0} & N_y^k & N_y^k & \bm{0} \\ -N_y^k & \bm{0} & \bm{0} & -N_y^k
\end{bmatrix}.   \nonumber
\end{align}
Similarly,~at the east and west boundaries,~they are
\begin{align}
\bm{\lambda}_1^k &= \dfrac{\textbf{u}_{\bm{n}}^k}{2} - \sqrt{\left( \dfrac{\textbf{u}_{\bm{n}}^k}{2} \right)^2 + 1}, \quad \bm{\lambda}_2^k =  \vec{0}, \quad \bm{\lambda}_3^k = \vec{0}, \quad \bm{\lambda}_4^k =  \dfrac{\textbf{u}_{\bm{n}}^k}{2} + \sqrt{\left( \dfrac{\textbf{u}_{\bm{n}}^k}{2} \right)^2 + 1}, \quad  k\in \{e,w\}, \label{eq56} \\
\textbf{X}_k &= \begin{bmatrix}
\text{diag}(\bm{\lambda}^k_1) & \bm{0} & \bm{0} & \text{diag}(\bm{\lambda}^k_4) \\ \bm{0} & N_x^k & \bm{0} & \bm{0} \\ N_x^k & \bm{0} & \bm{0} & N_x^k \\ \bm{0} & \bm{0} & N_x^k & \bm{0}
\end{bmatrix},   \nonumber
\end{align}
where $\vec{1}$,~$\vec{0}$ respectively denote vector of ones and zeros.
\begin{remark} The square-roots and multiplications in $\eqref{eq55}$ and $\eqref{eq56}$ should be interpreted element-wise.~For subsequent analysis,~all operations involving $\bm{\lambda}_i$ should also be interpreted element-wise.
\end{remark}
\noindent By substituting $\eqref{eq54}$ into $\eqref{eq53}$ and defining the discrete characteristic variables $\textbf{W}^k = \textbf{X}_k \text{Q}^T $,~$\eqref{eq53}$ becomes
\begin{align}
\textbf{BT} = \sum\limits_{k \in \{ s,e,n,w \}  }  \textbf{W}^{k,T} (I_4 \otimes \mathbb{P}_k )\bm{\Lambda}_k \textbf{W}^k,   \label{eq57}
\end{align}
which mimics $\eqref{eq14}$ discretely.~We further divide it in terms of the positive and negative components as before in $\eqref{eq17}$
\begin{align}
\textbf{BT} = -\sum\limits_{k \in \{ s,e,n,w \} } \begin{bmatrix}
\textbf{W}^k_+ \\ \textbf{W}^k_- 
\end{bmatrix}^T  (I_{2j} \otimes \mathbb{P}_k )\begin{bmatrix}
\bm{\Lambda}^k_+ & \bm{0} \\ \bm{0} & \bm{\Lambda}^k_-
\end{bmatrix}  \begin{bmatrix}
\textbf{W}^k_+ \\ \textbf{W}^k_- 
\end{bmatrix}.  \label{eq58}
\end{align}
The variables in $\eqref{eq58}$ are the discrete analogues of $\eqref{eq15}$ and $\eqref{eq16}$.~For the north and south boundaries,~they are defined as
\begin{align}
\textbf{W}^k_+ &= \begin{bmatrix}
\textbf{v} + N_y^k \textbf{p} \\ \bm{\lambda}_4^k \textbf{u} - \mu N_y^kD_y \textbf{u}
\end{bmatrix},\quad  \textbf{W}^k_- = \begin{bmatrix}
\bm{\lambda}_1^k \textbf{u} - \mu N_y^kD_y \textbf{u} \\ \textbf{v} - N_y^k \textbf{p} 
\end{bmatrix} ,\label{eq59}  \\  \bm{\Lambda}^k_+ &= \begin{bmatrix}
\text{diag}\big( \bm{\lambda}^k_3/2 \big) & \bm{0} \\ \bm{0} & \text{diag}\big(\bm{\lambda}^k_4/((\bm{\lambda}_4^k)^2 + 2) \big)
\end{bmatrix} , \quad  
\bm{\Lambda}^k_- = \begin{bmatrix}
\text{diag}\left(\bm{\lambda}^k_1/((\bm{\lambda}_1^k)^2 + 2) \right) & \bm{0} \\ \bm{0} &\text{diag}\left( \bm{\lambda}^k_2/2 \right)
\end{bmatrix}  , \nonumber
\end{align}
For the east and west boundaries,~we have
\begin{align}
\textbf{W}^k_+ &= \begin{bmatrix}
\bm{\lambda}^k_4\textbf{u} + N_x^k \textbf{p} 
\end{bmatrix},\quad  \textbf{W}^k_- = \begin{bmatrix} \bm{\lambda}^k_1 \textbf{u} + N_x^k \textbf{p} 
\end{bmatrix}, \quad  \bm{\Lambda}^k_+ = \begin{bmatrix} \text{diag}\big(\bm{\lambda}^k_4/((\bm{\lambda}_4^k)^2 + 2) \big)
\end{bmatrix}, \label{eq60} \\ \bm{\Lambda}^k_- &= \begin{bmatrix}
\text{diag}\big(\bm{\lambda}^k_1/((\bm{\lambda}_1^k)^2 + 2) \big) 
\end{bmatrix}  . \nonumber
\end{align}
\noindent Note that in $\eqref{eq58}$,~we adjusted the dimensions of the unit matrix to $2j$ since $\textbf{W} = [\textbf{W}_+,\textbf{W}_-]^T$ has two vectors and index $j$ is defined as before in Remark \ref{remark6}.~Next,~we define the discrete version of $\eqref{eq22}$
\begin{align*}
\bm{\mathcal{B}}\textbf{U} = \textbf{W}_- - \bm{\mathcal{S}}\textbf{W}_+,
\end{align*}
and rewrite $\textbf{PT}$ in $\eqref{eq46} $ as
\begin{align}
\textbf{PT} = \sum\limits_{k \in \{ s,e,n,w \}  } &\Big[ \textbf{U}^T \bm{\Sigma}_k (I_j \otimes \mathbb{P}_k)\big(\textbf{W}_-^k - \bm{\mathcal{S}}_k\textbf{W}_+^k - G_k)  \label{eq61} \\&+ \big(\textbf{U}^T \bm{\Sigma}_k (I_j \otimes \mathbb{P}_k)\big(\textbf{W}_-^k - \bm{\mathcal{S}}_k\textbf{W}_+^k - G_k) \big)^T \Big]. \nonumber
\end{align}
Equation $\eqref{eq61}$ is the discrete version of the penalty term in $\eqref{eq39}$.~Therefore,~we make a similar choice $\textbf{U}^T \bm{\Sigma}_k = (\textbf{W}_-^k)^T\bm{\Lambda}_-^k$ and rewrite the rhs of $\eqref{eq50}$ as
\begin{align}
\frac{d}{dt} ||\textbf{U}||_{\bm{\mathcal{I}} \textbf{P}}^2 + 2\mu||\textbf{D}_y \textbf{U} ||^2_{\bm{\mathcal{I}}\textbf{P}} &=- \begin{bmatrix}
\textbf{W}_+^k\\\textbf{W}_-^k 
\end{bmatrix}^T  (I_{2j} \otimes \mathbb{P}_k ) \begin{bmatrix}
\bm{\Lambda}_+^k & \bm{\Lambda}_-^k \bm{\mathcal{S}} \\ \bm{\mathcal{S}}^T \bm{\Lambda}_-^k & - \bm{\Lambda}_-^k 
\end{bmatrix}\begin{bmatrix}
\textbf{W}_+^k\\ \textbf{W}_-^k 
\end{bmatrix}.  \label{eq62}
\end{align}
Adding and subtracting $(\textbf{W}_+^k)^T(I_j \otimes \mathbb{P}_k) [\bm{\mathcal{S}}^T_k\bm{\Lambda}_-^k\bm{\mathcal{S}}_k]\textbf{W}_+^k$ on the RHS of $\eqref{eq62}$ transforms the energy rate to
\begin{align}
\frac{d}{dt} ||\textbf{U}||_{\bm{\mathcal{I}} \textbf{P}}^2 + 2\mu||\textbf{D}_y \textbf{U} ||^2_{\bm{\mathcal{I}}\textbf{P}} = &-(\textbf{W}_+^k)^T(I_j \otimes \mathbb{P}_k) \big[ \bm{\Lambda}^k_+ +  \bm{\mathcal{S}}^T_k\bm{\Lambda}_k\bm{\mathcal{S}}_k \big] \textbf{W}_+^k  \label{eq63} \\&+ \begin{bmatrix}
\textbf{W}_+^k\\\textbf{W}_-^k 
\end{bmatrix}^T  (I_{2j} \otimes \mathbb{P}_k ) \begin{bmatrix}
 \bm{\mathcal{S}}^T_k\bm{\Lambda}_-^k\bm{\mathcal{S}}_k & -\bm{\Lambda}_-^k \bm{\mathcal{S}}   \\ -\bm{\mathcal{S}}^T \bm{\Lambda}_-^k &  \bm{\Lambda}_-^k 
\end{bmatrix}\begin{bmatrix}
\textbf{W}_+^k\\ \textbf{W}_-^k 
\end{bmatrix}  \nonumber \\ 
=  &-(\textbf{W}_+^k)^T(I_j \otimes \mathbb{P}_k) \big[ \bm{\Lambda}^k_+ +  \bm{\mathcal{S}}^T_k\bm{\Lambda}_k\bm{\mathcal{S}}_k \big] \textbf{W}_+^k  \nonumber \\&+  (\textbf{W}_-^k - \bm{\mathcal{S}}_k\textbf{W}_+^k)^T(I_j \otimes \mathbb{P}_k)  \bm{\Lambda}^k_- (\textbf{W}_-^k - \bm{\mathcal{S}}_k\textbf{W}_+^k) .  \nonumber
\end{align}
which is similar to $\eqref{eq41}$.~The first term on the RHS of $\eqref{eq63}$ is negative if we can find $\bm{\mathcal{S}}_k$ such that
\begin{align}
\bm{\Lambda}^k_+ +  \bm{\mathcal{S}}^T_k\bm{\Lambda}_k\bm{\mathcal{S}}_k \geq 0 \label{eq64}
\end{align}
which imitates $\eqref{eq21}$ discretely.~In $\eqref{eq27}-\eqref{eq32}$,~we computed the continuous analogue of $\bm{\mathcal{S}}_k$ at each boundary satisfying the continuous version of $\eqref{eq64}$.~Without loss of generality,~we assume that they also hold in the discrete setting as well.~The last term in $\eqref{eq63}$ is clearly negative and hence the energy rate is bounded.~Therefore,~time integration lead to the energy estimate that resembles $\eqref{eq26}$
\begin{align}
||\textbf{U}||_{\bm{\mathcal{I}} \textbf{P}}^2 + 2\mu \int\limits_0^T ||\textbf{D}_y \textbf{U} ||^2_{\bm{\mathcal{I}}\textbf{P}} dt \leq ||F||_{\bm{\mathcal{I}} \textbf{P}}^2 .\label{eq65}
\end{align}
\par Lastly,~we digress and consider the penalty terms in $\eqref{eq46}$.~Using the penalty coefficients given in $\eqref{eq44}$,~we show that the boundary conditions $\eqref{eq23}$ also lead to stability in the discrete setting.~The penalty coefficients and boundary operators in $\eqref{eq46}$ leads to
\begin{align}
\textbf{PT} = &+ \big[ \textbf{u}^T(- \mathbb{P}_w) \text{diag}(\textbf{u}) \textbf{u} + \textbf{u}^T  (\alpha  \mathbb{P}_n - \mathbb{P}_s) \text{diag}(\textbf{v}) \textbf{u} +  2\textbf{u}^T  ( \mathbb{P}_e - \mathbb{P}_w) \textbf{p} + 2\textbf{v}^T  ( \mathbb{P}_n - \mathbb{P}_s) \textbf{p} \label{eq66} \\&- 2 \mu \textbf{u}^T  ( \mathbb{P}_n - \mathbb{P}_s) D_y \textbf{u} \big], \nonumber
\end{align}
which is analogue to $\eqref{eq37}$.~Therefore,~substituting $\eqref{eq51}$ and $\eqref{eq66}$ into $\eqref{eq50}$ leads to cancellation of several boundary terms and the energy rate becomes
\begin{align}
\frac{d}{dt} ||\textbf{U}||_{\bm{\mathcal{I}} \textbf{P}}^2 + 2\mu||\textbf{D}_y \textbf{U} ||^2_{\bm{\mathcal{I}}\textbf{P}} = & -(1-\alpha ) \textbf{u}^T\mathbb{P}_n  \text{diag}(\textbf{v}) \textbf{u} -\textbf{u}^T\mathbb{P}_e \text{diag}(\textbf{u})\textbf{u}  , \label{eq67}
\end{align}
which is discretely identical to $\eqref{eq38}$ and here,~$\alpha  \in \{0,1\}$ as before.~We recall that d$(\textbf{u})>0$ and d$(\textbf{v})>0$ at the east and west boundaries respectively,~and therefore time integration leads to the estimate $\eqref{eq65}$ which proves that the approximation $\eqref{eq44}$ is stable.
\subsection{Null-space of the discrete spatial operator}
We revisit the spatial operator $\bm{\mathcal{D}}$ $\eqref{eq45}$ in this section.~Without the inclusion of the boundary conditions,~$\bm{\mathcal{D}}$ is singular and leads to non-unique or spurious solutions.~This is the reason for the majority of incompressible flow schemes create augmented equations as listed in the introduction.~However,~in this work we avoid this via the imposition of weak boundary conditions.~We therefore now prove the effect of the energy stable boundary conditions in removing the null-space of $\bm{\mathcal{D}}$.~This would be the case if all eigenvalues of $\bm{\mathcal{D}}$ were positive.~Following what was done for INS in \cite{Nordstrom2020c},~we first show that we expect the real parts of all eigenvalues to be positive in the case of the BL equations.~Following this,~we will also demonstrate this clearly by computing the eigenvalues of $\bm{\mathcal{D}}$ with and with-out the developed boundary conditions.~We begin by formulating the generalized nonlinear eigenvalue problem
\begin{align}
\bm{\mathcal{D}}(\textbf{U})\textbf{U} = \lambda \textbf{U}, \label{eq68}
\end{align}
where $\lambda$ denote the complex eigenvalues of spatial operator $\bm{\mathcal{D}}$ and are independent of the solution.~Here,~$\bm{\mathcal{D}}$ is the same as $\eqref{eq45}$ but with the SAT homogeneous boundary conditions included.~It is given by
\begin{align*}
\bm{\mathcal{D}}(\textbf{U}) &= \frac{1}{2} ( \textbf{A} \textbf{D}_x  + \textbf{D}_x \textbf{A}  +\textbf{B} \textbf{D}_y  + \textbf{D}_y \textbf{B} )   - \mu \bm{\mathcal{I}} \textbf{D}_y^2  - \sum \limits_{k \in \{n,e,s,w\}}  \textbf{P}^{-1}  \bm{\Sigma}_k ( I_j \otimes  \mathbb{P}_k )\bm{\mathcal{B}}_k  \\ &= 
\renewcommand*{\arraystretch}{1.5} 
\left[\begin{array}{c|c|c}
\begin{array}{c}
\dfrac{1}{2} \big[  \text{diag}(\textbf{u})D_x + D_x\text{diag}(\textbf{u}) +   \text{diag}(\textbf{v})D_y \\ + D_y\text{diag}(\textbf{v})\big] - \mu D_y^2 
  +  \dfrac{\alpha}{2}\mathcal{P}^{-1} \mathbb{P}_n  \text{diag}(\textbf{v}) \\ - \mu\mathcal{P}^{-1} \mathbb{P}_n  D_y  -\dfrac{1}{2} \mathcal{P}^{-1} \mathbb{P}_w \text{diag}(\textbf{u})\\-\dfrac{1}{2} \mathcal{P}^{-1} \mathbb{P}_s \text{diag}(\textbf{v}) + \mu \mathcal{P}^{-1}  D_y^T \mathbb{P}_s)
\end{array} & \bm{0} & D_x  + \mathcal{P}^{-1} \mathbb{P}_e  \\
\hline \bm{0}  & \bm{0} & D_y  + \mathcal{P}^{-1} \mathbb{P}_n  \\
\hline  D_x  - \mathcal{P}^{-1} \mathbb{P}_w  & D_y - \mathcal{P}^{-1} \mathbb{P}_s & \bm{0}
\end{array}\right] .
\end{align*}
To determine the sign of $\lambda$,~we employ the discrete energy method.~By multiplying $\eqref{eq68}$ with $\textbf{U}^*\textbf{P}$ from the left and adding to its complex transpose,~we obtain
\begin{align}
\textbf{U}^*\big[ \textbf{P}\bm{\mathcal{D}} + (\textbf{P}\bm{\mathcal{D}})^T] \textbf{U} &= (\lambda + \bar{\lambda})\textbf{U}^*\textbf{P}\textbf{U} = 2\text{Re}(\lambda)||\textbf{U}||^2_{\textbf{P}} , \label{eq69}
\end{align} 
where $\textbf{U}^*$ is the complex conjugate transpose of $\textbf{U}$.~For stability,~the left-hand side (LHS) of $\eqref{eq69}$ must be nonnegative or equivalently,~Re$(\lambda) >0$.~In Section $\ref{sec4}$,~we considered the energy analysis of the semi-discrete problem which includes the LHS of $\eqref{eq69}$ and therefore,~we will drop the temporal term and reuse the results for the spatial terms.~The LHS of $\eqref{eq69}$ becomes
\begin{align}
\textbf{U}^*\big[ \textbf{P}\bm{\mathcal{D}} + (\textbf{P}\bm{\mathcal{D}})^T] \textbf{U} &= 2\mu||\mathcal{D}_y \textbf{U} ||^2_{\bm{\mathcal{I}}\textbf{P}} -\textbf{BT} - \textbf{PT} \label{eq70} \\&= 2\mu||\mathcal{D}_y \textbf{U} ||^2_{\bm{\mathcal{I}}\textbf{P}} + ( 1 - \alpha) \textbf{u}^T\mathbb{P}_n  \text{diag}(\textbf{v}) \textbf{u} +\textbf{u}^T\mathbb{P}_e \text{diag}(\textbf{u})\textbf{u} > 0 .  \nonumber
\end{align}
where $\textbf{BT}$ and $\textbf{PT}$ are given in $\eqref{eq51}$ and $\eqref{eq66}$,~respectively.~Moreover,~they are preceded by the negative signs here because initially,~they were sitting on the RHS of the energy rate $\eqref{eq50}$.~Equation $\eqref{eq70}$ implies that $\text{Re}(\lambda) > 0$ in $\eqref{eq69}$ i.e.~all the eigenvalues of $\bm{\mathcal{D}}$ are on the right side of the complex plane for all $\textbf{U} \neq 0$.~Furthermore,~their sign is independent of the order of accuracy of the SBP operators and the number of computational grid points.
\par
\begin{figure}[H]
\centering
 \begin{subfigure}[t]{0.5\textwidth}
\includegraphics[scale = 0.6]{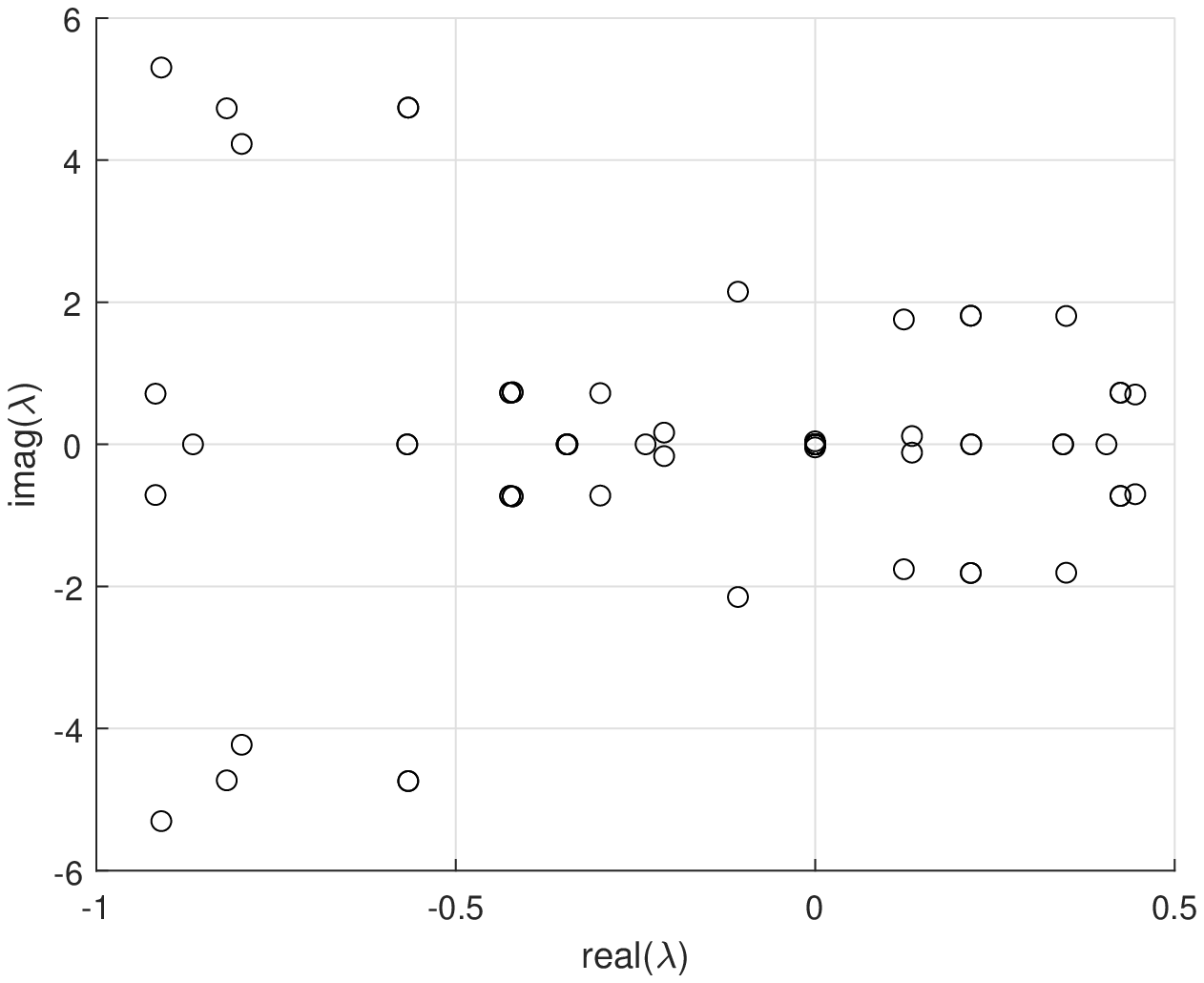}
\subcaption{}
\label{fig3a}
\end{subfigure}
\hspace{1cm}
\begin{subfigure}[t]{0.5\textwidth}
\includegraphics[scale = 0.6]{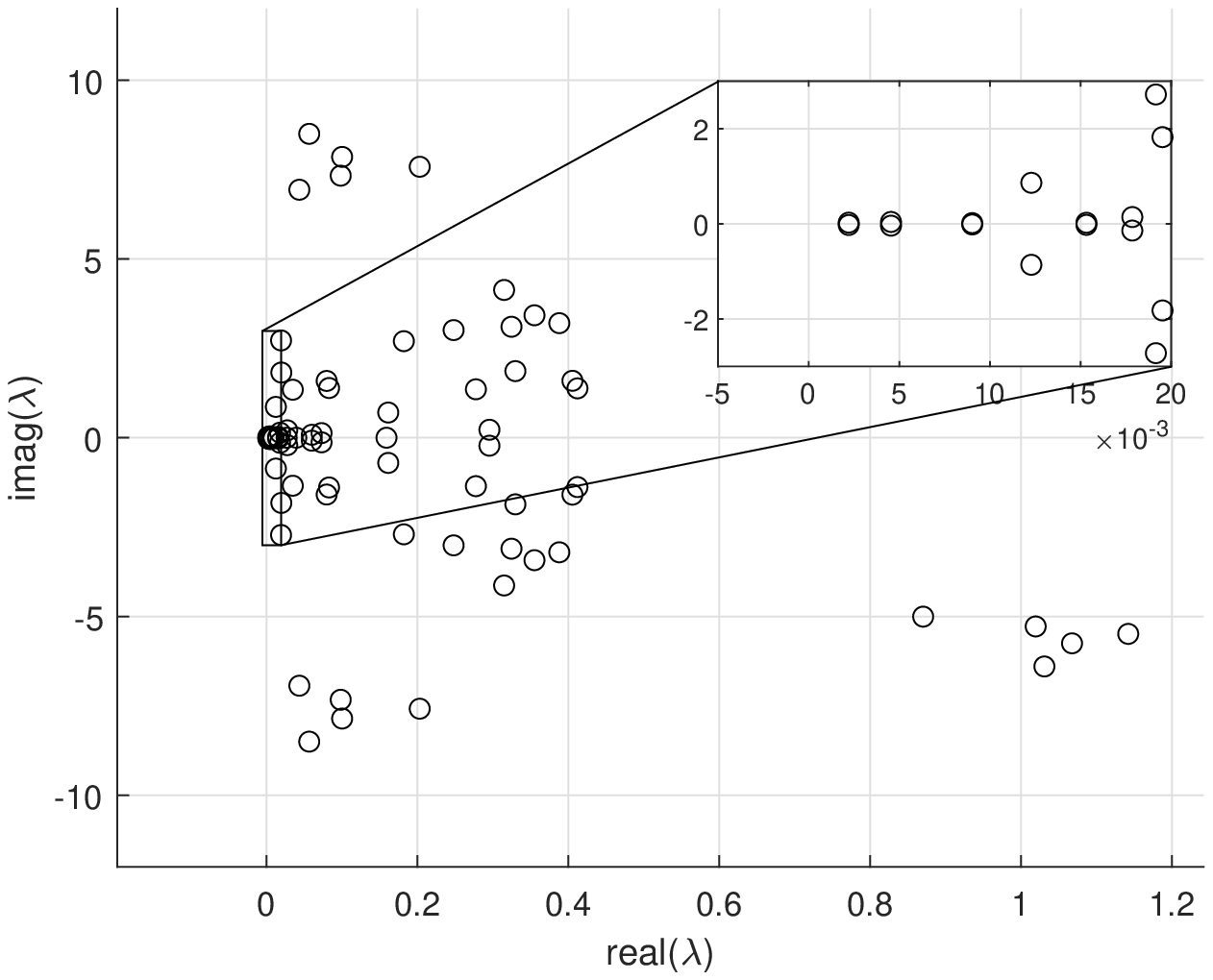}
\subcaption{}
\label{fig3b}
\end{subfigure}   
\caption{Eigenvalues of the spatial operator $\bm{\mathcal{D}}$ with (a) boundary conditions not included and (b) boundary conditions with $\alpha = 1$ are included.}
\label{fig3}
\end{figure}

 Next,~we inject the solution $\textbf{U} = [1,\hdots,1]^T$ in $\eqref{eq68}$ and use $4$th-order SBP operators to numerically compute the eigenvalues of $\bm{\mathcal{D}}$.~Further,~we choose $\alpha = 1$ in $\eqref{eq68}$,~this choice suffice to guarantee positive spectrum of $\mathcal{D}$ in $\eqref{eq70}$.~Setting $\alpha = 0$ will yield even more positive spectrum in $\eqref{eq70}$.~We consider two cases where:~the developed boundary conditions are (a) not included and (b) included in $\mathcal{D}$~.The eigenvalues resulting resulting from the first are depicted in Figure \ref{fig3a}.~As shown,~these contain both negative and zero values,~which will result in an unstable solution scheme.~As shown in Figure \ref{fig3b} however,~the addition of the developed BCs remedies the latter in full.~This is a key insight and contribution of this work.
\section{Temporal discretization and solution}\label{sec5}
To discretize the temporal derivative and progress the approximation $\eqref{eq44}$ in time,~we employ the first-order backward Euler method.~Let $\Delta t$ be the time-step size and $k$ denote the time-level.~The solution at two consecutive time-levels are denoted by $\textbf{U}^{k+1}$ and $\textbf{U}^k$.~The fully discrete approximation becomes
\begin{align}
\mathbb{F}(\textbf{U}^{k+1}) = \bm{\mathcal{I} } \dfrac{\textbf{U}^{k+1} - \textbf{U}^{k}}{\Delta t} +  \bm{\mathcal{D}}(\textbf{U}^{k+1})\textbf{U}^{k+1} = 0, \label{eq71}
\end{align}
where $\bm{\mathcal{D}}$ is the spatial operator (with the boundary conditions included) given in $\eqref{eq68}$.~Equation $\eqref{eq71}$ is a system of nonlinear equations which we linearise using Newtons method
\begin{align}
\textbf{U}^{k+1}= \textbf{U}^k - \textbf{J}_{\mathbb{F}(\textbf{U}^k)}^{-1} \mathbb{F}(\textbf{U}^k),  \label{eq72}
\end{align}
where $\textbf{J}_{\mathbb{F}}^{-1}$ is the inverse Jacobian matrix of $\mathbb{F}$.~Equation $\eqref{eq72}$ is solved iteratively till 
\begin{align*}
||\textbf{U}^{k+1} - \textbf{U}^k ||^2_{\textbf{P}}  < tol,
\end{align*}
where $tol$ is the specified tolerance.~The matrix $\textbf{J}_{\mathbb{F}}$ therefore comprises of the Jacobian matrix of the temporal term,~spatial terms,~and the boundary contributions i.e. 
\begin{align}
\textbf{J}_{\mathbb{F}} = \frac{\textbf{I}}{\Delta t} + \textbf{J}_{\bm{\mathcal{D}}_\Omega} - \textbf{J}_{\bm{\mathcal{D}}_\Gamma},  \label{eq73}
\end{align}
where
\begin{align*}
\textbf{J}_{\bm{\mathcal{D}}_\Omega} = 
\left[\begin{array}{ccc}
\begin{array}{c}
\dfrac{1}{2}\left(\text{diag}(\textbf{u})D_x + D_x \text{diag}(\textbf{u})  \right) \\ + 2D_x\text{diag}(\textbf{u}) + \text{diag}(\textbf{v}) D_y \\+ D_y\text{diag}(\textbf{v}) -\mu D_y^2 
\end{array} &  D_y\text{diag}(\textbf{u}) & D_x  \\ \bm{0} & \bm{0} & D_y \\ D_x & D_y & \bm{0}
\end{array} \right]
\end{align*}
\begin{align*}
\textbf{J}_{\bm{\mathcal{D}}_\Gamma} &=  \textbf{P}^{-1}(I_3 \otimes \mathbb{P}_n) \begin{bmatrix} \dfrac{\alpha}{2} \text{diag}(\textbf{v})-\mu D_y & \dfrac{\alpha}{2} \text{diag}(\textbf{u})& \bm{0} \\ \bm{0} & \bm{0} & \textbf{I}   \\ \bm{0} & \bm{0} & \bm{0} \end{bmatrix}  + \textbf{P}^{-1}(I_3 \otimes \mathbb{P}_e) \begin{bmatrix} \bm{0} & \bm{0} & \textbf{I} \\ \bm{0} & \bm{0} & \bm{0}   \\ \bm{0} & \bm{0} & \bm{0} \end{bmatrix} \\&+ \textbf{P}^{-1} \begin{bmatrix}
-\dfrac{1}{2}\text{diag}(\textbf{v}) + \mu D_y^T & -\dfrac{1}{2}\text{diag}(\textbf{u}) & \bm{0} \\ \bm{0} & \bm{0} &\bm{0} \\ \bm{0} & -\textbf{I} & \bm{0}
\end{bmatrix}(I_3 \otimes \mathbb{P}_s) \\& + \textbf{P}^{-1}(I_3 \otimes \mathbb{P}_w) \begin{bmatrix} -\text{diag}(\textbf{u}) & \bm{0} & \bm{0} \\ \bm{0} & \bm{0} & \bm{0}   \\ -\textbf{I} & \bm{0} & \bm{0} \end{bmatrix}.
\end{align*}
\section{Numerical experiments}\label{sec6}
We start by verifying the accuracy of the approximation scheme $\eqref{eq44}$,~and later move on to the comparison with the Blasius and INS solutions.
\subsection{Order of accuracy}
To compute the convergence rates,~we employ the method of manufactured solution \cite{Roache} on a compact domain $\Omega \in [0,1] \times [0,1]$.~The manufactured solution we choose is
\begin{align}
u &= \cosh(x)\sinh(y)e^{\mu t},  \quad
v = -\sinh(x)\cosh(y)e^{\mu t} , \quad 
p = \frac{1}{2} \sinh^2(x)e^{2\mu t}  \label{eq74},
\end{align}
where $\mu = 0.01$ and it satisfies $\eqref{eq2}$ exactly.~We impose Robin condition at the north boundary.~Further,~the boundary and initial data are sourced from $\eqref{eq74}$.~For temporal discretization,~we use first-order Backward Euler with time-step size $\Delta t  = 1\text{e}-04$ and set $tol = 1\text{e}-08$ for successive Newton's iterations.~We chose this small time-step to discard any temporal errors and the computations are terminated at $t = 1$.~Spatial derivatives are discretized using finite difference SBP $(2s,s)$-accurate operators where $s \in \{1,2,3\}$ is the accuracy near the boundaries.~The rate of convergence is computed as 
\begin{align*}
q = \log_{10} \left(\frac{||\textbf{e}^{h_1}||_{P_x \otimes P_y}^2}{||\textbf{e}^{h_2}||_{P_x \otimes P_y}^2} \right) \Big/ \log_{10} \left( \frac{h_1}{h_2} \right)
\end{align*}
where $||\textbf{e}||_{P_x \otimes P_y}^2$ is the $L_2$-norm of pointwise errors of the numerical and analytical solutions.~The mesh-spacing corresponding to the coarse and fine meshes are denoted by $h_1$ and $h_2$,~respectively.~The convergence rates for different orders of accuracy are presented in Tables~$\ref{Table1},\ref{Table2},\ref{Table3},\ref{Table4}$ and they coincide with the theoretical order of convergence. 
\begin{table}[H]
\centering
\begin{tabular}[t]{p{2cm} p{2cm} p{2cm} p{2cm} p{2cm} p{2cm} p{2cm} }
\hline
\multicolumn{7}{c}{$u$-velocity} \\
\hline
$N=M$&\multicolumn{2}{c}{SBP (2,1)}&\multicolumn{2}{c}{SBP (4,2) }&\multicolumn{2}{c}{SBP (6,3)}\\
\hline
& $||\textbf{e}||$ & $q$ & $||\textbf{e}||$ & $q$ & $||\textbf{e}||$ & $q$ \\ 
21 &0.0318 &- &0.0030 &- &5.9106e-04 &- \\ 41&0.0079&2.0829 &3.1441e-04&3.3769 &3.3657e-05 &4.2834\\ 61&0.0032 &2.2951 &7.3585e-05 &3.6553 &4.5945e-06 &5.0120\\81 &0.0016&2.4443 &2.4828e-05&3.8313 &1.0355e-06 &5.2543 \\ 
\hline 
Theoretical order  & & 2 & & 3 & & 4 \\
\hline
\end{tabular}
\caption{The $l_2$ norm of errors and the global order of accuracy of the approximation $\eqref{eq44}$ for the $u$-velocity using the different SBP operators.}
\label{Table1}
\end{table}
\begin{table}[H]
\centering
\begin{tabular}[t]{p{2cm} p{2cm} p{2cm} p{2cm} p{2cm} p{2cm} p{2cm} }
\hline
\multicolumn{7}{c}{$v$-velocity} \\
\hline
$N=M$&\multicolumn{2}{c}{SBP (2,1)}&\multicolumn{2}{c}{SBP (4,2) }&\multicolumn{2}{c}{SBP (6,3)}\\
\hline
& $||\textbf{e}||$ & $q$ & $||\textbf{e}||$ & $q$ & $||\textbf{e}||$ & $q$ \\ 
21 &0.0907 &- &0.0098 &- &0.0028 &- \\ 41&0.0210&2.1870 &0.0015&2.7727&2.0366e-04 &3.9257\\ 61&0.0092 &2.0861 &4.3612e-04 &3.1544 &3.0721e-05 &4.7609\\81 &0.0050&2.1503 &1.6813e-04&3.3613 &7.6274e-06 &4.9130 \\ 
\hline
Theoretical order & & 2 & & 3 & & 4 \\
\hline
\end{tabular}
\caption{The $l_2$ norm of errors and the global order of accuracy of the approximation $\eqref{eq44}$ for the $v$-velocity using the different SBP operators.}
\label{Table2}
\end{table}
\begin{table}[H]
\centering
\begin{tabular}[t]{p{2cm} p{2cm} p{2cm} p{2cm} p{2cm} p{2cm} p{2cm} }
\hline
\multicolumn{7}{c}{pressure} \\
\hline
$N=M$&\multicolumn{2}{c}{SBP (2,1)}&\multicolumn{2}{c}{SBP (4,2) }&\multicolumn{2}{c}{SBP (6,3)}\\
\hline
& $||\textbf{e}||$ & $q$ & $||\textbf{e}||$ & $q$ & $||\textbf{e}||$ & $q$ \\ 
21 &0.0159 &- &0.0021 &- &5.2357e-04 &- \\ 41&0.0038&2.1549 &1.9608e-04&3.5691 &1.6050e-04 &5.2088\\ 61&0.0016&2.1989 &4.7083e-05 &3.6000&2.1108e-06 &5.1061\\81 &8.4464e-04&2.2528 &1.6355e-05&3.7287&6.9588e-07 &3.9131\\ 
\hline
Theoretical order  & & 2 & & 3 & & 4 \\
\hline
\end{tabular}
\caption{The $l_2$ norm of errors and the global order of accuracy of the approximation $\eqref{eq44}$ for the pressure using the different SBP operators.}
\label{Table3}
\end{table}
\begin{table}[H]
\centering
\begin{tabular}[t]{p{2cm} p{2cm} p{2cm} p{2cm} p{2cm} p{2cm} p{2cm} }
\hline
\multicolumn{7}{c}{Solution vector U} \\
\hline
$N=M$&\multicolumn{2}{c}{SBP (2,1)}&\multicolumn{2}{c}{SBP (4,2) }&\multicolumn{2}{c}{SBP (6,3)}\\
\hline
& $||\textbf{e}||$ & $q$ & $||\textbf{e}||$ & $q$ & $||\textbf{e}||$ & $q$ \\ \hline 
21 &0.0029 &- &0.0104&- &0.0029 &- \\ 41&0.0227&2.1743&0.0016&2.8300 &2.0705e-04 &3.9576\\ 61&0.0098 &2.1125 &4.4478e-04 &3.1770 &3.1134e-05&4.7688\\81 &0.0053&2.1676&1.7074e-04&3.3763&7.7287e-06 &4.9135\\
\hline
Theoretical order & & 2 & & 3 & & 4 \\
\hline
\end{tabular}
\caption{The $l_2$ norm of errors and the global order of accuracy of the approximation $\eqref{eq44}$ for all variables $\text{U}=[u,v,p]^T$ using the different SBP operators.}
\label{Table4}
\end{table}
\subsection{Blasius boundary layer}
\begin{center}
\begin{figure}[H]
\captionsetup[figure]{justification=centering}
\centering
\includegraphics[scale=0.6]{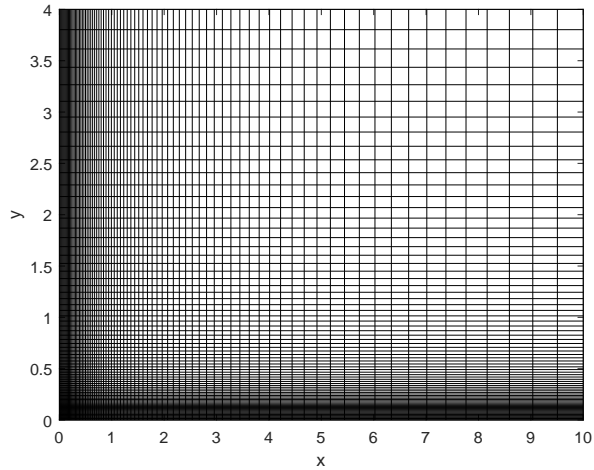}
\caption{Nonuniform computational mesh.}
\label{fig4}
\end{figure}
\end{center}

Viscous flow over a flat plate as illustrated in Figure $\ref{fig1}$ is finally modelled.~When encountering the plate's leading edge,~the fluid near the solid wall slows down due to the no slip condition.~Outside the boundary region,~the fluid's speed increases rapidly in the vertical direction until it reach the stream velocity leading to the formation of the boundary layer.~As a result,~the velocity gradients are the steepest near the leading edge and the plate's surface.~The thickness of this boundary layer $\delta(x)$ grows as a function of distance from the leading edge.~To resolve it effectively,~we employ nonuniform stretched meshes in the vicinity of the solid surface.~The use of SBP finite difference operators on nonuniform computational grids requires a consistent coordinate transformation that preserves the overall accuracy of the approximation scheme \cite{Petersson2015,Gong2007,Nordstrom2001}.~In \cite{Alund2019},~a simplified framework which encapsulates coordinate transformations into the SBP operators was developed.~This framework bypasses the need to first transform $\eqref{eq44}$ into curvilinear coordinates and subsequently apply the traditional SBP operators.~For coordinates stretching,~we use the continuous hyperbolic trigonometric functions such that the mesh is saturated in the region where the velocity gradients are the steepest as depicted in Figure $\ref{fig4}$.~We use
\begin{align*}
x(\xi,\eta) = x_0 + x_1 \frac{\sinh(\beta \xi)}{\sinh(\beta)}, \quad y(\xi,\eta) = y_0 + y_1 \frac{\sinh(\beta \eta)}{\sinh(\beta)},
\end{align*}
where $(\xi,\eta) \in [0,1]^2$ are the coordinates of the transformed regular domain,~$\beta$ is the stretching factor and we set it to $\beta = 4$.~Moreover,~$x_0,y_0$ and $x_1,y_1$ denotes the minimum and maximum of the values of the physical rectangular domain such that the pairs $(x_0,y_0)$ and $(x_0,y_1)$ are the coordinates of the left lower and upper corners.~Similarly,~the coordinates of the right lower and upper corners are $(x_1,y_0)$ and $(x_1,y_1)$,~respectively.
\par
\begin{center}
\begin{figure}[H]
\centering
\includegraphics[scale=0.6]{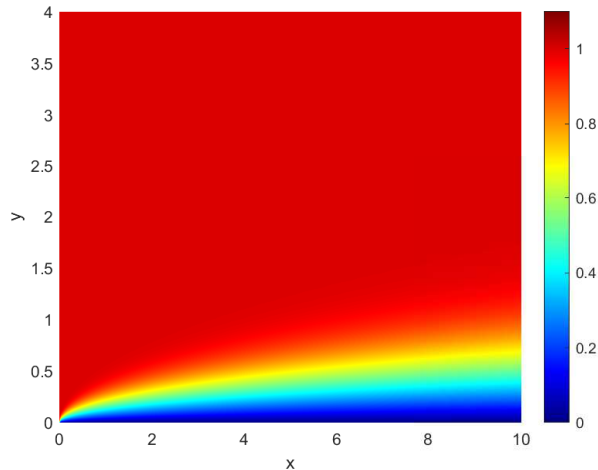}
\caption{The horizontal velocity distribution on the entire domain.}
\label{fig5}
\end{figure}
\end{center}

\begin{figure}[H]
\captionsetup[subfigure]{justification=centering}
 \begin{subfigure}[b]{0.5\textwidth}
 \centering
\includegraphics[scale=0.5]{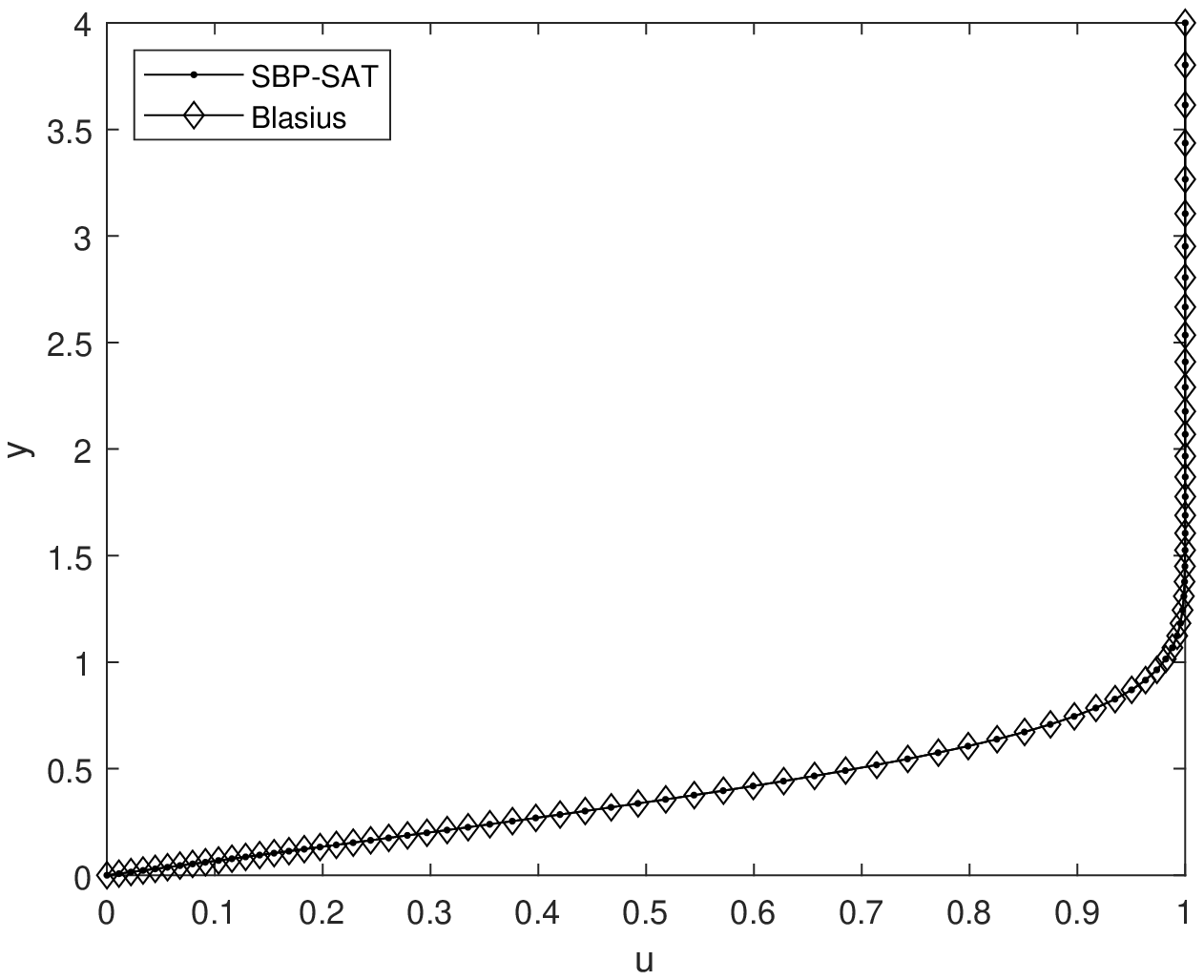}
\caption{}
     \label{fig6a}
     \end{subfigure}
     \hspace{1cm}
     \begin{subfigure}[b]{0.5\textwidth}
     \centering
\includegraphics[scale=0.5]{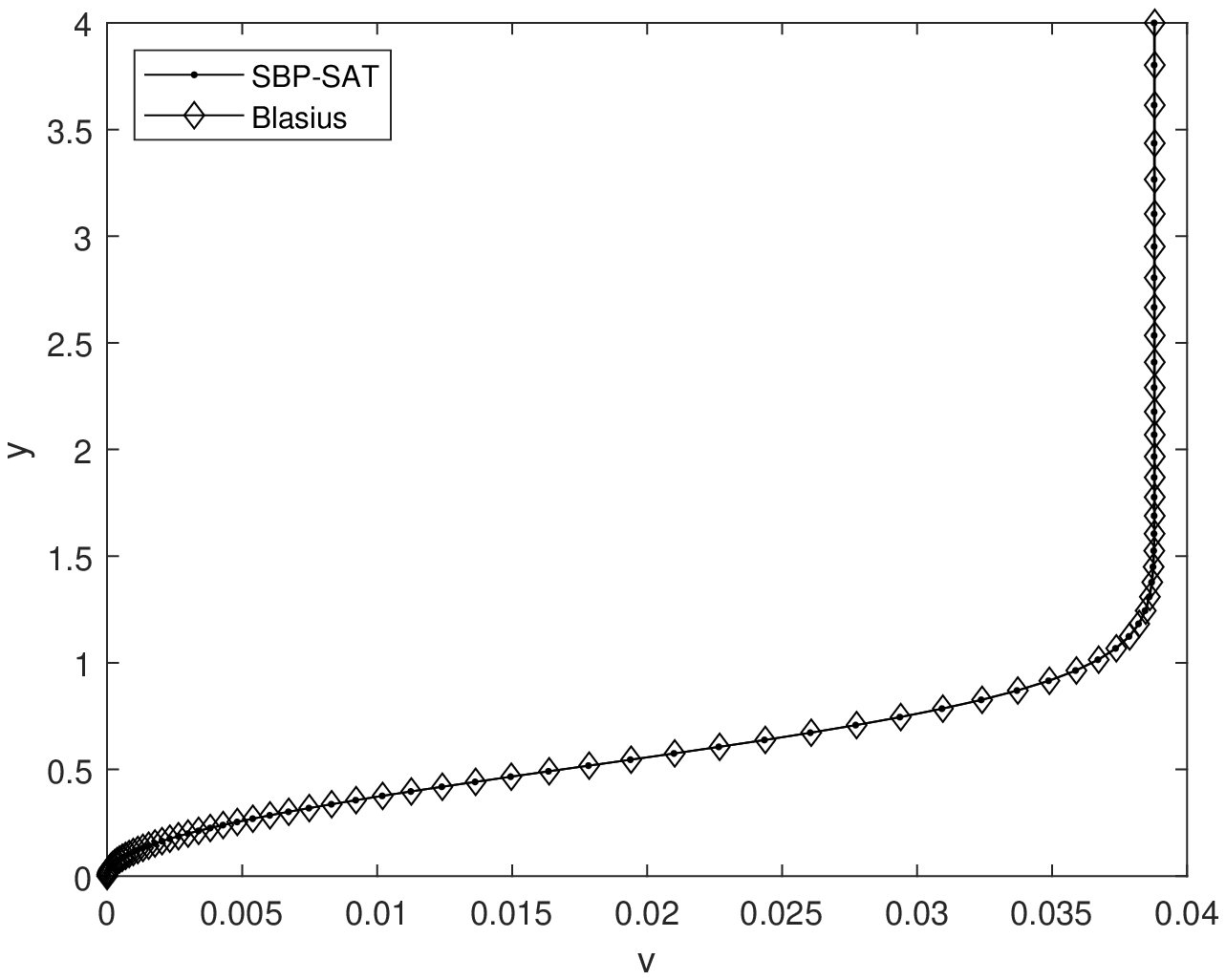}
     \caption{}
     \label{fig6b}
     \end{subfigure}   
     \caption{The SBP-SAT approximation for the IBL equations compared with Blasius solution along the line $x \approx 5$ with computations starting at the plate's leading edge.~(a) $u$-velocity profile and (b) $v$-velocity profile.}
     \label{fig6}
\end{figure}

 Let's consider the steady version of $\eqref{eq2}$ on the domain $\Omega \in [0,10] \times [0,4]$ and the stable approximation $\eqref{eq44}$.~We discretize $\Omega$ using $N = M = 80$ points as depicted in Figure.~$\ref{fig4}$ and set boundary data $U_\infty = 1$,~$p_\infty = 0$.~At the north boundary,~we consider Neumann boundary condition (i.e.~set $\alpha = 0$ in $\eqref{eq44}$) since we only know the $u$ velocity gradient in the freestream.~This choice however does not affect the positive definiteness of the resulting coefficient matrix as shown in $\eqref{eq70}$.~The continuous derivatives are approximated using 3rd-order accurate SBP operators.~Starting with the initial guess $\textbf{U}^1 =[\textbf{u}^1,\textbf{v}^1,\textbf{p}^1]^T =  [1,\hdots,1,0,\hdots,0,0\hdots,0]^T$,~we iterate $\eqref{eq72}$ progressively until we reach the steady state solution which is measured by
\begin{align*}
||Res_k||_{\textbf{P}}^2\leq  10^{-8} ||Res_1||_{\textbf{P}}^2
\end{align*}
where $Res_k$ is the residual (comprising of the spatial terms) at the $k$th time-level.~There are no restrictions on the time-step size and we set it to $\Delta t = 0.01$.~Moreover,~we set $\mu = 0.01$ such that $\delta \ll l$.~Figure~$\ref{fig5}$ shows the velocity distribution on the entire computational domain,~with a fully developed boundary layer.
Equation $\eqref{eq2}$ with (\textit{i}) $p_x = 0$,~(\textit{ii}) boundary conditions $\eqref{eq23}$ with $\alpha = 0$,~and (\textit{iii}) $U_\infty= \text{constant}$ has a well-known time-independent solution called the Blasius solution.~This similarity method-based solution reduces $\eqref{eq2}$ to a nonlinear ordinary differential equation which is then solved numerically (see \ref{Appendix}).~Therefore,~we will use this case to validate the SBP-SAT approximation $\eqref{eq44}$.~We compute $\eqref{eq71}$ till we reach the steady solution and then compare the results with the Blasius solution along particular vertical cross-sections on the domain.~Figure $\ref{fig6a}$ and $\ref{fig6b}$ shows the velocity profile along $x \approx 5$.~As shown,~our numerical approximation compares very well with the Blasius solution.
\begin{center}
\begin{figure}[H]
\captionsetup[subfigure]{justification=centering}
 \begin{subfigure}[b]{0.5\textwidth}
\includegraphics[scale=0.5]{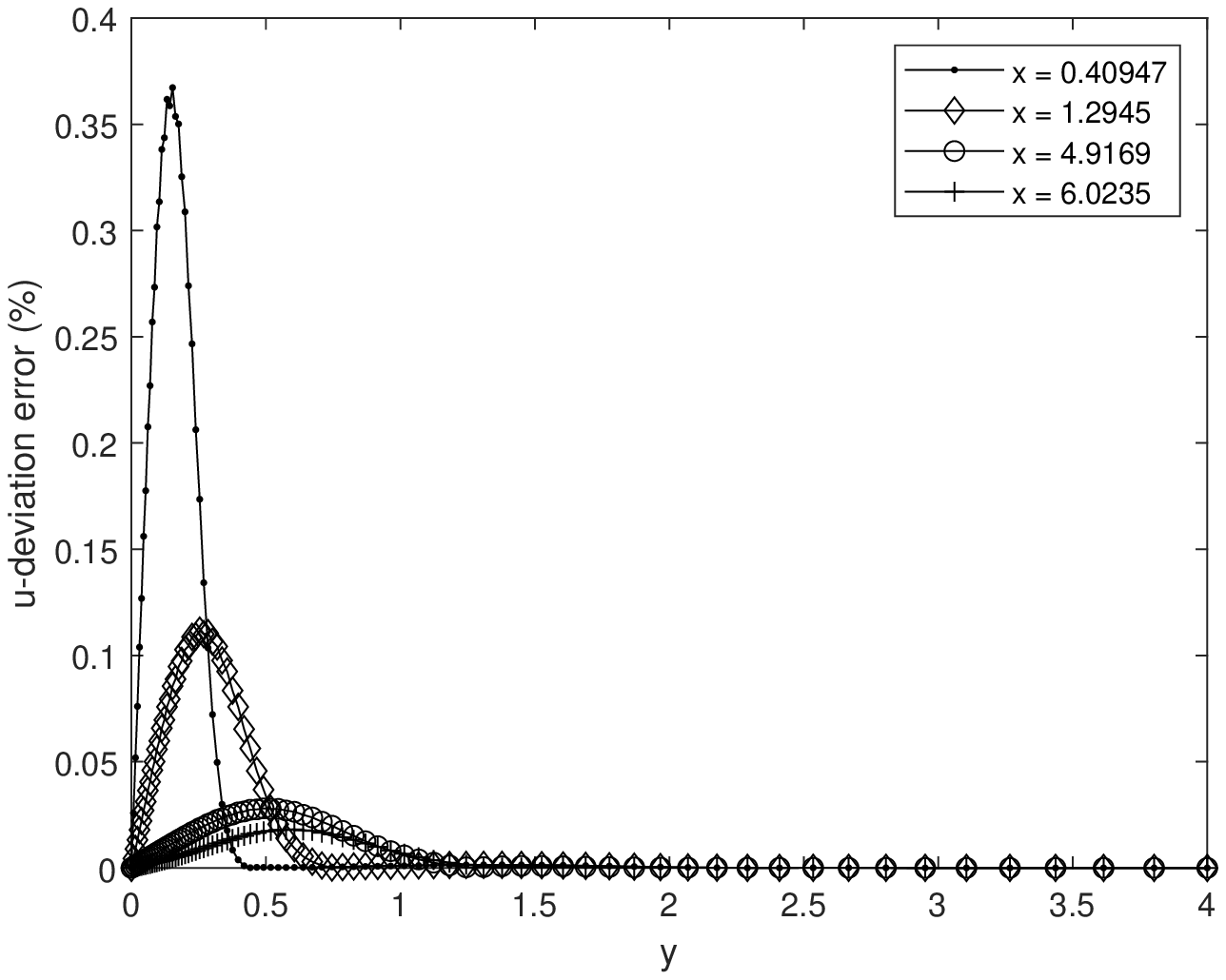}
\caption{}
     \label{fig7a}
     \end{subfigure}
     \hspace{1cm}
     \begin{subfigure}[b]{0.5\textwidth}
     \centering
\includegraphics[scale=0.5]{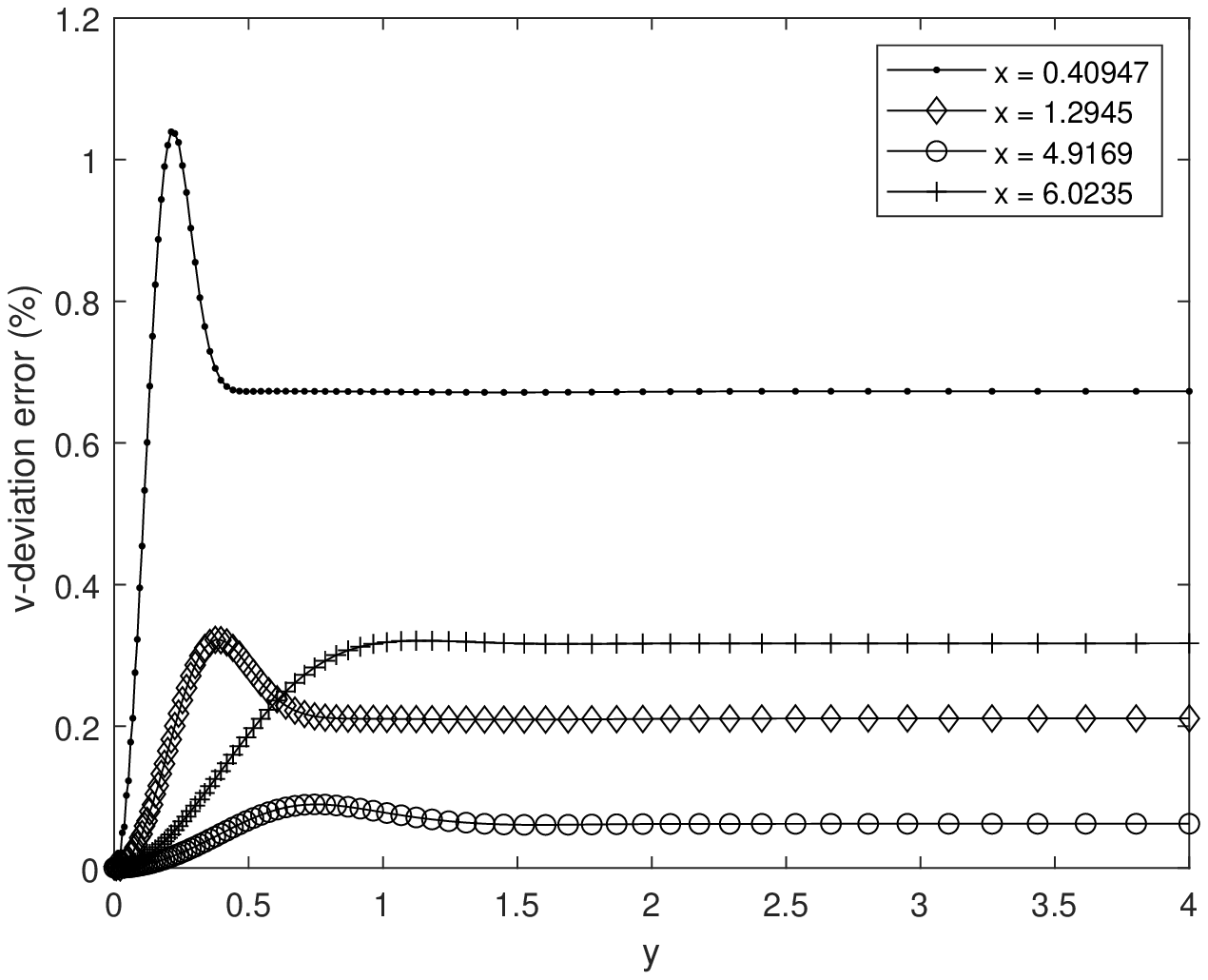}
     \caption{}
     \label{fig7b}
     \end{subfigure}   
     \caption{The deviation errors of the (a) $u$-velocity and (b) $v$-velocity at various $x$ position along the plate with computations starting at the leading edge.}
     \label{fig7}
\end{figure}
\end{center}
Next,~we include more vertical cross-sections across the domain and compare the errors between the two solution as shown in Figure~$\ref{fig7a}$ and $\ref{fig7a}$.~The point-wise errors are computed as
\begin{align}
\text{error}|_u = \frac{|u - u_{\text{B}}|}{|\max(u_{\text{B}})|} \times 100\%.  \label{eq75}
\end{align}
where $u$ and $u_{\text{B}}$ denote the SBP-SAT approximation and the Blasius solution.~Similarly,~we use $\eqref{eq75}$ to compute the deviations for $v$-velocity.~As shown on both profiles,~the errors are more dominant towards the leading edge and they dissipate downstream.~This is however thought to be due to the singularity at the leading edge of the plate i.e.~gradients in $u$ tend to infinity here \cite{Capatina2021}.~To overcome this,~we truncate $\Omega$ such that it excludes the tip of the plate and start the computations at a point $x_0$ on the domain as illustrated in Figure~$\ref{fig8a}$.~Instead of using $U_\infty$ as inflow data,~we use the Blasius solution evaluated at $x_0$ such that the gradients with respect to x are not large.~We choose $x_0 = 2$ as shown in Figure~$\ref{fig8b}$.
\begin{center}
\begin{figure}[H]
\captionsetup[subfigure]{justification=centering}
\centering
 \begin{subfigure}[b]{0.3\textwidth}
 \centering
 \psfrag{t0}{$y$} 
 \psfrag{t1}{$\delta(x)$} 
\psfrag{t2}{$x_0$}
\psfrag{t3}{$0$}
\psfrag{t5}{$l$}
\psfrag{t6}{$x$}
\includegraphics[scale=0.7]{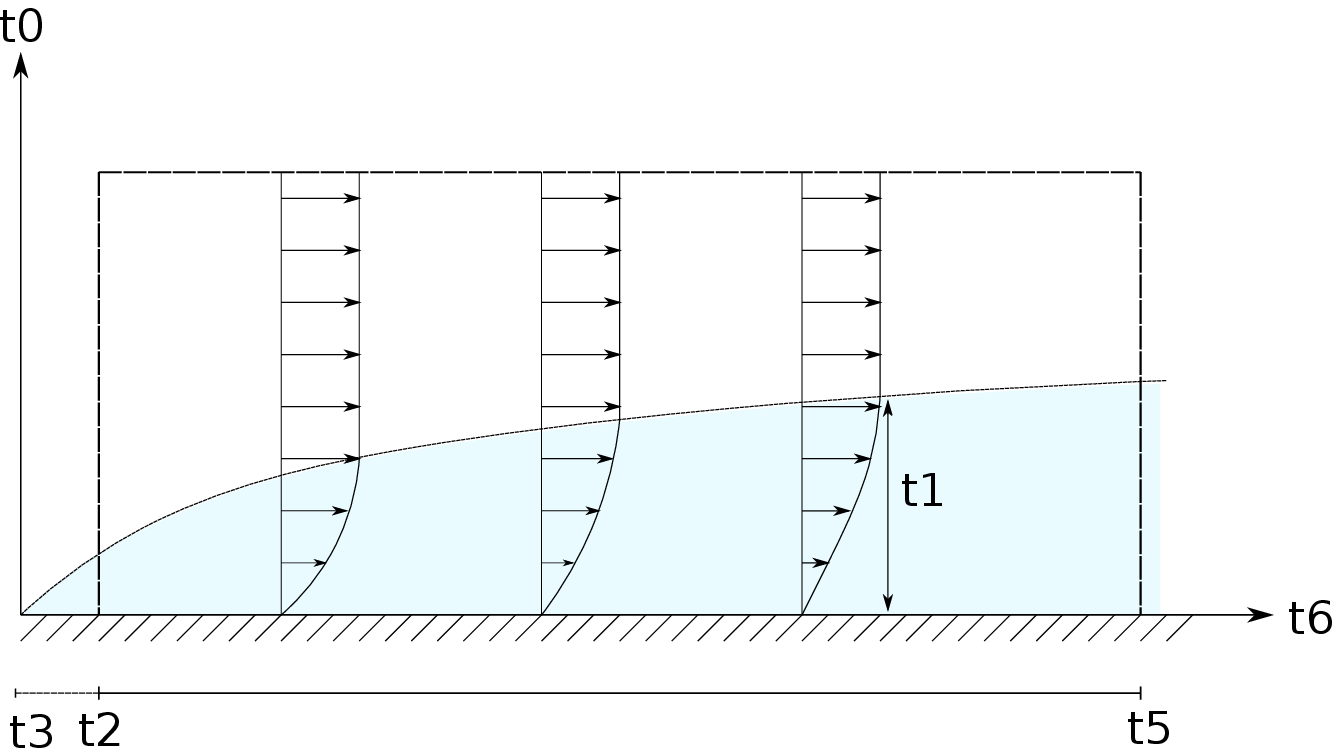}
\caption{}
         \label{fig8a}
     \end{subfigure}
     \hfill
     \begin{subfigure}[b]{0.5\textwidth}
        \centering
\includegraphics[scale=0.5]{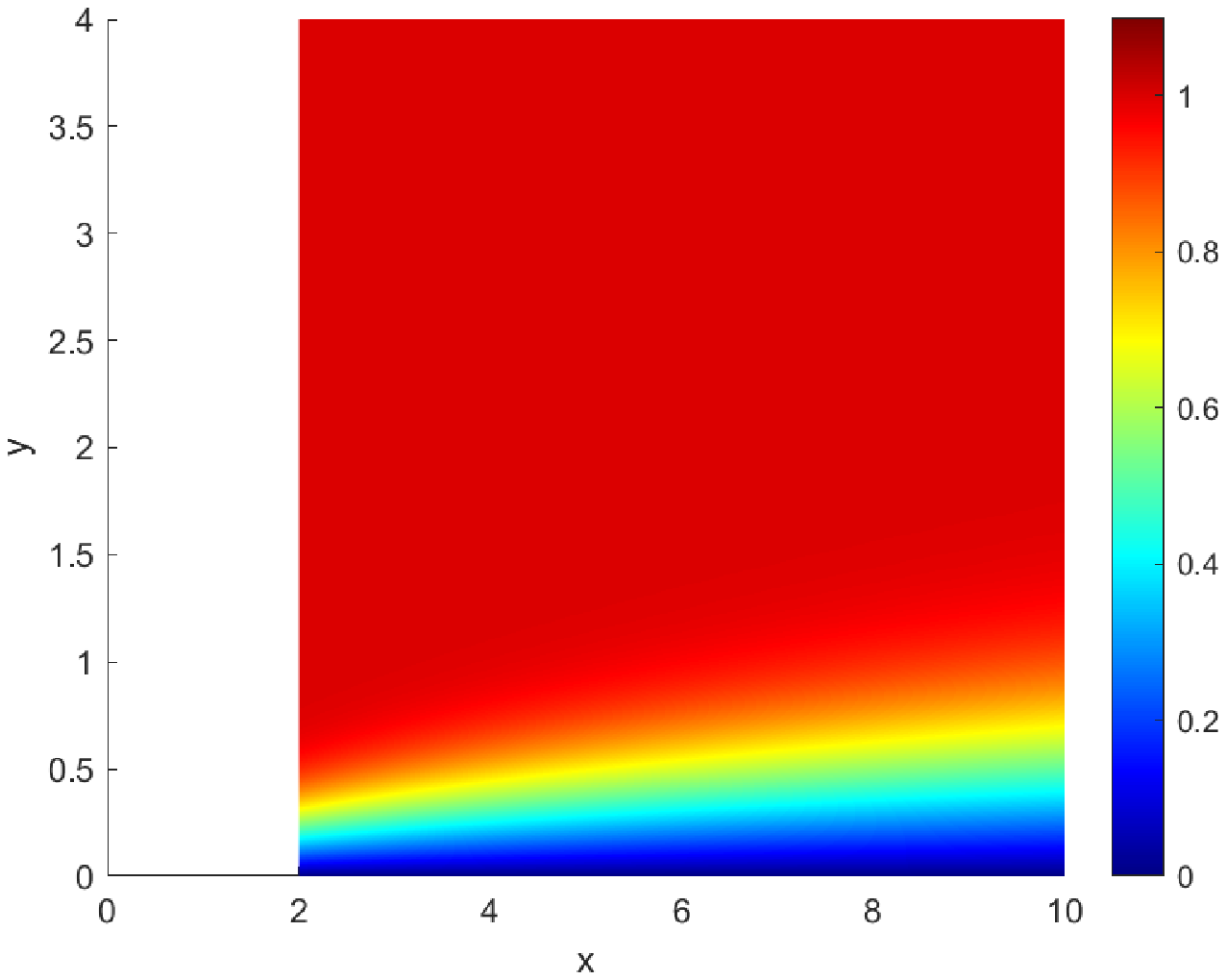}
         \caption{}
         \label{fig8b}
     \end{subfigure}   
        \caption{The illustration of the (a) truncated domain that exclude the leading edge and (b) the $u$-velocity distribution.}
         \label{fig8}
\end{figure}
\end{center}
\begin{center}
\begin{figure}[H]
\centering
 \begin{subfigure}[b]{0.3\textwidth}
 \centering
\includegraphics[scale=0.6]{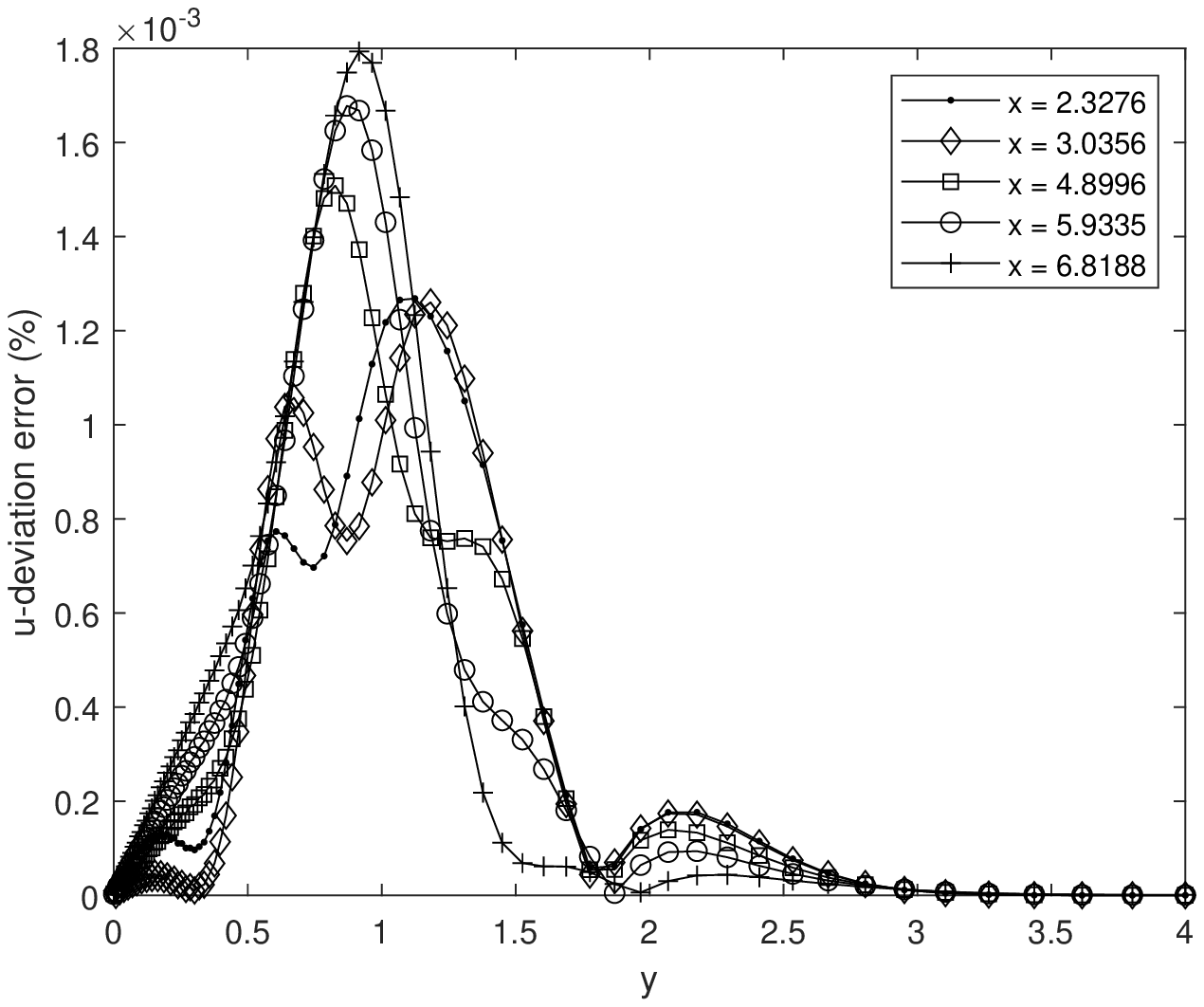}
\caption{}
     \label{fig9a}
     \end{subfigure}
     \hfill
     \begin{subfigure}[b]{0.5\textwidth}
     \centering
\includegraphics[scale=0.6]{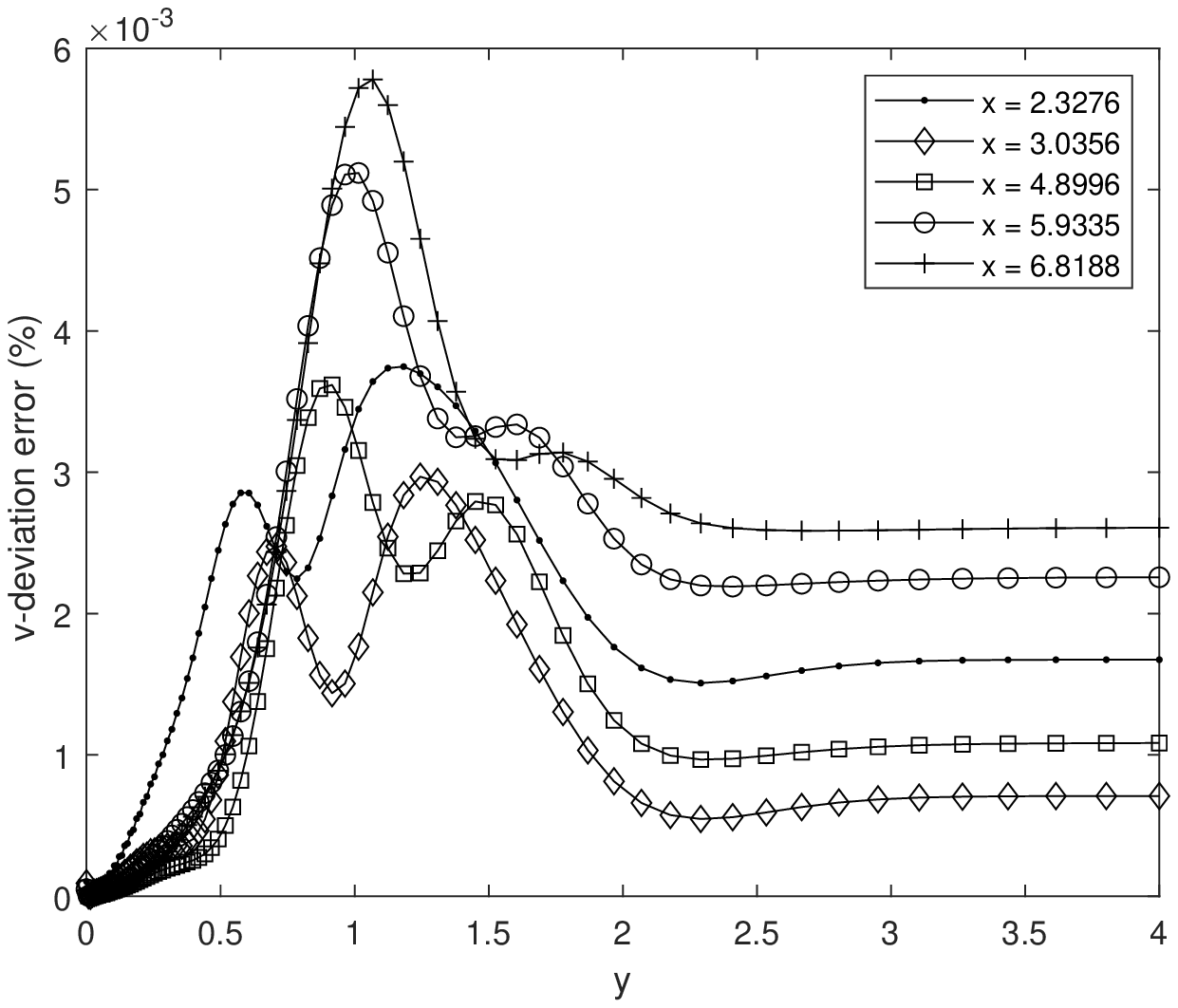}
     \caption{}
     \label{fig9b}
     \end{subfigure}   
     \caption{The deviation errors of the (a) $u$-velocity and (b) $v$-velocity at various $x$ position along the plate on the truncated domain.}
     \label{fig9}
\end{figure}
\end{center}
Similar to the full domain case,~we compare the deviations between the two solutions at various $x-$positions on the truncated domain \textcolor{blue}{(see Figure~$\ref{fig9}$)}.~A notable observation is the decrease in the magnitude of the errors.~Further,~though small,~the errors in $v$ do not drop to zero with increasing $y$ but asymptote to a value.~This was similarly the case for the full domain (Figure \ref{fig8b}).~The reason for this will be investigated as part of future work.~Next,~we compute the wall shear $\tau_w$ along the plate and compare to that computed from the Blasius solution (see \ref{Appendix}).~Note that at the tip of the plate,~the velocity gradients are infinite and moreover,~the Blasius solution does not exist at $x =  0$.~Therefore,~we will next consider the truncated domain.~The wall shear is computed as
\begin{align*}
\tau_w = \mu D_y \textbf{u}|_{y=0}
\end{align*}
and an accurate solution is achieved as shown in Figure~$\ref{fig10}$.
\begin{center}
\begin{figure}[H]
\centering
\includegraphics[scale=0.6]{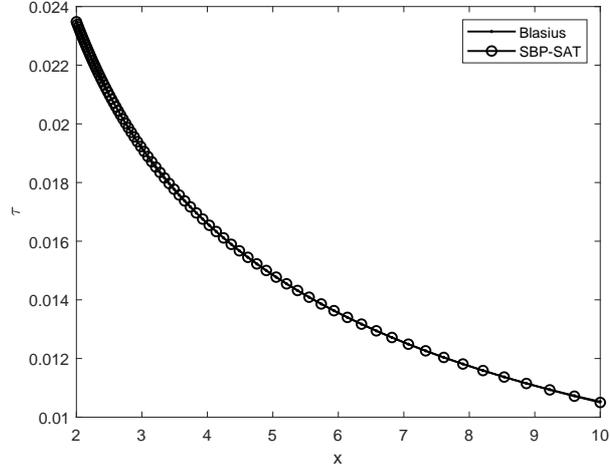} 
\caption{The wall shear along plate computed from the SBP-SAT approximation and the Blasius solution.}
\label{fig10}
\end{figure}
\end{center}
As a last test case,~we consider the stable INS approximation in \cite{Nordstrom2018,Nordstrom2020c} on a truncated domain and use the Blasius solution as inflow data.~At large Reynolds number,~we note that $\eqref{eq1}$ reduced to $\eqref{eq2}$.~We demonstrate this numerically by comparing the INS solution with the Blasius solution at different Reynolds numbers.~As depicted in Figure~$\ref{fig11}-\ref{fig12}$ and $\ref{fig13}-\ref{fig14}$,~the variation between the two approximation decreases as the Reynolds number increases.
\begin{figure}[H]
\captionsetup[subfigure]{justification=centering}
 \begin{subfigure}[b]{0.5\textwidth}
 \centering
\includegraphics[scale=0.5]{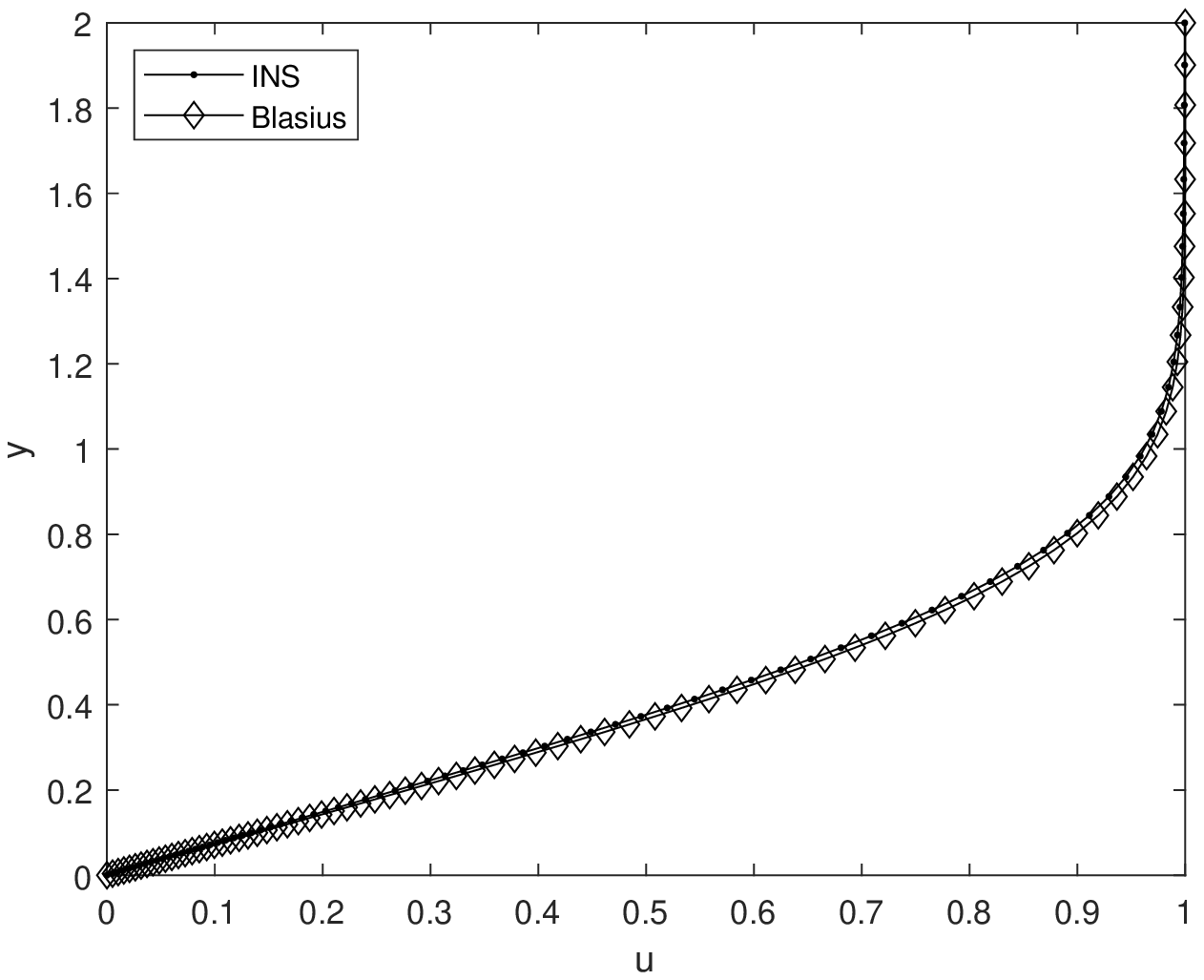}
\caption{}
     \label{fig11a}
     \end{subfigure}
     \hspace{1cm}
     \begin{subfigure}[b]{0.5\textwidth}
     \centering
\includegraphics[scale=0.5]{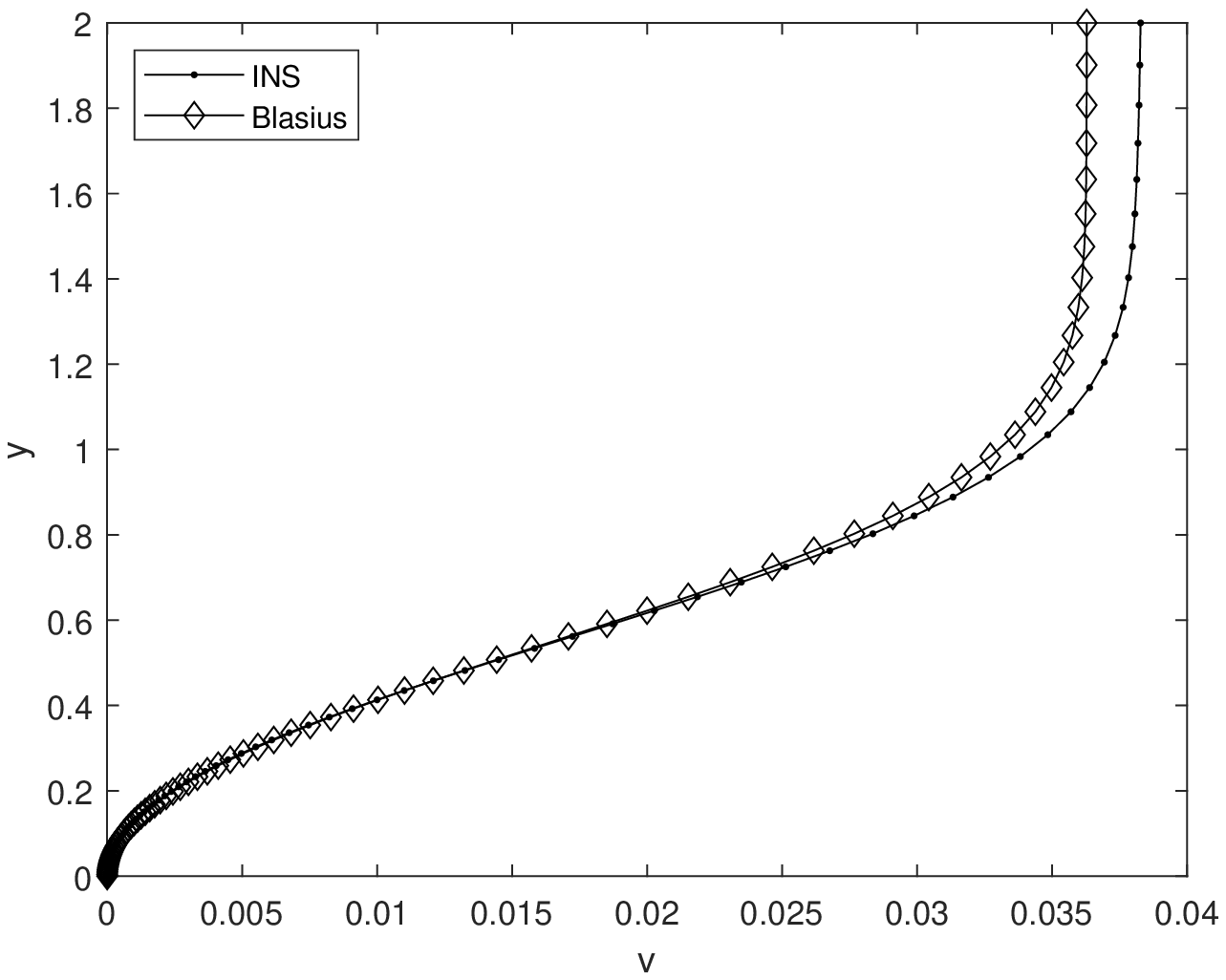}
     \caption{}
     \label{fig11b}
     \end{subfigure}   
     \caption{The SBP-SAT approximation for the INS equations compared with Blasius solution along the line $x = 5.63$ on the truncated domain at Re = 1 000 ($\mu = 1\text{e}-02$,~$\rho = 1$,~$l = 10$).~(a) $u$-velocity profile and (b) $v$-velocity profile.}
     \label{fig11}
\end{figure}
\begin{center}
\begin{figure}[H]
\captionsetup[subfigure]{justification=centering}
 \begin{subfigure}[b]{0.5\textwidth}
\includegraphics[scale=0.5]{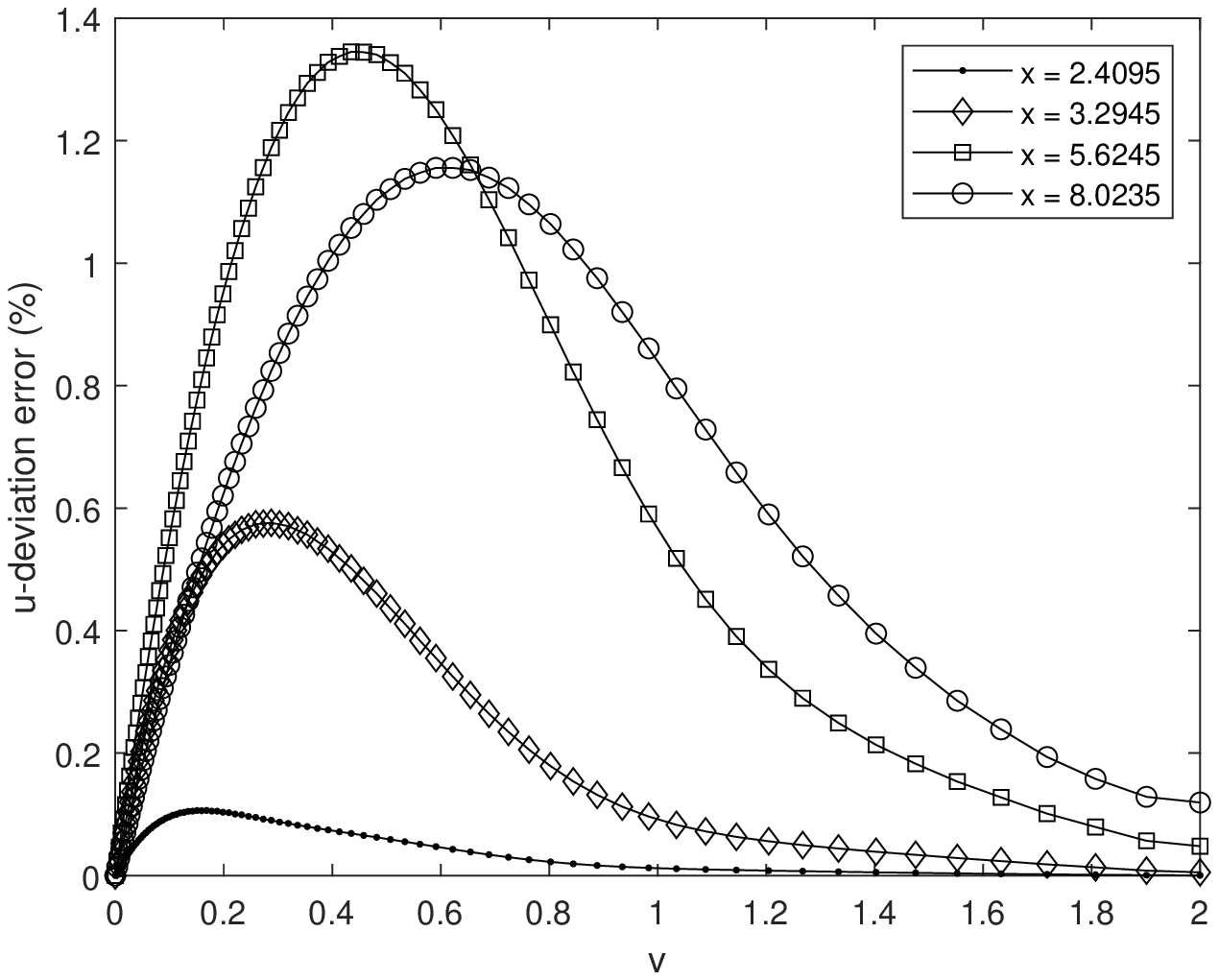}
\caption{}
     \label{fig12a}
     \end{subfigure}
     \hspace{1cm}
     \begin{subfigure}[b]{0.5\textwidth}
     \centering
\includegraphics[scale=0.5]{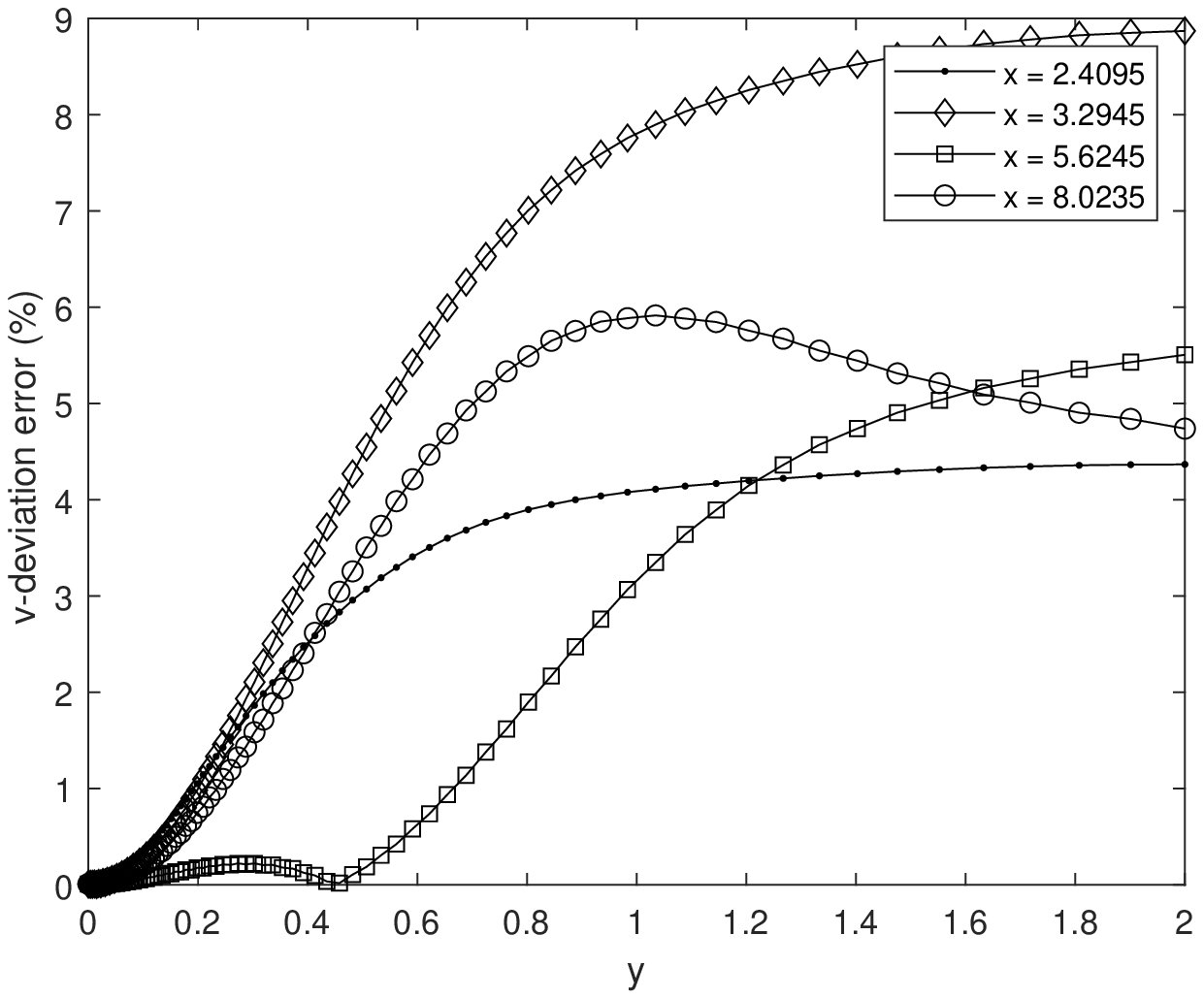}
     \caption{}
     \label{fig12b}
     \end{subfigure}   
     \caption{The deviation errors between the SBP-SAT approximation for the INS equations and the Blasius solution at Re = 1 000 ($\mu = 1\text{e}-02$,~$\rho = 1$,~$l = 10$).~We compare the (a) $u$-velocity and (b) $v$-velocity at various $x$ on the truncated domain.}
     \label{fig12}
\end{figure}
\end{center}
\begin{figure}[H]
\captionsetup[subfigure]{justification=centering}
 \begin{subfigure}[b]{0.5\textwidth}
 \centering
\includegraphics[scale=0.5]{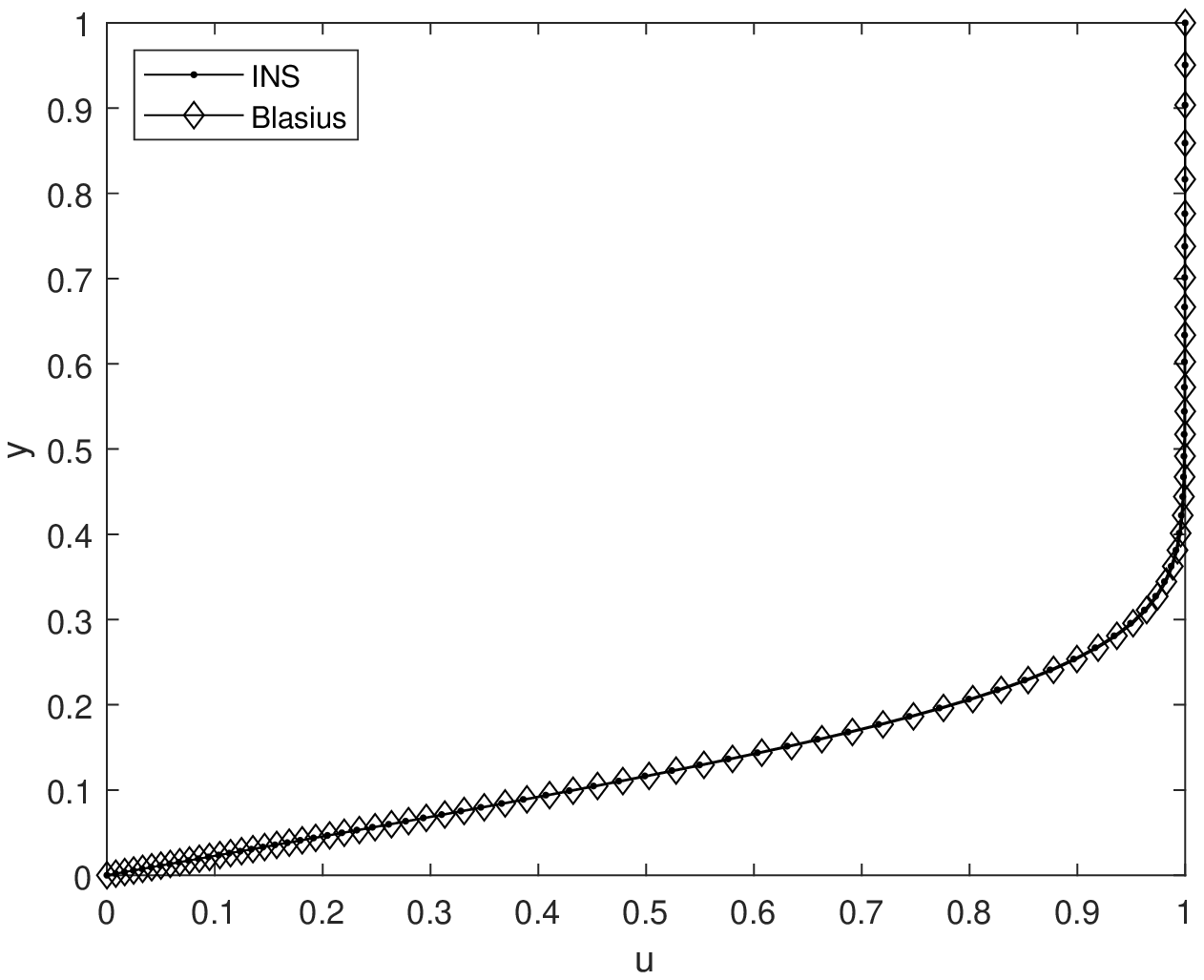}
\caption{}
     \label{fig13a}
     \end{subfigure}
     \hspace{1cm}
     \begin{subfigure}[b]{0.5\textwidth}
     \centering
\includegraphics[scale=0.5]{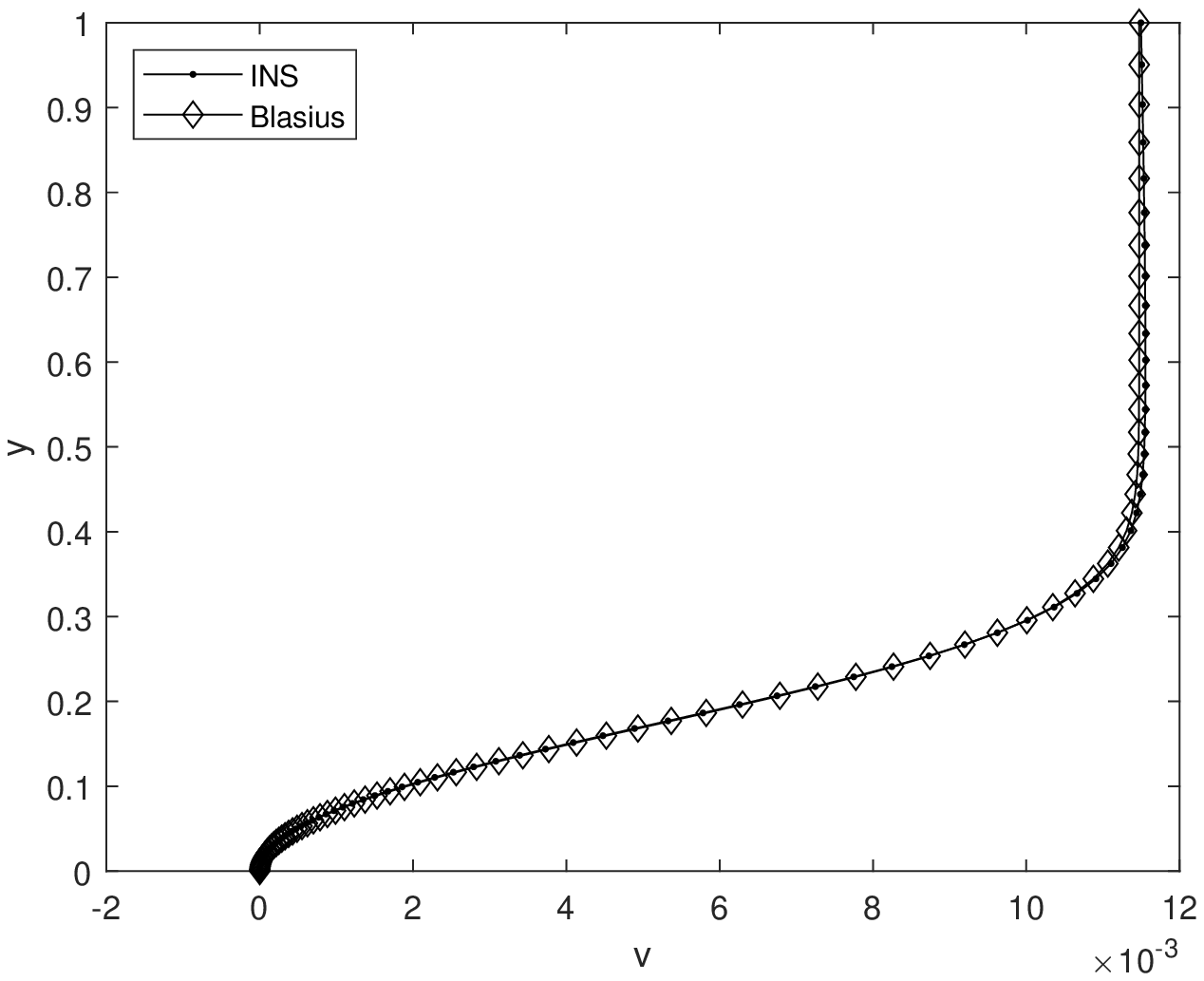}
     \caption{}
     \label{fig13b}
     \end{subfigure}   
     \caption{The SBP-SAT approximation for the INS equations compared with Blasius solution along the line $x = 5.63$ on the truncated domain at Re = 10 000 ($\mu = 1\text{e}-03$,~$\rho = 1$,~$l = 10$).~(a) $u$-velocity profile and (b) $v$-velocity profile.}
     \label{fig13}
\end{figure}
\begin{center}
\begin{figure}[H]
\captionsetup[subfigure]{justification=centering}
 \begin{subfigure}[b]{0.5\textwidth}
\includegraphics[scale=0.5]{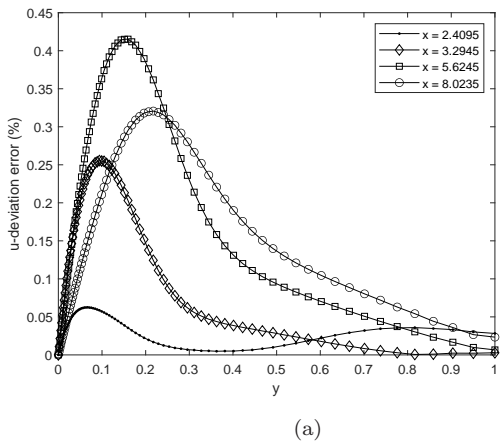}
\caption{}
     \label{fig14a}
     \end{subfigure}
     \hspace{1cm}
     \begin{subfigure}[b]{0.5\textwidth}
     \centering
\includegraphics[scale=0.5]{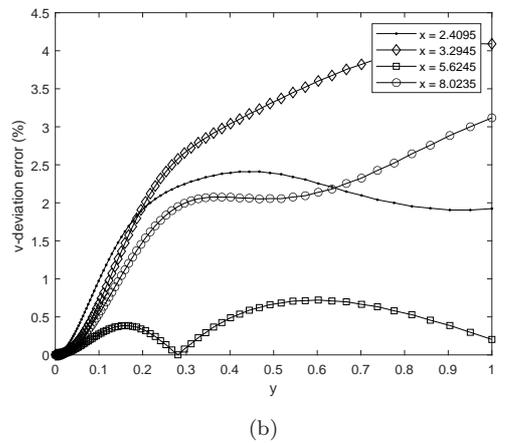}
     \caption{}
     \label{fig14}
     \end{subfigure}   
     \caption{The deviation errors between the SBP-SAT approximation for the INS equations and the Blasius solution at Re = 10 000 ($\mu = 1\text{e}-03$,~$\rho = 1$,~$l = 10$).~We compare the (a) $u$-velocity and (b) $v$-velocity at various $x$ on the truncated domain.}
     \label{fig14}
\end{figure}
\end{center}
\section{Summary and conclusion}\label{sec7}
This study was concerned with the development of a high-order accurate and stable finite difference approximation scheme for the incompressible laminar boundary layer equations.~We proposed a set of energy stable boundary conditions specific to the flat-plate boundary layer and obtained the energy estimate.~By mimicking the continuous analysis discretely,~we formulated SBP-SAT approximation scheme and proved stability.~Thus,~we obtained discrete estimates that resembled the continuous counterparts.
\par Stability of the resulting numerical scheme is proven by computing the eigenvalues of the resulting discrete spatial matrix.~A numerical application study proves expected high order spatial accuracy using an MMS.~This is followed by demonstrating an accurate correlation between our computed solution to the boundary layer equations and the celebrated Blasius similarity solution.~The study is concluded by demonstrating the Reynolds number effect on the solution of the incompressible Navier Stokes equations as compared to the Blasius solution.
\section*{CRediT authorship contribution statement}
\noindent Mojalefa Nchupang:~Conceptualization,~Methodology,~Software,~Writing-original draft,~Writing-review \& editing.\\ 
Arnaud Malan:~Conceptualization,~Methodology,~Writing-review \& editing.\\
Fredrik Laur\'en:~Methodology,~Software,~Writing-review \& editing.\\
Jan Nordstr\"om:~Conceptualization,~Methodology,~Writing-review \& editing.
\section*{Declaration of computing interests}
The authors declare that they have no known competing financial interests or personal relationships that could have appeared to influence the work reported in this paper.
\section*{Acknowledgments}
This work is based on research partly supported by the National Research Foundation of South Africa (Grant Numbers:~89916).~The opinions,~findings and conclusions or recommendations expressed are that of the authors alone,~and the NRF accepts no liability whatsoever in this regard.~The research leading to these results has also received funding from the European Union’s Horizon 2020 research and innovation programme under grant agreement No 815044,~the SLOshing Wing Dynamics (SLOWD) project.~The statements made herein do not necessarily have the consent or agreement of the SLOWD consortium and represent the opinion and findings of the author(s).
\par Jan Nordstr\"om was supported by Vetenskapsr\aa det,~Sweden [award no.~2018-05084 VR and 2021-0584].
\appendix
\section{Blasius similarity solution} \label{Appendix}
As with many PDEs governing fluid dynamics problems,~there is no known close-form solution of $\eqref{eq2}$.~Instead,~there is a well-known approximation solution called the Blasius similarity solution.~This solution method is based on the observation that,~the boundary layer model exhibits self-similar solution across the plate.~Therefore,~the PDE $\eqref{eq2}$ can be reduced to an ordinary differential equations (ODE) by introducing a similarity variable
\begin{align}
\eta = y \sqrt{\rho \frac{U_\infty}{\mu x}}. \label{eqA1}
\end{align}
Note that $\eta$ is not defined at the plate's leading edge (if it situated at the origin).~Further,~the Blasius solution method assumes zero pressure gradient in $\eqref{eq2}$.
\par Instead of computing the flow variables explicitly,~$u$,~$v$ are defined by stream functions 
\begin{align}
u &= \frac{\partial \psi}{\partial y}, \qquad  v = -\frac{\partial \psi}{\partial x}  \label{eqA2}
\end{align}
where $\psi = \sqrt{U_\infty \varepsilon x}f(\eta)$ and $f(\eta)$ is unknown.~In $\eqref{eqA2}$,~$u$,~$v$ further simplifies to
\begin{align}
u  = U_\infty f'(\eta),\quad v = \frac{1}{2} \sqrt{\varepsilon \frac{U_\infty }{x} } (\eta f'(\eta) -f(\eta))   \label{eqA3}.
\end{align}
Note that $\eqref{eqA3}$ satisfies the divergence relation since $\psi_{xy} = \psi_{yx}$.~By substituting $\eqref{eqA3}$ into $\eqref{eq2}$ with $p_x= 0$,~$\eqref{eq2}$ transforms to a nonlinear ODE 
\begin{align}
2 f'''(\eta) + f(\eta)f''(\eta) = 0. \label{eqA4}
\end{align}
Lastly,~to solve for $f$,~we need at least one boundary condition for $f$,~$f'$,~and $f''$.~It follows from the no-slip velocity condition that
\begin{align*}
u(x,0) &= U_\infty f'(0) = 0 \quad \Rightarrow f'(0) = 0, \\
v(x,0) &= \frac{1}{2} \sqrt{\varepsilon \frac{U_\infty }{x} } (0f'(0) -f(0)) = 0 \quad \Rightarrow f(0) = 0 .
\end{align*}
We also know that as $y \rightarrow \infty$,~ $u = U_\infty$ and using $\eqref{eqA2}$,~we get that $u(x,y\rightarrow \infty ) = U_\infty f'(\eta \rightarrow \infty) = U_\infty$.~Therefore,
\begin{align*}
f'(\eta \Rightarrow \infty) = 1.
\end{align*} 
Lastly,~by using the nonlinear shooting method \cite{Ha2001},~we determine $f''(0) = 0.332$.~Equation $\eqref{eqA4}$ can now be solved numerically using the 4th-order Runge-Kutta method.
\par By differentiating $\eqref{eqA2}$ with respect to $y$,~we compute the Newtonian shear stress on the plate in terms of the Blasius variables
\begin{align}
\tau_w = \mu U_\infty \sqrt{\rho \frac{U_\infty}{\mu x}} f''(0).  \label{A5}
\end{align}

\bibliography{Masterspapers}

\begin{thebibliography}{10}
\providecommand{\url}[1]{#1}
\csname url@samestyle\endcsname
\providecommand{\newblock}{\relax}
\providecommand{\bibinfo}[2]{#2}
\providecommand{\BIBentrySTDinterwordspacing}{\spaceskip=0pt\relax}
\providecommand{\BIBentryALTinterwordstretchfactor}{4}
\providecommand{\BIBentryALTinterwordspacing}{\spaceskip=\fontdimen2\font plus
\BIBentryALTinterwordstretchfactor\fontdimen3\font minus
  \fontdimen4\font\relax}
\providecommand{\BIBforeignlanguage}[2]{{%
\expandafter\ifx\csname l@#1\endcsname\relax
\typeout{** WARNING: IEEEtran.bst: No hyphenation pattern has been}%
\typeout{** loaded for the language `#1'. Using the pattern for}%
\typeout{** the default language instead.}%
\else
\language=\csname l@#1\endcsname
\fi
#2}}
\providecommand{\BIBdecl}{\relax}
\BIBdecl

\bibitem{Hume2022}
S.~Hume, J.~M.~I. Tshimanga, P.~Geoghegan, A.~G. Malan, W.~H. Ho, and M.~N.
  Ngoepe, ``{Effect of Pulsatility on the Transport of Thrombin in an Idealized
  Cerebral Aneurysm Geometry},'' \emph{Symmetry}, vol.~14, no.~1, pp. 1--18,
  2022.

\bibitem{YullPark2007}
J.~{Yull Park}, C.~{Young Park}, C.~{Mo Hwang}, K.~Sun, and B.~{Goo Min},
  ``{Pseudo-organ boundary conditions applied to a computational fluid dynamics
  model of the human aorta},'' \emph{Computers in Biology and Medicine},
  vol.~37, no.~8, pp. 1063--1072, 2007.

\bibitem{Cerrolaza1996}
M.~Cerrolaza, M.~Herrera, R.~Berrios, and W.~Annichiaricco, ``{A comparison of
  the hydrodynamical behaviour of three heart aortic prostheses by numerical
  methods},'' \emph{Journal of Medical Engineering and Technology}, vol.~20,
  no.~6, pp. 219--228, 1996.

\bibitem{Alam2015}
M.~F. Alam, D.~S. Thompson, and D.~K. Walters, ``{Hybrid Reynolds-averaged
  Navier-Stokes/large-eddy simulation models for flow around an iced wing},''
  \emph{Journal of Aircraft}, vol.~52, no.~1, pp. 244--256, 2015.

\bibitem{Kurzin2009}
V.~B. Kurzin and V.~A. Yudin, ``{Aerodynamic characteristics of a thin airfoil
  cascade in an ideal incompressible flow with separation from the leading
  edges},'' \emph{Fluid Dynamics}, vol.~44, no.~2, pp. 178--188, 2009.

\bibitem{Haddadpour2008}
H.~Haddadpour, M.~A. Kouchakzadeh, and F.~Shadmehri, ``{Aeroelastic instability
  of aircraft composite wings in an incompressible flow},'' \emph{Composite
  Structures}, vol.~83, no.~1, pp. 93--99, 2008.

\bibitem{Marshall1997b}
J.~Marshall, A.~Adcroft, C.~Hill, L.~Perelman, and C.~Heisey, ``{A
  finite-volume, incompressible navier stokes model for, studies of the ocean
  on parallel computers},'' \emph{Journal of Geophysical Research C: Oceans},
  vol. 102, no.~C3, pp. 5753--5766, 1997.

\bibitem{Teixeira2013}
\BIBentryALTinterwordspacing
P.~R. Teixeira, D.~P. Davyt, E.~Didier, and R.~Ramalhais, ``{Numerical
  simulation of an oscillating water column device using acode based on
  Navier-Stokes equations},'' \emph{Energy}, vol.~61, pp. 513--530, 2013.
  [Online]. Available: \url{http://dx.doi.org/10.1016/j.energy.2013.08.062}
\BIBentrySTDinterwordspacing

\bibitem{Malan2022}
L.~C. Malan, C.~Pilloton, A.~Colagrossi, and A.~G. Malan, ``{Numerical
  Calculation of Slosh Dissipation},'' \emph{Applied Sciences (Switzerland)},
  vol.~12, no.~23, pp. 1--31, 2022.

\bibitem{Mowat2014}
A.~G. Mowat, A.~G. Malan, L.~H. {Van Zyl}, and J.~P. Meyer, ``{Hybrid
  finite-volume reduced-order model method for nonlinear aeroelastic
  modeling},'' \emph{Journal of Aircraft}, vol.~51, no.~6, pp. 1805--1812,
  2014.

\bibitem{Changfoot2019}
D.~M. Changfoot, A.~G. Malan, and J.~Nordstr{\"{o}}m, ``{Hybrid
  computational-fluid-dynamics platform to investigate aircraft trailing
  vortices},'' \emph{Journal of Aircraft}, vol.~56, no.~1, pp. 344--355, 2019.

\bibitem{Patankar2018}
\BIBentryALTinterwordspacing
S.~V. Patankar, \emph{{Numerical Heat Transfer and Fluid Flow}}.\hskip 1em plus
  0.5em minus 0.4em\relax CRC Press, oct 2018. [Online]. Available:
  \url{https://www.taylorfrancis.com/books/9781482234213}
\BIBentrySTDinterwordspacing

\bibitem{Chorin1967}
A.~J. Chorin, ``{A numerical method for solving incompressible viscous flow
  problems},'' \emph{Journal of Computational Physics}, vol.~2, no.~1, pp.
  12--26, 1967.

\bibitem{Malan2002}
A.~G. Malan, R.~W. Lewis, and P.~Nithiarasu, ``{An improved unsteady,
  unstructured, artificial compressibility, finite volume scheme for viscous
  incompressible flows: Part I. Theory and implementation},''
  \emph{International Journal for Numerical Methods in Engineering}, vol.~54,
  no.~5, pp. 695--714, 2002.

\bibitem{Malan2013}
\BIBentryALTinterwordspacing
A.~G. Malan and O.~F. Oxtoby, ``{An accelerated, fully-coupled, parallel 3D
  hybrid finite-volume fluid-structure interaction scheme},'' \emph{Computer
  Methods in Applied Mechanics and Engineering}, vol. 253, pp. 426--438, 2013.
  [Online]. Available: \url{http://dx.doi.org/10.1016/j.cma.2012.09.004}
\BIBentrySTDinterwordspacing

\bibitem{Merrick2018}
\BIBentryALTinterwordspacing
D.~G. Merrick, A.~G. Malan, and J.~A. van Rooyen, ``{A novel finite volume
  discretization method for advection–diffusion systems on stretched
  meshes},'' \emph{Journal of Computational Physics}, vol. 362, pp. 220--242,
  2018. [Online]. Available: \url{https://doi.org/10.1016/j.jcp.2018.02.025}
\BIBentrySTDinterwordspacing

\bibitem{KREISS1974}
\BIBentryALTinterwordspacing
H.-O. KREISS and G.~SCHERER, ``{Finite Element and Finite Difference Methods
  for Hyperbolic Partial Differential Equations},'' in \emph{Mathematical
  Aspects of Finite Elements in Partial Differential Equations}.\hskip 1em plus
  0.5em minus 0.4em\relax Elsevier, 1974, pp. 195--212. [Online]. Available:
  \url{https://linkinghub.elsevier.com/retrieve/pii/B9780122083501500121}
\BIBentrySTDinterwordspacing

\bibitem{Gustafsson2008}
\BIBentryALTinterwordspacing
B.~Gustafsson, \emph{{High Order Difference Methods for Time Dependent PDE}},
  ser. Springer Series in Computational Mathematics.\hskip 1em plus 0.5em minus
  0.4em\relax Berlin, Heidelberg: Springer Berlin Heidelberg, 2008, vol.~38.
  [Online]. Available: \url{http://link.springer.com/10.1007/978-3-540-74993-6}
\BIBentrySTDinterwordspacing

\bibitem{Gustafsson2013}
\BIBentryALTinterwordspacing
B.~Gustafsson, H.-O. Kreiss, and J.~Oliger, \emph{{Time-Dependent Problems and
  Difference Methods}}.\hskip 1em plus 0.5em minus 0.4em\relax Hoboken, NJ,
  USA: John Wiley {\&} Sons, Inc., sep 2013. [Online]. Available:
  \url{http://doi.wiley.com/10.1002/9781118548448}
\BIBentrySTDinterwordspacing

\bibitem{Carpenter1994}
\BIBentryALTinterwordspacing
M.~H. Carpenter, D.~Gottlieb, and S.~Abarbanel, ``{Time-Stable Boundary
  Conditions for Finite-Difference Schemes Solving Hyperbolic Systems:
  Methodology and Application to High-Order Compact Schemes},'' \emph{Journal
  of Computational Physics}, vol. 111, no.~2, pp. 220--236, apr 1994. [Online].
  Available:
  \url{https://linkinghub.elsevier.com/retrieve/pii/S0021999184710576}
\BIBentrySTDinterwordspacing

\bibitem{Strikwerda1977}
\BIBentryALTinterwordspacing
J.~C. Strikwerda, ``{Initial boundary value problems for incompletely parabolic
  systems},'' \emph{Communications on Pure and Applied Mathematics}, vol.~30,
  no.~6, pp. 797--822, nov 1977. [Online]. Available:
  \url{https://onlinelibrary.wiley.com/doi/10.1002/cpa.3160300606}
\BIBentrySTDinterwordspacing

\bibitem{Nordstrom2020b}
\BIBentryALTinterwordspacing
J.~Nordstr{\"{o}}m and T.~M. Hagstrom, ``{The Number of Boundary Conditions for
  Initial Boundary Value Problems},'' \emph{SIAM Journal on Numerical
  Analysis}, vol.~58, no.~5, pp. 2818--2828, jan 2020. [Online]. Available:
  \url{https://epubs.siam.org/doi/10.1137/20M1322571}
\BIBentrySTDinterwordspacing

\bibitem{Kreiss1970a}
\BIBentryALTinterwordspacing
H.-O. Kreiss, ``{Initial boundary value problems for hyperbolic systems},''
  \emph{Communications on Pure and Applied Mathematics}, vol.~23, no.~3, pp.
  277--298, may 1970. [Online]. Available:
  \url{https://onlinelibrary.wiley.com/doi/10.1002/cpa.3160230304}
\BIBentrySTDinterwordspacing

\bibitem{Nordstrom2017a}
J.~Nordstr{\"{o}}m, ``{A Roadmap to Well Posed and Stable Problems in
  Computational Physics},'' \emph{Journal of Scientific Computing}, vol.~71,
  no.~1, pp. 365--385, 2017.

\bibitem{Lauren2022}
\BIBentryALTinterwordspacing
F.~Laur{\'{e}}n and J.~Nordstr{\"{o}}m, ``{Energy stable wall modeling for the
  Navier-Stokes equations},'' \emph{Journal of Computational Physics}, vol.
  457, p. 111046, 2022. [Online]. Available:
  \url{https://doi.org/10.1016/j.jcp.2022.111046}
\BIBentrySTDinterwordspacing

\bibitem{Nordstrom2018}
\BIBentryALTinterwordspacing
J.~Nordstr{\"{o}}m and C.~{La Cognata}, ``{Energy stable boundary conditions
  for the nonlinear incompressible Navier–Stokes equations},''
  \emph{Mathematics of Computation}, vol.~88, no. 316, pp. 665--690, aug 2018.
  [Online]. Available:
  \url{http://www.ams.org/mcom/2019-88-316/S0025-5718-2018-03375-0/}
\BIBentrySTDinterwordspacing

\bibitem{Nordstrom2022}
\BIBentryALTinterwordspacing
J.~Nordstr{\"{o}}m and A.~R. Winters, ``{A linear and nonlinear analysis of the
  shallow water equations and its impact on boundary conditions},''
  \emph{Journal of Computational Physics}, vol. 463, p. 111254, 2022. [Online].
  Available: \url{https://doi.org/10.1016/j.jcp.2022.111254}
\BIBentrySTDinterwordspacing

\bibitem{Nordstrom2015a}
J.~Nordstr{\"{o}}m and S.~Ghader, ``{A new well-posed vorticity divergence
  formulation of the shallow water equations},'' \emph{Ocean Modelling},
  vol.~93, pp. 1--6, 2015.

\bibitem{Nordstrom2022a}
\BIBentryALTinterwordspacing
J.~Nordstr{\"{o}}m and F.~Laur{\'{e}}n, ``{A stable and conservative nonlinear
  interface coupling for the incompressible Euler equations},'' \emph{Applied
  Mathematics Letters}, vol. 132, p. 108171, 2022. [Online]. Available:
  \url{https://doi.org/10.1016/j.aml.2022.108171}
\BIBentrySTDinterwordspacing

\bibitem{Nordstrom2022b}
\BIBentryALTinterwordspacing
J.~Nordstr{\"{o}}m, ``{A skew-symmetric energy and entropy stable formulation
  of the compressible Euler equations},'' \emph{Journal of Computational
  Physics}, vol. 470, p. 111573, 2022. [Online]. Available:
  \url{https://doi.org/10.1016/j.jcp.2022.111573}
\BIBentrySTDinterwordspacing

\bibitem{Benzi2005}
M.~Benzi, G.~H. Golubt, and J.~Liesen, ``{Numerical solution of saddle point
  problems},'' \emph{Acta Numerica}, vol.~14, pp. 1--137, 2005.

\bibitem{Nordstrom2020c}
\BIBentryALTinterwordspacing
J.~Nordstr{\"{o}}m and F.~Laur{\'{e}}n, ``{The spatial operator in the
  incompressible Navier–Stokes, Oseen and Stokes equations},'' \emph{Computer
  Methods in Applied Mechanics and Engineering}, vol. 363, p. 112857, 2020.
  [Online]. Available: \url{https://doi.org/10.1016/j.cma.2020.112857}
\BIBentrySTDinterwordspacing

\bibitem{Lauren2021}
\BIBentryALTinterwordspacing
F.~Laur{\'{e}}n and J.~Nordstr{\"{o}}m, ``{Spectral properties of the
  incompressible Navier-Stokes equations},'' \emph{Journal of Computational
  Physics}, vol. 429, p. 110019, 2021. [Online]. Available:
  \url{https://doi.org/10.1016/j.jcp.2020.110019}
\BIBentrySTDinterwordspacing

\bibitem{Frank2006b}
W.~M. Frank, \emph{{Viscous Fluid Flow}}.\hskip 1em plus 0.5em minus
  0.4em\relax McGraw-Hill, 2006.

\bibitem{Rahman2001}
M.~M. Rahman and T.~Siikonen, ``{An artificial compressibility method for
  incompressible flows},'' \emph{Numerical Heat Transfer, Part B:
  Fundamentals}, vol.~40, no.~5, pp. 391--409, 2001.

\bibitem{Vreman2014a}
\BIBentryALTinterwordspacing
A.~W. Vreman, ``{The projection method for the incompressible Navier-Stokes
  equations: The pressure near a no-slip wall},'' \emph{Journal of
  Computational Physics}, vol. 263, pp. 353--374, 2014. [Online]. Available:
  \url{http://dx.doi.org/10.1016/j.jcp.2014.01.035}
\BIBentrySTDinterwordspacing

\bibitem{Matsui2021}
\BIBentryALTinterwordspacing
K.~Matsui, ``{A projection method for Navier-Stokes equations with a boundary
  condition including the total pressure},'' pp. 1--30, 2021. [Online].
  Available: \url{http://arxiv.org/abs/2105.13014}
\BIBentrySTDinterwordspacing

\bibitem{Oxtoby2012a}
\BIBentryALTinterwordspacing
O.~F. Oxtoby and A.~G. Malan, ``{A matrix-free, implicit, incompressible
  fractional-step algorithm for fluid-structure interaction applications},''
  \emph{Journal of Computational Physics}, vol. 231, no.~16, pp. 5389--5405,
  2012. [Online]. Available: \url{http://dx.doi.org/10.1016/j.jcp.2012.04.037}
\BIBentrySTDinterwordspacing

\bibitem{OReilly2017a}
\BIBentryALTinterwordspacing
O.~O'Reilly, T.~Lundquist, E.~M. Dunham, and J.~Nordstr{\"{o}}m, ``{Energy
  stable and high-order-accurate finite difference methods on staggered
  grids},'' \emph{Journal of Computational Physics}, vol. 346, pp. 572--589,
  2017. [Online]. Available: \url{http://dx.doi.org/10.1016/j.jcp.2017.06.030}
\BIBentrySTDinterwordspacing

\bibitem{Kress2003}
W.~Kress and J.~Nilsson, ``{Boundary conditions and estimates for the
  linearized Navier-Stokes equations on staggered grids},'' \emph{Computers and
  Fluids}, vol.~32, no.~8, pp. 1093--1112, 2003.

\bibitem{Gustafsson2000}
B.~Gustafsson and J.~Nilsson, ``{Boundary conditions and estimates for the
  steady Stokes equations on staggered grids},'' \emph{Journal of Scientific
  Computing}, vol.~15, no.~1, pp. 29--59, 2000.

\bibitem{Manzanero2019}
\BIBentryALTinterwordspacing
J.~Manzanero, G.~Rubio, D.~A. Kopriva, E.~Ferrer, and E.~Valero, ``{A
  free-energy stable nodal discontinuous Galerkin approximation with
  summation-by-parts property for the Cahn-Hilliard equation},'' \emph{Journal
  of Computational Physics}, vol.~1, p. 109072, 2019. [Online]. Available:
  \url{http://arxiv.org/abs/1902.08089}
\BIBentrySTDinterwordspacing

\bibitem{Chan2018}
\BIBentryALTinterwordspacing
J.~Chan, ``{On discretely entropy conservative and entropy stable discontinuous
  Galerkin methods},'' \emph{Journal of Computational Physics}, vol. 362, pp.
  346--374, 2018. [Online]. Available:
  \url{https://doi.org/10.1016/j.jcp.2018.02.033}
\BIBentrySTDinterwordspacing

\bibitem{Yamaleev2017}
\BIBentryALTinterwordspacing
N.~K. Yamaleev and M.~H. Carpenter, ``{A family of fourth-order entropy stable
  nonoscillatory spectral collocation schemes for the 1-D Navier–Stokes
  equations},'' \emph{Journal of Computational Physics}, vol. 331, pp. 90--107,
  2017. [Online]. Available: \url{http://dx.doi.org/10.1016/j.jcp.2016.11.039}
\BIBentrySTDinterwordspacing

\bibitem{Abgrall2020}
\BIBentryALTinterwordspacing
R.~Abgrall, J.~Nordstr{\"{o}}m, P.~{\"{O}}ffner, and S.~Tokareva, ``{Analysis
  of the SBP-SAT Stabilization for Finite Element Methods Part I: Linear
  Problems},'' \emph{Journal of Scientific Computing}, vol.~85, no.~2, 2020.
  [Online]. Available: \url{https://doi.org/10.1007/s10915-020-01349-z}
\BIBentrySTDinterwordspacing

\bibitem{Ham2006}
F.~Ham, K.~Mattsson, and G.~Iaccarino, ``{Accurate and stable finite volume
  operators for unstructured flow solvers},'' \emph{Center for Turbulence
  Research Annual Research Briefs}, pp. 243--261, 2006.

\bibitem{Nordstrom2003}
J.~Nordstr{\"{o}}m, K.~Forsberg, C.~Adamsson, and P.~Eliasson, ``{Finite volume
  methods, unstructured meshes and strict stability for hyperbolic problems},''
  \emph{Applied Numerical Mathematics}, vol.~45, no.~4, pp. 453--473, 2003.

\bibitem{Nordstrom2006a}
J.~Nordstr{\"{o}}m, ``{Conservative finite difference formulations, variable
  coefficients, energy estimates and artificial dissipation},'' \emph{Journal
  of Scientific Computing}, vol.~29, no.~3, pp. 375--404, 2006.

\bibitem{Nordstrom2023a}
\BIBentryALTinterwordspacing
------, ``{Nonlinear Boundary Conditions for Energy and Entropy Stable Initial
  Boundary Value Problems in Computational Fluid Dynamics},'' pp. 1--13, 2023.
  [Online]. Available: \url{http://arxiv.org/abs/2301.04568}
\BIBentrySTDinterwordspacing

\bibitem{Box2005}
J.~Sudirham, J.~{Van Der Vegt}, and R.~{Van Damme}, ``{A study on Discontinuous
  Galerkin finite elements methods for elliptic problems},'' \emph{Memorandum
  No. 1690, University of Twente, Faculty of EEMCS}, no. 1690, 2003.

\bibitem{Arnold2001}
D.~N. Arnold, F.~Brezzi, B.~Cockburn, and L.~{Donatella Marini}, ``{Unified
  analysis of discontinuous Galerkin methods for elliptic problems},''
  \emph{SIAM Journal on Numerical Analysis}, vol.~39, no.~5, pp. 1749--1779,
  2001.

\bibitem{Mattsson2004}
K.~Mattsson and J.~Nordstr{\"{o}}m, ``{Summation by parts operators for finite
  difference approximations of second derivatives},'' \emph{Journal of
  Computational Physics}, vol. 199, no.~2, pp. 503--540, 2004.

\bibitem{Svard2006}
M.~Sv{\"{a}}rd and J.~Nordstr{\"{o}}m, ``{On the order of accuracy for
  difference approximations of initial-boundary value problems},''
  \emph{Journal of Computational Physics}, vol. 218, no.~1, pp. 333--352, 2006.

\bibitem{Svard2014}
\BIBentryALTinterwordspacing
------, ``{Review of summation-by-parts schemes for initial-boundary-value
  problems},'' \emph{Journal of Computational Physics}, vol. 268, pp. 17--38,
  2014. [Online]. Available: \url{http://dx.doi.org/10.1016/j.jcp.2014.02.031}
\BIBentrySTDinterwordspacing

\bibitem{Roache}
\BIBentryALTinterwordspacing
P.~J. Roache, ``{The Method of Manufactured Solutions for Code Verification},''
  in \emph{Computer Simulation Validation. Simulation Foundations, Methods and
  Applications}, 2019, pp. 295--318. [Online]. Available:
  \url{http://link.springer.com/10.1007/978-3-319-70766-2{\_}12}
\BIBentrySTDinterwordspacing

\bibitem{Petersson2015}
\BIBentryALTinterwordspacing
N.~A. Petersson and B.~Sj{\"{o}}green, ``{Wave propagation in anisotropic
  elastic materials and curvilinear coordinates using a summation-by-parts
  finite difference method},'' \emph{Journal of Computational Physics}, vol.
  299, pp. 820--841, 2015. [Online]. Available:
  \url{http://dx.doi.org/10.1016/j.jcp.2015.07.023}
\BIBentrySTDinterwordspacing

\bibitem{Gong2007}
J.~Gong and J.~Nordstr{\"{o}}m, ``{A stable and efficient hybrid scheme for
  viscous problems in complex geometries},'' \emph{Journal of Computational
  Physics}, vol. 226, no.~2, pp. 1291--1309, 2007.

\bibitem{Nordstrom2001}
J.~Nordstr{\"{o}}m and M.~H. Carpenter, ``{High-order finite difference
  methods, multidimensional linear problems, and curvilinear coordinates},''
  \emph{Journal of Computational Physics}, vol. 173, no.~1, pp. 149--174, 2001.

\bibitem{Alund2019}
O.~{\AA}lund and J.~Nordstr{\"{o}}m, ``{Encapsulated high order difference
  operators on curvilinear non-conforming grids},'' \emph{Journal of
  Computational Physics}, vol. 385, pp. 209--224, 2019.

\bibitem{Capatina2021}
D.~Capatina, D.~Capatina, D.~Graebling, and D.~Trujillo, ``{Velocity overshoot
  for incompressible flows past a semi-infinite flat plate},'' no. October,
  2021.

\bibitem{Ha2001}
S.~N. Ha, ``{A nonlinear shooting method for two-point boundary value
  problems},'' \emph{Computers and Mathematics with Applications}, vol.~42, no.
  10-11, pp. 1411--1420, 2001.

\end{thebibliography}
\bibliographystyle{IEEEtran}
\end{document}